\documentclass[review,onefignum,onetabnum]{siamart190516}

\usepackage{amssymb}
\usepackage{verbatim}
\usepackage{lipsum}
\usepackage{amsfonts}
\usepackage{graphicx}
\usepackage{subfigure}
\usepackage{epstopdf}
\usepackage{algorithmic}
\ifpdf
\DeclareGraphicsExtensions{.eps,.pdf,.png,.jpg}
\else
\DeclareGraphicsExtensions{.eps}
\fi
\renewcommand\theequation{\thesection.\arabic{equation}}
\theoremstyle{plain}

\newtheorem{remark}{\textbf{Remark}}[section]

\newcommand{\eps}{\epsilon}

\newcommand{\bm}{\boldsymbol}

\newcommand{\Grad}[1]{\nabla #1}

\newcommand{\be}{\begin{equation}}
\newcommand{\ee}{\end{equation}}

\newcommand{\bse}{\begin{subequations}}
\newcommand{\ese}{\end{subequations}}
\def\benl{\begin{eqnarray*}}
\def\eenl{\end{eqnarray*}}

\def\be{\bm{e}}

\def\bx{\bm{x}}

\def\bz{\bm{z}}

\def\bmu1{\bm{\mu_1}}

\newcommand{\ole}{\overline{e}}
\newcommand{\te}{\tilde{e}}
\newcommand{\he}{\hat{e}}

\newcommand{\ben}{\begin{eqnarray}}
\newcommand{\een}{\end{eqnarray}}
\newcommand{\beq}{\begin{equation}}
\newcommand{\eeq}{\end{equation}}
\newcommand{\bea}{\begin{array}}
\newcommand{\eea}{\end{array}}
\newcommand{\bef}{\begin{figure}[H]}
\newcommand{\eef}{\end{figure}}

\newsiamremark{hypothesis}{Hypothesis}
\crefname{hypothesis}{Hypothesis}{Hypotheses}
\newsiamthm{claim}{Claim}

\title{ A new class  of  energy dissipative, mass conserving and positivity/bound-preserving schemes    for Keller-Segel equations}

\author{Zexiong Fang \thanks{Department of Mathematics, Tongji University, Shanghai 200092, China (2311746$@$tongji.edu.cn). Key Laboratory of Intelligent Computing and Applications (Tongji University), Ministry of Education, China. }
\and Qing Cheng  \thanks{Department of Mathematics, Tongji University, Shanghai 200092, China (qingcheng$@$tongji.edu.cn). Key Laboratory of Intelligent Computing and Applications (Tongji University), Ministry of Education, China. }
}
  
\usepackage{amsopn}

\begin{document}
\bibliographystyle{plain}
\graphicspath{ {Figures_porous/} }
\maketitle

\begin{abstract}
In this paper, we   improve  the original Lagrange multiplier approach \cite{ChSh22,ChSh_II22} and  introduce  a new energy correction approach to construct a class of robust,  positivity/bound-preserving, mass conserving and energy dissipative  schemes for  Keller-Segel equations which only need to solve several  linear Poisson like equations. To be more specific,  we use  a predictor-corrector approach to construct a class of positivity/bound-preserving and mass conserving schemes which can be implemented with negligible cost. Then a energy correction step is introduced to construct schemes which are also energy dissipative, in addition to positivity/bound-preserving and mass conserving.   This new approach is not restricted to any particular spatial discretization and can be  combined with various time discretization to achieve high-order accuracy in time. 
We show stability  results for mass-conservative, positivity/bound-preserving and energy dissipative schemes   for two different Keller-Segel systems. A error analysis is presented for a second-order, bound-preserving, mass-conserving and energy dissipative scheme for the second-type of Keller-Segel equations. Ample  numerical experiments are shown  to validate the stability and accuracy of our approach. 
\end{abstract}




\section{Introduction}
The paper is focused on designing new positivity/bound-preserving schemes for two types  of  Keller-Segel equations \cite{bellomo2015toward,horstmann20031970,liu2018positivity,shen2020unconditionally,hillen2001global}  describing the phenomenon of chemotaxis, in which cells approach the chemically favorable environments according to the chemical substance generated by cells.  The first type  of 2D Keller-Segel equations \cite{zhou2017finite,filbet2006finite,chertock2008second} in domain $\Omega$  can be written into the following form:
\begin{eqnarray}
&& \partial_t\rho = \gamma\Delta \rho -\chi\Grad\cdot(\rho\Grad c),\label{de:keller:1}\\
&& \eps\partial_t c= \mu\Delta c +\rho,\label{de:keller:2}\\
&&\rho(\bx,0)=f(\bx),\; c(\bx,0)=g(\bx), \quad \bx \in \Omega, \label{de:keller:3}
\end{eqnarray}
where $f(\bx)$ and $g(\bx)$ are given initial conditions. The second type of Keller-Segel equations \cite{dolak2005keller} with chemotactic sensitivity function  $\eta(\rho)=\frac{\rho(M-\rho)}{M}$   are formulated as
\begin{eqnarray}
&& \partial_t\rho = \gamma\Delta \rho -\chi\Grad\cdot(\eta(\rho)\Grad c),\label{keller:1}\\
&& \eps\partial_t c= \mu\Delta c +\rho,\label{keller:2}\\
&& \rho(\bx,0)=f(\bx),\; c(\bx,0)=g(\bx), \label{keller:3}
\end{eqnarray} 
 where  $M$ is the  maximal density.  Parameters  $\gamma$, $\chi$, $\mu$ and $\eps$ are all  positive constants. Specifically, $\chi$ represents chemotactic sensitivity and $\eps$ describes how fast the chemoattractant concentration reacts to the organism concentration. For the second type of Keller-Segel equations  which are first investigated by Hillen and Painter in \cite{hillen2001global}, the chemotactic response of the cells is shut off when a maximal density is reached. 
 
 Suitable boundary conditions are  periodic  or no-flux boundary conditions on $\rho(\bx,t)$ and the Neumann boundary condition on $c(\bx,t)$ for two types of Keller-Segel equations. Here $\rho(\bx,t)$ denotes the density distribution of cells and $c(\bx,t)$ is the chemical concentration. Mathematically, the Keller-Segel equations \eqref{de:keller:1}-\eqref{de:keller:3} and \eqref{keller:1}-\eqref{keller:3}  describe the competition beteen diffusion and the nonlocal aggregation in biology or social science systems. When $\eps>0$, the system \eqref{de:keller:1}-\eqref{de:keller:3}  and \eqref{keller:1}-\eqref{keller:3}   are parabolic-parabolic models, whereas when $\eps=0$,  the models are parabolic-elliptic. When $\eps \ll 1$, the models are in a transition regime between the parabolic-parabolic and parabolic-elliptic.

In recent years, a large effort has been devoted to  construct energy dissipative,  positivity or bound-preserving schemes \cite{chen2019positivity,chen2022positivity,wang2022fully,li2021unconditionally,tong2024positivity,Shen_2019} for PDEs.
For Keller-Segel equations, the existing approaches can be roughly classified into the following categories: the first approach is  reformulating the problem so that solution of the corresponding discrete problem  is always positive, see \cite{huang2021bound}. The second approach is to construct discrete maximum principle preserving schemes (see, for instance, \cite{liu2018positivity} and the references therein):  these schemes are usually based on reformulating Keller-Segel equations into form of  heat equation $\partial_t \rho =\Grad\cdot(e^c\Grad (\frac{\rho}{e^c}))$.  The  spatial and time  discretization preserving discrete maximum principle should be adopted to achieve positivity. The third approach is to rewrite Keller-Segel equations into a Wasserstein gradient flows \cite{shen2020unconditionally}, the positivity can be guaranteed by  logarithmic type  potentials. The Finite Volume Method, Finite Difference Method, Finite Element Method and Discontinuous Galerkin Method are also proposed to capture the correct behavior of the solution when it is effectively smooth and when it blows-up, see \cite{filbet2006finite,chertock2008second,carrillo2012cross,calvez2006volume,epshteyn2009fully,chertock2018high}.

In this work, we aim to develop numerical schemes which preserve positivity/bound, mass and energy dissipative law. The key ingredient in our approach is  to reformulate the Keller-Segel equation \eqref{de:keller:1}-\eqref{de:keller:3} into an expanded system by using Lagrange multipliers to enforce the constraints of positivity-preserving, mass conserving and energy dissipating.  To efficiently solve the expanded systems, we shall employ an operator splitting approach to decouple Lagrange multipliers from variables of  concentration and density.
A energy correction approach in final step is also introduced to preserve a discrete energy dissipative law.   We would like to point out that  numerical schemes  proposed in this paper  which enjoy the following advantages:
\begin{itemize}
 \item It can be applied directly to any finite difference discretization or other spatial discretization with a Lagrangian basis;
 \item It can be uniquely solved, and it also has essentially the same computational cost as the corresponding semi-implicit  scheme with the same spatial discretization;
 \item It is easy to construct high-order positivity or bound-preserving  schemes;
 \item It has good stability property.
\end{itemize}
To be more specific,  we make the following contributions in this paper:
\begin{itemize}
\item a new class of  positivity/bound-preserving , mass conservative and energy dissipative schemes for two types of Keller-Segel equations are proposed;

\item the second-order, energy dissipative, mass conserving and positivity-preserving  schemes based on Crank-Nicolson and Backward Difference Formulas (BDF2) for the first type of Keller-Segel equations are unconditional stable under small initial condition;

\item  the second-order, bound-preserving, mass-conserving and energy dissipative schemes 
are unconditionally stable for the second type of Keller-Segel equations;

\item the error analysis for the second-order, bound-preserving, mass-conserving and energy dissipative scheme based on Crank-Nicolson is presented.

\end{itemize}

The rest of this paper is organized  as follows. We construct high-order positivity and bound-preserving schemes for Keller-Segel equations \eqref{de:keller:1}-\eqref{de:keller:3} and \eqref{keller:1}-\eqref{keller:3} in Section 2, 3. We make stability analysis for  positivity schemes and bound-preserving schemes in Section 4. Error estimate is  presented  for the second-order, bound-preserving, energy dissipative  scheme based on Crank-Nicolson version for Keller-Segel equations \eqref{keller:1}-\eqref{keller:3} in Section 5. Some numerical simulations are presented in Section 6 to validate the accuracy, stability  and positivity, bound-preserving properties for both types Keller-Segel equations. In Section 7, we give some concluding remarks.

\section{ Structure-preserving  schemes}
We start with a general description of the spatial discretization, followed by the positivity/bound-preserving, mass conserving and energy dissipative  schemes. 

\subsection{Spatial discretization}

We consider a Galerkin type discretization with finite-element method, spectral method or  finite-difference method with summation-by-parts in the subspace $X_h\subset X$, and 
 define a discrete inner product, i.e. numerical integration, on  $\Sigma_h=\{\bm{z}\}$ in $\bar \Omega$: 
\begin{equation}\label{numint}
  [u,v]=\sum_{\bm{z}\in \Sigma_h} \beta_{\bm{z}}u(\bm{z})v(\bm{z}), 
\end{equation}
where we require that the weights $\beta_{\bm{z}}>0$. We also denote the induced norm by $\|u\|=[u,u]^{\frac 12}$. For finite element methods, the sum should be understood as  $\sum_{K\subset \mathcal{T}}\sum_{\bm{z}\in Z(K)}$  where $\mathcal{T}$ is a given triangulation. 
And we assume that there is an unique function $\psi_{\bm{z}}(\bm {x})$ in $X_h$ satisfying $\psi_{\bm{z}}(\bm{z'})=\delta_{\bm{z}\bm{z'}}$ for $\bm{z},\bm{z'}\in \Sigma_h$.  
We denote by $< u_h,v_h>$ the bilinear form on $X_h\times X_h$ based on  the discrete inner product after suitable integration by part.

To construct energy dissipative schemes for two types of Keller-Segel equations,  we denote the free energy of the Keller-Segel equation \eqref{de:keller:1}-\eqref{de:keller:3} as
 \begin{equation}
 E(\rho,c)=\int_{\Omega} \rho\log\rho -\rho-\rho c+\frac 12|\Grad c|^2 d\bx, 
 \end{equation}
 and the free energy of \eqref{keller:1}-\eqref{keller:3} as 
 \begin{equation}
E(\rho,c)=\int_{\Omega} \rho\log\rho +(M-\rho)\log(1-\rho/M)-\rho c+\frac 12|\Grad c|^2 d\bx.
 \end{equation}
 Then the systems \eqref{de:keller:1}-\eqref{de:keller:3} and \eqref{keller:1}-\eqref{keller:3}  satisfy the following energy  dissipative law
 \begin{equation}
 \frac{d}{dt}E(\rho,c)+\int_{\Omega} [\rho|\Grad(\log \rho-c)|^2+|\partial_t c|^2] d\bx=0,
 \end{equation}
 and 
 \begin{equation}
\frac{d}{dt}E(\rho,c)+\int_{\Omega} [\eta(\rho)|\Grad(\log \rho-c)|^2+|\partial_t c|^2] d\bx=0.
 \end{equation}

Firstly, consider a generic spatial discretization of \eqref{de:keller:1}-\eqref{de:keller:3},  to construct energy dissipative, positivity-preserving and mass conserving schemes for \eqref{de:keller:1}-\eqref{de:keller:3}, we can introduce Lagrange multipliers $\lambda_h(\bx,t)$, $\xi_h(t)$ and $\eta_h(t)$ and solve the following expanded system, see also \cite{ChSh22,ChSh_II22}:
\begin{equation}\label{full:dis:1}
\begin{split}
&\eps \partial_tc_h = \mu\Delta c_h + \rho_h+\eta_h,\\
&\partial_t \rho_h -\gamma\Delta \rho_h+\chi\Grad\cdot(\rho_h\Grad c_h)=\lambda_h+\xi_h,\\
&\lambda_h \ge 0,\; \rho_h\ge 0,\;\lambda_h \rho_h=0, \;[\rho_h,1]=[\rho_h^0,1],\\
&\frac{d}{dt}E(\rho_h,c_h)+ [\rho_h|\Grad(\log \rho_h-c_h)|^2+|\partial_t c_h|^2,1] =0,
\end{split}
\end{equation}
where the energy is defined by
\begin{equation}
E(\rho_h,c_h)=[ \rho_h\log\rho_h -\rho_h-\rho_h c_h+\frac 12|\Grad c_h|^2 ,1 ]. 
\end{equation}
Notice that  $\lambda_h$ and $\xi_h$ in \eqref{full:dis:1} are Lagrange multipliers  enforcing physical constraints of positivity-preserving and mass conserving, see \cite{ChSh22,ChSh_II22}.  Scalar auxiliary variable $\eta_h(t)$ is introduced to preserve the energy dissipative law, see \cite{Shen_2019}. Notice that the inequality in \eqref{full:dis:1}  represents the well-known KKT conditions enforcing the constraint of  positivity \cite{ito2008lagrange,facchinei2007finite,harker1990finite,bergounioux1999primal} for constrained systems.

\subsection{Energy correction approach based on CN}
We shall first introduce the energy correction approach to construct structure-preserving schemes based on CN (Crank-Nicolson)  version.  Denote $\delta t$ as  the time step, and $t^n=n\delta t$ for $n=0,1,2,\cdots, \frac{T}{\delta t}$ where $T$ is the final computational time and define
\begin{equation}
\rho_h^{n+\frac 12}=\frac{ \rho_h^{n+1}+\rho^n_h}{2},\quad  \tilde \rho_h^{n+\frac 12}=\frac{\tilde \rho_h^{n+1}+\rho^n_h}{2},\quad  c_h^{n+\frac 12}=\frac{ c_h^{n+1}+c^n_h}{2},\quad  \tilde c_h^{n+\frac 12}=\frac{ \tilde c_h^{n+1}+c^n_h}{2}.
\end{equation}
Firstly, we shall  develop a numerical scheme which satisfies the property of energy dissipating, mass-conserving and positivity-preserving for the Keller-Segel equations \eqref{de:keller:1}-\eqref{de:keller:3}, then consider bound-preserving schemes for \eqref{keller:1}-\eqref{keller:3}  in next subsection.

 Given $c_h^n$ and $\rho_h^n$ and $\lambda_h^n$, we proceed as follows:

{\bf Step 1:} (Predictor) solve $\tilde c_h^{n+1}$ from
\begin{equation}\label{en:positivity:lag:1}
\eps\frac{\tilde c_h^{n+1}(\bz)-c_h^n(\bz)}{\delta t}= \mu\Delta \tilde c_h^{n+\frac 12}(\bz) +\frac 32\rho_h^{n}(\bz)-\frac 12\rho_h^{n-1}(\bz), \quad\forall {\bm z} \in \Sigma_h,
\end{equation}

{\bf Step 2:} (Predictor) solve $\tilde{\rho}_h^{n+1}$ from
\begin{equation}\label{en:positivity:lag:2}
\frac{\tilde{\rho}_h^{n+1}(\bz)-\rho_h^n(\bz)}{\delta t} =\gamma\Delta\tilde{\rho}_h^{n+\frac 12}(\bz)-\chi\Grad\cdot(\tilde{\rho}_h^{n+\frac 12}(\bz)\Grad \tilde c_h^{n+\frac 12}(\bz))+\lambda_h^n(\bz)+\xi_h^n,\quad\forall {\bm z} \in \Sigma_h,
\end{equation}

{\bf Step 3:} (Correction) solve $\rho_h^{n+1}$ and $\lambda_h^{n+1}$ from
\begin{eqnarray}
&&\frac{\rho_h^{n+1}(\bz)-\tilde{\rho}_h^{n+1}(\bz)}{\delta t} =\frac{\lambda_h^{n+1}(\bz)+\xi_h^{n+1}-\lambda_h^n(\bz)-\xi_h^n}{2},\quad\forall {\bm z} \in \Sigma_h,\label{en:positivity:lag:3}\\
&&\rho_h^{n+1}(\bz) \ge 0,\;\lambda_h^{n+1}(\bz)\ge 0,\; \lambda_h^{n+1}(\bz)\rho_h^{n+1}(\bz)=0,\quad\forall {\bm z} \in \Sigma_h,\label{en:positivity:lag:4}\\
&& [\rho^{n+1}_h(\bz),1]=[\rho_h^0(\bz),1].\label{en:positivity:lag:5}
\end{eqnarray}

{\bf Step 4:} (Energy correction) solve $c_h^{n+1}$ and $\eta_h^{n+1}$ from
\begin{eqnarray}
&& \eps\frac{ c_h^{n+1}(\bz)-\tilde c_h^{n+1}(\bz)}{\delta t} = \eta_h^{n+\frac 12},\label{en:positivity:lag:6}\\
&&\frac{E_h^{n+1}-E_h^n}{\delta t}=-[ \rho_h^{n+\frac 12}|\Grad(\log (\rho_h^{n+\frac 12}+\sigma)-\tilde c_h^{n+\frac 12})|^2+|\frac{\tilde  c_h^{n+1}-\tilde c_h^n}{\delta t}|^2,1],\label{en:positivity:lag:7}
\end{eqnarray}
where the energy $E_h^{n+1}$  is defined by 
\begin{equation}\label{en:def}
\begin{split}
E_h^{n+1}=[\rho_h^{n+1}\log( \rho_h^{n+1}+\sigma)-\rho_h^{n+1}-\rho_h^{n+1}  c_h^{n+1}+\frac 12|\Grad  c_h^{n+1}|^2,1].
\end{split}
\end{equation}
The parameter $\sigma \approx 0 $ in \eqref{en:positivity:lag:7} and \eqref{en:def} is introduced to make  $\rho_h^{n+1} +\sigma>0$ which  should be a small positive number.

\begin{lemma}\label{pos:lemma}
For the second-order scheme \eqref{en:positivity:lag:1}-\eqref{en:positivity:lag:7}, Lagrange multipliers $\lambda_h^{n+1}$ and $\xi_h^{n+1}$ satisfy
\begin{equation}
\lambda_h^{n+1} \ge 0, \quad \xi_h^{n+1} \le 0, \quad n=0,1,2,3,\cdots, \frac{T}{\delta t}.
\end{equation}
\end{lemma}
\begin{proof}
Summing up equations \eqref{en:positivity:lag:2} and \eqref{en:positivity:lag:3}, we obtain
\begin{equation}\label{lem:eq:1}
\frac{\rho_h^{n+1}(\bz)-\rho_h^n(\bz)}{\delta t} =\gamma\Delta\tilde{\rho}_h^{n+\frac 12}(\bz)-\chi\Grad\cdot(\tilde{\rho}_h^{n+\frac 12}(\bz)\Grad c_h^{n+\frac 12}(\bz))+\lambda_h^{n+1}(\bz)+\xi_h^{n+1},\quad\forall {\bm z} \in \Sigma_h.
\end{equation}
Taking discrete inner product of \eqref{lem:eq:1} with $1$, we obtain
\begin{equation}
[\lambda_h^{n+1}(\bz)+\xi_h^{n+1},1]=0.
\end{equation}
Noticing $\lambda_h^{n+1} \ge 0 $ from \eqref{en:positivity:lag:4}, we obtain
\begin{equation}
\xi_h^{n+1}=-\frac{[\lambda_h^{n+1},1]}{|\Omega|} \leq 0. 
\end{equation}

\end{proof}

\begin{remark}
Lagrange multiplier $\lambda_h$ defined in \eqref{en:positivity:lag:4} is merely  introduced to enforce the constraint of  postivity-preserving, a situation reminiscent to the pressure which can be reviewed as a Lagrange multiplier to enforce divergence free condition in the Naiver-Stokes equations \cite{Chor68,GMS06,Tema69}.
\end{remark}

\begin{remark}
We introduce local Lagrange multipliers $\lambda(\bx,t)$ instead of global Lagrange multipliers \cite{cheng2020global,cheng2020new,cheng2018multiple}, since the properties of  positivity or bound-preserving are satisfied at every collocation points in the domain.
\end{remark}

\begin{theorem}\label{th:positivity}
Assume that $[\rho_h^{0},1] \neq 0$, then the numerical scheme \eqref{en:positivity:lag:1}-\eqref{en:positivity:lag:7} is unique solvable and satisfies the property of mass-conserving and postivity-preserving for $\rho_h^{n+1} \in X_h$,
\begin{equation}
\rho_h^{n+1}(\bz) \ge 0,\quad  [\rho^{n+1}_h(\bz),1]=[\rho_h^0(\bz),1], \quad \forall {\bm z} \in \Sigma_h.
\end{equation}
It is also energy dissipative in the sense that
\begin{equation}
\frac{E_h^{n+1}-E_h^n}{\delta t}=-[\rho_h^{n+\frac 12}|\Grad(\log (\rho_h^{n+\frac 12}+\sigma)-\tilde c_h^{n+\frac 12})|^2+|\frac{\tilde  c_h^{n+1}-\tilde c_h^n}{\delta t}|^2,1],
\end{equation}
where the energy is defined as \eqref{en:def}.
\end{theorem}
\begin{proof}
The unique solvability of  \eqref{en:positivity:lag:1}-\eqref{en:positivity:lag:2}  is obviously from above steps, since both {\bf Step 1} and {\bf Step 2} are linear schemes for $\tilde c_h^{n+1}$ and $\tilde{\rho}_h^{n+1}$. The total computational cost is exactly the same with IMEX schemes \cite{wang2015stability}.

Introducing $\rho^{\star,n+1}=\tilde \rho^{n+1}-\frac{\delta t}{2}(\lambda_h^n(\bz)+\xi_h^n)$, then equation \eqref{en:positivity:lag:3} can be rewritten into
\begin{equation}\label{eq:bdf1}
\begin{split}
&\frac{\rho_h^{n+1}(\bz)-\rho_h^{\star,n+1}(\bz)}{\delta t} =\frac{\lambda_h^{n+1}(\bz)+\xi_h^{n+1}}{2},\\
&\rho_h^{n+1}(\bz) \ge 0,\;\lambda_h^{n+1}(\bz)\ge 0,\; \lambda_h^{n+1}(\bz)\rho_h^{n+1}(\bz)=0,\quad [\rho^{n+1}_h(\bz),1]=[\rho_h^0(\bz),1].
\end{split}
\end{equation}
From \eqref{eq:bdf1},  we observe that the {\bf Step 3} is equivalent to the following minimization problem
\begin{equation}\label{convex}
\begin{split}
&\min_{\rho_h^{n+1}\in X_h} \frac{\|\rho_h^{n+1}-\rho_h^{\star,n+1}\|^2}{\delta t},\\
&\rho_h^{n+1}(\bz) \ge 0,\;\lambda_h^{n+1}(\bz)\ge 0,\; \lambda_h^{n+1}(\bz)\rho_h^{n+1}(\bz)=0,\quad [\rho^{n+1}_h(\bz),1]=[\rho_h^0(\bz),1].
\end{split}
\end{equation}
We find that the problem \eqref{convex} is a convex minimization problem which admits a unique solution for $\rho_h^{n+1}$.  The constraints of mass conserving and positivity-preserving is enforced exactly in \eqref{en:positivity:lag:4} and \eqref{en:positivity:lag:5}.

 Below we  mainly show how to determine $c_h^{n+1}$ and $\eta_h^{n+1}$ uniquely from \eqref{en:positivity:lag:6}-\eqref{en:positivity:lag:7}. We rewrite \eqref{en:positivity:lag:6} into
\begin{equation}\label{c:eq}
c_h^{n+1}=\tilde c_h^{n+1}+\frac{\delta t}{\eps}\eta_h^{n+\frac 12}. 
\end{equation}
From \eqref{c:eq}, since $\eta_h^{n+\frac 12}$ is a scalar variable,  we obtain
\begin{equation}
[\frac 12|\Grad  c_h^{n+1}|^2,1]=[\frac 12|\Grad  \tilde c_h^{n+1}|^2,1]
\end{equation}
We plug \eqref{c:eq} into \eqref{en:positivity:lag:7} and obtain a  linear  algebraic equation
\begin{equation}\label{th:eq:1}
\begin{split}
&[\rho_h^{n+1}\log( \rho_h^{n+1}+\sigma)-\rho_h^{n+1}-\rho_h^{n+1}  c_h^{n+1}+\frac 12|\Grad  \tilde c_h^{n+1}|^2,1]\\&= E_h^n-\delta t[ \rho_h^{n+\frac 12}|\Grad(\log (\rho_h^{n+\frac 12}+\sigma)-\tilde c_h^{n+\frac 12})|^2+|\frac{\tilde  c_h^{n+1}-\tilde c_h^n}{\delta t}|^2,1]. 
\end{split}
\end{equation}
Then we derive from \eqref{th:eq:1} 
\begin{equation}\label{th:eq:2}
\begin{split}
&[\rho_h^{n+1}c_h^{n+1},1]=E_h^n-\delta t[ \rho_h^{n+\frac 12}|\Grad(\log (\rho_h^{n+\frac 12}+\sigma)-\tilde c_h^{n+\frac 12})|^2+|\frac{\tilde  c_h^{n+1}-\tilde c_h^n}{\delta t}|^2,1]\\&-[\rho_h^{n+1}\log( \rho_h^{n+1}+\sigma)-\rho_h^{n+1}+\frac 12|\Grad  \tilde c_h^{n+1}|^2,1].
\end{split}
\end{equation}
Plugging \eqref{c:eq} into \eqref{th:eq:2}, we obtain
\begin{equation}
\begin{split}
&[\rho_h^{n+1}\tilde c_h^{n+1}+\frac{\delta t}{\eps}\rho_h^{n+1}\eta_h^{n+1},1]=E_h^n-\delta t[ \rho_h^{n+\frac 12}|\Grad(\log (\rho_h^{n+\frac 12}+\sigma)-\tilde c_h^{n+\frac 12})|^2+|\frac{\tilde  c_h^{n+1}-\tilde c_h^n}{\delta t}|^2,1]\\&-[\rho_h^{n+1}\log( \rho_h^{n+1}+\sigma)-\rho_h^{n+1}+\frac 12|\Grad  \tilde c_h^{n+1}|^2,1].
\end{split}
\end{equation}
Finally, notice that $\rho_h^{n+\frac 12}$ is known,  we solve $\eta_h^{n+\frac 12}$ from the following linear equation 
\begin{equation}\label{v:eta}
\begin{split}
&\eta_h^{n+\frac 12}=\frac{E_h^n-\delta t[ \rho_h^{n+\frac 12}|\Grad(\log (\rho_h^{n+\frac 12}+\sigma)-\tilde c_h^{n+\frac 12})|^2+|\frac{\tilde  c_h^{n+1}-\tilde c_h^n}{\delta t}|^2,1]}{[\frac{\delta t}{\eps}\rho_h^{n+1},1]}\\&-\frac{[\rho_h^{n+1}\log( \rho_h^{n+1}+\sigma)-\rho_h^{n+1}-\rho_h^{n+1}\tilde c_h^{n+1}+\frac 12|\Grad  \tilde c_h^{n+1}|^2,1]}{[\frac{\delta t}{\eps}\rho_h^{n+1},1]}.
\end{split}
\end{equation}
From the property of mass-conserving, we have
\begin{equation}
[\rho_h^{n+1},1] =[\rho_h^{0},1] \neq 0. 
\end{equation}
The energy dissipative law is easily observed from the scheme \eqref{en:positivity:lag:1}-\eqref{en:positivity:lag:7}.

\end{proof}

\begin{lemma}
From Theorem \ref{th:positivity}, we observe that the second-order scheme \eqref{en:positivity:lag:1}-\eqref{en:positivity:lag:7} can be viewed as  operator splitting scheme and {\bf Step 3} is regarded as projection step, .ie. projecting $\tilde{\rho}_h^{n+1}$ in $L^2$ norm into positivity-preserving space, 
 a situation reminiscent to the projection into divergence free step for  Naiver-Stokes equation \cite{shen1992pressure,GMS06,shen1992error,shen1992errornew}.
\end{lemma}

\begin{lemma}
For the Crank-Nicolson  scheme \eqref{en:positivity:lag:1}-\eqref{en:positivity:lag:7},  the Lagrange multiplier $\eta_h^{n+\frac 12}$  is  second-order  approximation to $0$ in the sense that
\begin{equation}
\eta_h^{n+\frac 12} = c\delta t^2, \quad n=0,1,2,3,\cdots .
\end{equation}
where $c$ is a constant independent of $\delta t$. 
\end{lemma}
\begin{proof}
According to  Taylor expansion for \eqref{en:positivity:lag:7}, we can easily derive 
\begin{equation}
\begin{split}
&E_h^n-[\rho_h^{n+1}\log( \rho_h^{n+1}+\sigma)-\rho_h^{n+1}-\rho_h^{n+1}\tilde c_h^{n+1}+\frac 12|\Grad  \tilde c_h^{n+1}|^2,1] \\&-\delta t[\tilde \rho_h^{n+\frac 12}|\Grad(\log (\rho_h^{n+\frac 12}+\sigma)-\tilde c_h^{n+\frac 12})|^2+|\frac{\tilde  c_h^{n+1}-\tilde c_h^n}{\delta t}|^2 =O(\delta t^3). 
\end{split}
\end{equation}
Then from \eqref{v:eta}, we obtain
\begin{equation}
\eta_h^{n+\frac 12}=O(\delta t^2).
\end{equation}

\end{proof}

\subsection{Energy correction approach based on BDFk}
We shall construct high-order energy dissipative, positivity-preserving and mass conserving  schemes by using Backward Difference Formulas (BDF) \cite{ChSh22,hou2023linear}. More precisely,  a $k$th order IMEX schemes based on BDFk \cite{hou2023linear} and Adam-Bashforth are proceeded as follows:

given $A_k(c_h^n)$, $A_k(\rho_h^n)$, $C_k(\rho_h^n)$ and $B_{k-1}(\lambda_h^n)$, $B_{k-1}(\xi_h^n)$, 

{\bf Step 1:} (Predictor) solve $\tilde c_h^{n+1}$ from
\begin{equation}\label{high:positivity:lag:1}
\eps\frac{\alpha_k \tilde c_h^{n+1}(\bz)-A_k(c_h^n(\bz))}{\delta t}= \mu\Delta \tilde c_h^{n+1}(\bz) + C_k(\rho_h^{n}(\bz)), \quad\forall {\bm z} \in \Sigma_h,
\end{equation}

{\bf Step 2:} (Predictor) solve $\tilde{\rho}_h^{n+1}$ from
\begin{equation}\label{high:positivity:lag:2}
\begin{split}
&\frac{\alpha_k\tilde{\rho}_h^{n+1}(\bz)-A_k(\rho_h^n(\bz))}{\delta t} =\gamma\Delta\tilde{\rho}_h^{n+1}(\bz)\\&-\chi\Grad\cdot(\tilde{\rho}_h^{n+1}(\bz)\Grad \tilde c_h^{n+1}(\bz)) +B_{k-1}(\lambda_h^n+\xi_h^n), \quad\forall {\bm z} \in \Sigma_h,
\end{split}
\end{equation}
 
 {\bf Step 3:} (Correction) solve $\rho_h^{n+1}$ and $\lambda_h^{n+1}$, $\xi_h^{n+1}$ from
\begin{eqnarray}
&&\frac{\alpha_k(\rho_h^{n+1}(\bz)-\tilde{\rho}_h^{n+1}(\bz))}{\delta t} =\lambda_h^{n+1}(\bz)+\xi_h^{n+1}-B_{k-1}(\lambda_h^n+\xi_h^n),\quad\forall {\bm z} \in \Sigma_h,\label{high:positivity:lag:3}\\
&&\rho_h^{n+1}(\bz) \ge 0,\;\lambda_h^{n+1}(\bz)\ge 0,\; \lambda_h^{n+1}(\bz)\rho_h^{n+1}(\bz)=0,[\rho^{n+1}_h(\bz),1]=[\rho_h^0(\bz),1],\quad\forall {\bm z} \in \Sigma_h,\label{high:positivity:lag:4}
\end{eqnarray}

{\bf Step 4:} (Energy correction) solve $c_h^{n+1}$ and $\eta_h^{n+1}$ from
\begin{eqnarray}
&& \eps\frac{ \alpha_k c_h^{n+1}(\bz)-\alpha_k\tilde c_h^{n+1}(\bz)}{\delta t} = \eta_h^{n+1},\label{high:positivity:lag:5}\\
&&\frac{\alpha_k E_h^{n+1}-A_k(E_h^n)}{\delta t}=-[ \rho_h^{n+1}|\Grad(\log (\rho_h^{n+1}+\sigma)-\tilde c_h^{n+1})|^2+|\frac{\alpha_k\tilde  c_h^{n+1}-A_k(\tilde c_h^n)}{\delta t}|^2,1],\label{high:positivity:lag:6}
\end{eqnarray}
where the energy $E_h^{n+1}$  is defined by 
\begin{equation}\label{high:en:def}
\begin{split}
E_h^{n+1}=[\rho_h^{n+1}\log( \rho_h^{n+1}+\sigma)-\rho_h^{n+1}-\rho_h^{n+1}  c_h^{n+1}+\frac 12|\Grad  c_h^{n+1}|^2,1].
\end{split}
\end{equation}
 Notice that for any function $u_h^n$, the  parameter $\alpha_k,$ and  the operators $A_k$, $C_k$ $(k=1,2,3,4)$, $B_k$ $(k=0,1,2,3)$ in \eqref{high:positivity:lag:1}-\eqref{high:positivity:lag:6}  are defined  by:

\noindent {\bf First-order:}
\begin{equation}\label{eq:positivity:bdf1}
\alpha_1=1, \quad A_1(u_h^n)=u_h^n,\quad B_0(u_h^n)=u_h^n,\quad C_1(u_h^n)=u_h^n;
\end{equation}
\noindent {\bf Second-order:}
\begin{equation}\label{eq:positivity:bdf2}
\alpha_2=\frac{3}{2}, \quad A_2(u_h^n)=2u_h^n-\frac{1}{2}u_h^{n-1},\quad B_1(u_h^n)=u_h^n,\quad C_2(u_h^n)=2 u_h^n-u_h^n;
\end{equation}
\noindent {\bf Third-order:}
\begin{equation}\label{eq:positivity:bdf3}
\begin{split}
&\alpha_3=\frac{11}{6}, \quad A_3(u_h^n)=3u_h^n-\frac{3}{2}u_h^{n-1}+\frac{1}{3}u_h^{n-2},\\& B_2(u_h^n)=2u_h^n-u_h^{n-1},\quad C_3(u_h^n)=3 c_h^n-3 c_h^{n-1} + c_h^{n-2};
\end{split}
\end{equation}
\noindent {\bf Fourth-order:}
\begin{equation}\label{eq:positivity:bdf4}
\begin{split}
&\alpha_4=\frac{25}{12}, \; A_4(u_h^n)=4u_h^n-3u_h^{n-1}+\frac{4}{3}u_h^{n-2}-\frac{1}{4}u_h^{n-3},\;\\& B_3(u_h^n)=3u_h^n-3u_h^{n-1}+u_h^{n-2},\quad C_4(u_h^n)=4 u_h^n-6 u_h^{n-1}+4 u_h^{n-2} -u_h^{n-3},
\end{split}
\end{equation}
for any function $u_h^n$.
The formulas of BDF$k$ for $k=5,6$ can be derived similarly with Taylor expansions.

\begin{theorem}
The  $kth$ order  schemes \eqref{high:positivity:lag:1}-\eqref{high:positivity:lag:6} for $k=1,2,3,4,5$ are unique solvable and satisfy the property of mass conserving and positivity-preserving for $\rho_h^{n+1} \in X_h$, .i.e.
\begin{equation}
\rho_h^{n+1}(\bz) \ge 0,\quad [\rho_h^{n+1}(\bz),1]=[\rho_h^0(\bz),1],\quad  \forall {\bm z} \in \Sigma_h. 
\end{equation}
It also satisfies the following discrete energy dissipative law
\begin{equation}
\begin{split}
\frac{\alpha_k E_h^{n+1}-A_k(E_h^n)}{\delta t}=-[ \rho_h^{n+1}|\Grad(\log (\rho_h^{n+1}+\sigma)-\tilde c_h^{n+1})|^2+|\frac{\alpha_k\tilde  c_h^{n+1}-A_k(\tilde c_h^n)}{\delta t}|^2,1],
\end{split}
\end{equation}
where $E_h^{n+1}$ is defined as \eqref{high:en:def}. 
\end{theorem}
\begin{proof}
We  observe {\bf Step 1} and {\bf Step 2} are usual $k$th-order IMEX schemes \cite{wang2015stability} which can be solved uniquely and efficiently.  Introducing $\rho_h^{\star,n+1}=\tilde\rho_h^{n+1}-\frac{\delta t}{\alpha_k}B_{k-1}(\lambda_h^n+\xi_h^n)$, we rewrite \eqref{high:positivity:lag:3} into 
\begin{equation}\label{re:2}
\frac{\alpha_k(\rho_h^{n+1}(\bz)-\rho_h^{\star,n+1}(\bz))}{\delta t} =\lambda_h^{n+1}(\bz)+\xi_h^{n+1}
\end{equation}
Using  \eqref{re:2} ,  similarly, we observe that the  {\bf Step 3} \eqref{high:positivity:lag:3}-\eqref{high:positivity:lag:4}   is equivalent to the following convex minimization problem
\begin{equation}\label{convex:2}
\begin{split}
&\min_{\rho_h^{n+1}\in X_h} \frac{\alpha_k\|\rho_h^{n+1}-\rho_h^{\star,n+1}\|^2}{2\delta t},\\
&\rho_h^{n+1}(\bz) \ge 0,\;\lambda_h^{n+1}(\bz)\ge 0,\; \lambda_h^{n+1}(\bz)\rho_h^{n+1}(\bz)=0,\quad [\rho^{n+1}_h(\bz),1]=[\rho_h^0(\bz),1].
\end{split}
\end{equation}
Then the uniqueness of $\rho_h^{n+1}$ can be easily obtained. Below, we show how to determine $\eta_h^{n+1}$ and $c_h^{n+1}$ uniquely from \eqref{high:positivity:lag:5}-\eqref{high:positivity:lag:6}.

Similar to the Crank-Nicolson scheme \eqref{en:positivity:lag:1}-\eqref{en:positivity:lag:7}, we can derive the following algebraic equation
\begin{equation}
\begin{split}
&[\rho_h^{n+1}\log( \rho_h^{n+1}+\sigma)-\rho_h^{n+1}-\rho_h^{n+1}  c_h^{n+1}+\frac 12|\Grad  \tilde c_h^{n+1}|^2,1]\\&= \frac{1}{\alpha_k}A_k(E_h^n)-\frac{\delta t}{\alpha_k}[ \rho_h^{n+1}|\Grad(\log (\rho_h^{n+1}+\sigma)-\tilde c_h^{n+1})|^2+|\frac{\alpha_k\tilde  c_h^{n+1}-A_k(\tilde c_h^n)}{\delta t}|^2,1]. 
\end{split}
\end{equation}
We solve $\eta_h^{n+1}$ from the following linear equation 
\begin{equation}
\begin{split}
&\eta_h^{n+1}=\frac{\frac{1}{\alpha_k}A_k(E_h^n)-\frac{\delta t}{\alpha_k}[ \rho_h^{n+\frac 12}|\Grad(\log (\rho_h^{n+1}+\sigma)-\tilde c_h^{n+1})|^2+|\frac{\alpha_k\tilde  c_h^{n+1}-A_k(\tilde c_h^n)}{\delta t}|^2,1]}{[\frac{\delta t}{\alpha_k\eps}\rho_h^{n+1},1]}\\&-\frac{[\rho_h^{n+1}\log( \rho_h^{n+1}+\sigma)-\rho_h^{n+1}-\rho_h^{n+1}\tilde c_h^{n+1}+\frac 12|\Grad  \tilde c_h^{n+1}|^2,1]}{[\frac{\delta t}{\alpha_k\eps}\rho_h^{n+1},1]}.
\end{split}
\end{equation}

\end{proof}
The stability result for the  $kth$ order  scheme \eqref{high:positivity:lag:1}-\eqref{high:positivity:lag:6} will be demonstrated  in next section. 

\begin{remark}
For the  $kth$ order  schemes \eqref{high:positivity:lag:1}-\eqref{high:positivity:lag:6},  Lagrange multipliers $\lambda_h^{n+1}$ and $\xi_h^{n+1}$ also  satisfy
\begin{equation}
\lambda_h^{n+1} \ge 0, \quad \xi_h^{n+1} \le 0, \quad n=0,1,2,3,\cdots, \frac{T}{\delta t}.
\end{equation}
The proof is similar to the Crank-Nicolson scheme \eqref{en:positivity:lag:1}-\eqref{en:positivity:lag:7}. We omit it here. 
\end{remark}

\section{Bound-preserving schemes}
In this section, we shall construct  bound-preserving, mass conserving and energy dissipative schemes for Keller-Segel equations \eqref{keller:1}-\eqref{keller:3}. Similar to positivity-preserving schemes, consider a generic spatial discretization of \eqref{keller:1}-\eqref{keller:3}, we rewrite it 
 into the following equivalent  expanded system, see also \cite{ChSh_II22}
\begin{equation}\label{full:dis:2}
\begin{split}
&\eps \partial_tc_h = \mu\Delta c_h + \rho_h+\eta_h(t),\\
&\partial_t \rho_h -\gamma\Delta \rho_h+\chi\Grad\cdot(\eta(\rho_h)\Grad c_h)=\lambda_hg'(\rho_h)+\xi_h(t),\\
&\lambda_h \ge 0,\; g(\rho_h)\ge 0,\;\lambda_h g(\rho_h)=0,\;[\rho_h(t),1]=[\rho_h^0,1],\\
&\frac{d}{dt}E(\rho_h,c_h)+ [\eta(\rho_h)|\Grad(\log \rho_h-c_h)|^2+|\partial_t c_h|^2,1] =0,
\end{split}
\end{equation}
where the energy $E_h(\rho_h,c_h)$ is defined by 
\begin{equation}\label{law:2:con}
E_h(\rho_h,c_h)=[\rho_h\log\rho_h +(M-\rho_h)\log(1-\rho_h/M)-\rho_h c_h+\frac 12|\Grad c_h|^2 ,1].
 \end{equation}
To preserve $0\leq \rho_h\leq M$, we introduce the function $g(\rho_h)$ to be $g(\rho_h)=\rho_h(M-\rho_h)$ in \eqref{full:dis:2}.  The third equation in \eqref{full:dis:2} is the KKT-condition used to enforce the constraints of  bound-preserving and mass conserving, see \cite{ChSh_II22}. Lagrange mutliplier $\eta_h(t)$ is introduced to preserve the energy dissipative law \eqref{law:2:con} in the energy correction step. 
Similar to positivity-preserving schemes, we shall construct bound-preserving schemes based on the expanded system \eqref{full:dis:2}.

\subsection{Mass conserving and bound-preserving schemes with BDFk}
In this subsection, we construct high-order bound-preserving  schemes for Keller-Segel equations \eqref{keller:1}-\eqref{keller:3} based on the  predictor-correction approach.

Given $A_k(c_h^n)$, $A_k(\rho_h^n)$, $C_k(\rho_h^n)$ and $B_{k-1}(\lambda_h^n)$, $B_{k-1}(\xi_h^n)$, we proceed as follows:

{\bf Step 1:} (Prediction) solve $\tilde c_h^{n+1}$ from
\begin{equation}\label{keller:lag:1}
\eps\frac{\alpha_k \tilde c_h^{n+1}(\bz)-A_k(c_h^n(\bz))}{\delta t}= \mu\Delta \tilde c_h^{n+1}(\bz) + C_k(\rho_h^{n}(\bz)), \quad\forall {\bm z} \in \Sigma_h,
\end{equation}

{\bf Step 2:} (Prediction) solve $\tilde{\rho}_h^{n+1}$ from
\begin{equation}\label{keller:lag:2}
\begin{split}
&\frac{\alpha_k\tilde{\rho}_h^{n+1}(\bz)-A_k(\rho_h^n(\bz))}{\delta t} =\gamma\Delta\tilde{\rho}_h^{n+1}(\bz)\\&-\chi\Grad\cdot(\eta(C_k(\rho_h^{n}(\bz)))\Grad C_k(c_h^{n}(\bz))) +B_{k-1}(\lambda_h^ng'(\rho_h^n)+\xi_h^{n}), \quad\forall {\bm z} \in \Sigma_h,
\end{split}
\end{equation}
 
 {\bf Step 3:} (Correction) solve $\rho_h^{n+1}$ and $\lambda_h^{n+1}, \xi_h^{n+1}$ from
\begin{eqnarray}
&&\frac{\alpha_k(\rho_h^{n+1}(\bz)-\tilde{\rho}_h^{n+1}(\bz))}{\delta t} =\lambda_h^{n+1}(\bz)g'(\rho_h^{n+1}(\bz))+\xi_h^{n+1}-B_{k-1}(\lambda_h^ng'(\rho_h^n)+\xi_h^n),\quad\forall {\bm z} \in \Sigma_h,\label{keller:lag:3}\\
&&g(\rho_h^{n+1}(\bz)) \ge 0,\;\lambda_h^{n+1}(\bz)\ge 0,\; \lambda_h^{n+1}(\bz)g(\rho_h^{n+1}(\bz))=0,\; [\rho^{n+1}_h(\bz),1]=[\rho_h^0(\bz),1], \quad\forall {\bm z} \in \Sigma_h,\label{keller:lag:4}
\end{eqnarray}

{\bf Step 4:} (Energy correction) solve $c_h^{n+1}$ and $\eta_h^{n+1}$ from
\begin{eqnarray}
&& \eps\frac{ \alpha_k c_h^{n+1}(\bz)-\alpha_k\tilde c_h^{n+1}(\bz)}{\delta t} = \eta_h^{n+1},\label{keller:lag:5}\\
&&\frac{\alpha_k E_h^{n+1}-A_k(E_h^n)}{\delta t}=-[ \eta(\rho_h^{n+1})|\Grad(\log (\rho_h^{n+1}+\sigma)-\tilde c_h^{n+1})|^2+|\frac{\alpha_k\tilde  c_h^{n+1}-A_k(\tilde c_h^n)}{\delta t}|^2,1],\label{keller:lag:6}
\end{eqnarray}
where the energy $E_h^{n+1}$  is defined by 
\begin{equation}\label{keller:en:def}
\begin{split}
E_h^{n+1}=[\rho_h^{n+1}\log(\rho_h^{n+1}+\sigma) +(M-\rho^{n+1}_h)\log(1-\rho^{n+1}_h/M+\sigma)-\rho^{n+1}_h c^{n+1}_h+\frac 12|\Grad c^{n+1}_h|^2 ,1].
\end{split}
\end{equation}
Notice that  $\alpha_k$, $A_k$, $B_k$ and $C_k$ are defined in \eqref{eq:positivity:bdf1}-\eqref{eq:positivity:bdf4}.
We also introduce a small positive parameter $\sigma \approx 0$ in \eqref{keller:en:def} to satisfy that  $\rho_h^{n+1}+\sigma>0$ and $1-\rho^{n+1}_h/M+\sigma>0$. 

\begin{theorem}\label{lem2:f:bound}
The  $kth$ order  schemes \eqref{keller:lag:1}-\eqref{keller:lag:4} for $k=1,2,3,4,5$ are unique solvable and satisfy the property of mass conserving and bound-preserving for $\rho_h^{n+1} \in X_h$,
\begin{equation}
0\leq \rho_h^{n+1}(\bz) \le M,\quad [\rho_h^{n+1}(\bz),1]=[\rho_h^0(\bz),1],\quad  \forall {\bm z} \in \Sigma_h.
\end{equation}
It also satisfies the discrete energy dissipative law in the sense that
\begin{equation}
\frac{\alpha_k E_h^{n+1}-A_k(E_h^n)}{\delta t}=-[ \eta(\rho_h^{n+1})|\Grad(\log (\rho_h^{n+1}+\sigma)-\tilde c_h^{n+1})|^2+|\frac{\alpha_k\tilde  c_h^{n+1}-A_k(\tilde c_h^n)}{\delta t}|^2,1].
\end{equation}
\end{theorem}
\begin{proof}
 Introducing $\rho_h^{\star,n+1}=\tilde\rho_h^{n+1}-\frac{\delta t}{\alpha_k}B_{k-1}(\lambda_h^ng'(\rho_h^n)+\xi_h^n)$, we rewrite \eqref{keller:lag:3} into 
\begin{equation}\label{re:2:bound}
\frac{\alpha_k(\rho_h^{n+1}(\bz)-\rho_h^{\star,n+1}(\bz))}{\delta t} =\lambda_h^{n+1}(\bz)g'(\rho_h^{n+1}(\bz))+\xi_h^{n+1}.
\end{equation}
Using  \eqref{re:2:bound},  it is easy to observe that  {\bf Step 3} \eqref{keller:lag:3}-\eqref{keller:lag:4}   is equivalent to the following convex minimization problem
\begin{equation}\label{convex:2:bound}
\begin{split}
&\min_{\rho_h^{n+1}\in X_h} \frac{\alpha_k\|\rho_h^{n+1}-\rho_h^{\star,n+1}\|^2}{2\delta t},\\
&g(\rho_h^{n+1}(\bz)) \ge 0,\;\lambda_h^{n+1}(\bz)\ge 0,\; \lambda_h^{n+1}(\bz)g(\rho_h^{n+1}(\bz))=0,\quad [\rho^{n+1}_h(\bz),1]=[\rho_h^0(\bz),1].
\end{split}
\end{equation}
Then the uniqueness can be easily obtained. The discrete energy dissipative law can be easily observed from \eqref{keller:en:def}. 

\end{proof}

\subsection{Bound-preserving scheme based on CN}
We can also construct a second-order,  energy dissipative, mass conserving and bound-preserving scheme for \eqref{keller:1}-\eqref{keller:3} based on Crank-Nicolson version. Given $c_h^n, c_h^{n-1}$ and $\rho_h^n,\rho_h^{n-1}$, the CN scheme is proceeded as follows

{\bf Step 1:} (Prediction) solve $\tilde{c}_h^{n+1}$ from
\begin{equation}\label{en:keller:lag:1}
\eps\frac{\tilde{c}_h^{n+1}(\bz)-c_h^n(\bz)}{\delta t}= \mu\Delta \tilde{c}_h^{n+\frac 12}(\bz) + \frac32\rho_h^{n}(\bz)-\frac 12\rho_h^{n-1}(\bz), \quad\forall {\bm z} \in \Sigma_h,
\end{equation}

{\bf Step 2:} (Prediction) solve $\tilde{\rho}_h^{n+1}$ from
\begin{equation}\label{en:keller:lag:2}
\begin{split}
&\frac{\tilde{\rho}_h^{n+1}(\bz)-\rho_h^n(\bz)}{\delta t} =\gamma\Delta\tilde{\rho}_h^{n+\frac 12}(\bz)\\&-\chi\Grad\cdot(\eta(\frac 32 \rho_h^{n}(\bz)-\frac 12\rho_h^{n-1}(\bz))\Grad (\frac 32c_h^{n}-\frac 12c_h^{n-1})) +\lambda_h^ng'(\rho_h^n)+\xi_h^n, \quad\forall {\bm z} \in \Sigma_h,
\end{split}
\end{equation}
 
 {\bf Step 3:} (Bound-correction) solve $\rho_h^{n+1}$ and $\lambda_h^{n+1}, \xi_h^{n+1}$ from
\begin{eqnarray}
&&\frac{\rho_h^{n+1}(\bz)-\tilde{\rho}_h^{n+1}(\bz)}{\delta t} =\frac{\lambda_h^{n+1}(\bz)g'(\rho_h^{n+1}(\bz))+\xi_h^{n+1}-\lambda_h^ng'(\rho_h^n)-\xi_h^n}{2},\label{en:keller:lag:3}\\
&&g(\rho_h^{n+1}(\bz)) \ge 0,\;\lambda_h^{n+1}(\bz)\ge 0,\; \lambda_h^{n+1}(\bz)g(\rho_h^{n+1}(\bz))=0,\quad\forall {\bm z} \in \Sigma_h,\label{en:keller:lag:4}\\
&&[\rho^{n+1}_h(\bz),1]=[\rho_h^0(\bz),1], \quad\forall {\bm z} \in \Sigma_h.
\end{eqnarray}

{\bf Step 4:} (Energy-correction)  solve $c_h^{n+1}$ and $\eta_h^{n+\frac 12}$ from
\begin{eqnarray}
&& \eps\frac{ c_h^{n+1}(\bz)-\tilde c_h^{n+1}(\bz)}{\delta t} = \eta_h^{n+\frac 12},\label{en:keller:lag:5}\\
&&\frac{E_h^{n+1}-E_h^n}{\delta t}=-[ \eta(\rho_h^{n+\frac 12})|\Grad(\log (\rho_h^{n+\frac 12}+\sigma)-\tilde c_h^{n+\frac 12})|^2+|\frac{\tilde  c_h^{n+1}-\tilde c_h^n}{\delta t}|^2,1],\label{en:keller:lag:6}
\end{eqnarray}
where the energy $E_h^{n+1}$  is defined by \eqref{keller:en:def}.

\begin{theorem}\label{th:bound:u}
Assume that $[\rho_h^{0},1] \neq 0$, then the numerical scheme \eqref{en:keller:lag:1}-\eqref{en:keller:lag:6} is  unique solvable and satisfies the property of mass-conserving and bound-preserving for $\rho_h^{n+1} \in X_h$,
\begin{equation}
0\le \rho_h^{n+1}(\bz) \le M,\quad  [\rho^{n+1}_h(\bz),1]=[\rho_h^0(\bz),1], \quad \forall {\bm z} \in \Sigma_h.
\end{equation}
It is also energy dissipative in the sense that
\begin{equation}
\frac{E_h^{n+1}-E_h^n}{\delta t}=-[ \eta(\rho_h^{n+\frac 12})|\Grad(\log (\rho_h^{n+\frac 12}+\sigma)-\tilde c_h^{n+\frac 12})|^2+|\frac{\tilde  c_h^{n+1}-\tilde c_h^n}{\delta t}|^2,1],
\end{equation}
where the energy is defined by \eqref{keller:en:def}.
\end{theorem}
\begin{proof}
Since the proof is similar to energy dissipative and positivity-preserving scheme proposed above, we omit the proof  here. 
\end{proof}

\section{Stability results}

In this section,  we shall derive the stability result for positivity-preserving and bound-preserving schemes we proposed in last section.  The positivity-preserving and energy dissipative schemes proposed in above section only imply the bound for  energy, it does not imply the $L^2$ norm or $H^1$ norm for $\rho_h^{n+1}$ and $\lambda^{n+1}$. While the bound-preserving schemes constructed  automatically ensure the $L^\infty$ bound for $\{\rho_h^{n+1}\}$, it does not imply any bound on the $\{c_h^{n+1}\}$ and $\lambda_h^{n+1}$. In this section, we shall use the energy estimates to derive a bound on the  norm for $\{c_h^{n+1}\}$ as well as a bound on the Lagrange multipliers. 

We shall frequently use the following discrete  Gronwall Lemma  (P.431 in  \cite{shen2011spectral}):
 \begin{lemma}\label{Gron2}
 Let $a_n,\,b_n,\,c_n,$ and $d_n$ be four nonnegative sequences satisfying \begin{equation*}
		a_m+\tau \sum_{n=1}^{m} b_n \le \tau \sum_{n=0}^{m}a_n d_n +\tau \sum_{n=0}^{m} c_n+ C, \mbox{with} \;\delta t\sum\limits_{n=0}^{T/\delta t} d_n\leq M,\; \forall 0\leq m \leq T/\delta t,
	\end{equation*}
	where $C$ and $\tau$ are two positive constants and  we assume $\delta td_n <1$ and  let $\sigma =\max_{0\leq n\leq \frac{T}{\delta t}}(1-\delta td_n)^{-1}$.
	Then
	\begin{equation*}
		a_m+\tau \sum_{n=1}^{m} b_n \le \exp\big(\sigma M \big)\big(\tau \sum_{n=0}^{m}c_n+C \big),\; \forall m \le \frac{T}{\delta t}.
	\end{equation*}
\end{lemma}

For simplicity and without loss of generality, we consider periodic boundary condition and set $\eps=\mu=\gamma=\chi=1$ for both types of Keller-Segel equations when we make stability results and error estimate for our proposed schemes.

\subsection{Stability of   positivity preserving schemes}

\begin{theorem}\label{bdf2:stab:result}
Given  $\eps=\gamma=\chi=\mu=1$ and the final time $T$, then for $m\delta t\leq T$, the second-order, energy dissipative, mass conserving and  positivity preserving scheme \eqref{high:positivity:lag:1}-\eqref{high:positivity:lag:6} for Keller-Segel equation \eqref{de:keller:1}-\eqref{de:keller:3}  is stable in 2D  if the small data condition 
 \begin{equation}
4\|\rho_h^{1}\|^2+\|2\rho_h^{1}-\rho_h^{0}\|^2+\|\Grad c_h^{1}\|^2 +\|2\Grad c_h^{1}-\Grad c_h^0\|^2
\leq e^{-5T},
\end{equation}
is satisfied.
\end{theorem}
\begin{proof}
Taking inner product of equation \eqref{high:positivity:lag:2} with $4\delta t \tilde{\rho}_h^{n+1}$,we obtain
\begin{equation}\label{bdf2:stab:1}
\begin{split}
[3\tilde{\rho}_h^{n+1}-4\rho_h^n+\rho_h^{n-1},2\tilde{\rho}_h^{n+1}]
+4\delta t\|\Grad\tilde{\rho}_h^{n+1}\|^2=-2\delta t[(\tilde{\rho}_h^{n+1})^2,\Delta \tilde c_h^{n+1}]+4\delta t[\lambda_h^n+\xi_h^n,\tilde{\rho}_h^{n+1}],
\end{split}
\end{equation}
where we used 
\begin{equation}
[\Grad\cdot(\tilde \rho_h^{n+1}\Grad \tilde c_h^{n+1}), \tilde \rho_h^{n+1}]= -[\tilde \rho_h^{n+1}\Grad \tilde c_h^{n+1}, \Grad \tilde \rho_h^{n+1}]=[\tilde \rho_h^{n+1}\Grad \tilde c_h^{n+1}, \Grad \tilde \rho_h^{n+1}]+[(\tilde{\rho}_h^{n+1})^2,\Delta \tilde c_h^{n+1}].
\end{equation}

Now we consider left term in \eqref{bdf2:stab:1},
\begin{equation}\label{bdf2:stab:2}
\begin{split}
&[3\tilde{\rho}_h^{n+1}-4\rho_h^n+\rho_h^{n-1},2\tilde{\rho}_h^{n+1}]=
[3\rho_h^{n+1}-4\rho_h^n+\rho_h^{n-1},2\rho_h^{n+1}]\\&
+2[3\rho_h^{n+1}-4\rho_h^n+\rho_h^{n-1},\tilde{\rho}_h^{n+1}-\rho_h^{n+1}]+6[\tilde{\rho}_h^{n+1}-\rho_h^{n+1},\tilde{\rho}_h^{n+1}].
\end{split}
\end{equation}
We calculate each term in \eqref{bdf2:stab:2},
\begin{equation}\label{bdf2:stab:2:e}
\begin{split}
&[3\rho_h^{n+1}-4\rho_h^n+\rho_h^{n-1},2\rho_h^{n+1}]=
\|\rho_h^{n+1}\|^2-\|\rho_h^n\|^2\\&+\|2\rho_h^{n+1}-\rho_h^n\|^2-\|2\rho_h^n-\rho_h^{n-1}\|^2+\|\rho_h^{n+1}-2\rho_h^n+\rho_h^{n-1}\|^2.
\end{split}
\end{equation}
Similarly, we have 
\begin{equation}\label{bdf2:stab:3}
\begin{split}
&2[3\rho_h^{n+1}-4\rho_h^n+\rho_h^{n-1},\tilde{\rho}_h^{n+1}-\rho_h^{n+1}]
=2[\rho_h^{n+1}-2\rho_h^n+\rho_h^{n-1},\tilde{\rho}_h^{n+1}-\rho_h^{n+1}]\\&+
4[\rho_h^{n+1}-\rho_h^n,\tilde{\rho}_h^{n+1}-\rho_h^{n+1}].
\end{split}
\end{equation}
We consider the first term in right hand side of \eqref{bdf2:stab:3}, using Young inequality
\begin{equation}\label{bdf2:stab:3:e}
\begin{split}
2[\rho_h^{n+1}-2\rho_h^n+\rho_h^{n-1},\tilde{\rho}_h^{n+1}-\rho_h^{n+1}]\leq \|\rho_h^{n+1}-2\rho_h^n+\rho_h^{n-1}\|^2+\|\tilde{\rho}_h^{n+1}-\rho_h^{n+1}\|^2.
\end{split}
\end{equation}
We consider the last term in \eqref{bdf2:stab:3}, using \eqref{high:positivity:lag:3}
\begin{equation}
\begin{split}
&4[\rho_h^{n+1}-\rho_h^n,\tilde{\rho}_h^{n+1}-\rho_h^{n+1}]= -\frac{8\delta t}{3}[\rho_h^{n+1}-\rho_h^n,\lambda_h^{n+1}+\xi_h^{n+1}-\lambda_h^n-\xi_h^n]\\&=\frac{8\delta t}{3}\big([\rho_h^{n+1},\lambda_h^n]+[\rho_h^n,\lambda_h^{n+1}]\big)\ge 0,
\end{split}
\end{equation}
where have used $\rho_h^{n+1}\lambda_h^{n+1}=\rho_h^{n}\lambda_h^{n}=0$ and $\lambda_h^{n} \ge 0,\; \rho_h^n \ge 0$ for all $n$. Also, notice that 
\begin{equation}
[\rho_h^{n+1}-\rho_h^n, \xi_h^{n+1}-\xi_h^n]=( \xi_h^{n+1}-\xi_h^n)[\rho_h^{n+1}-\rho_h^n,1]=0.
\end{equation}

We calculate the term in \eqref{bdf2:stab:2}
\begin{equation}\label{e:bdf2:f}
6[\tilde{\rho}_h^{n+1}-\rho_h^{n+1},\tilde{\rho}_h^{n+1}]=3(\|(\tilde{\rho}_h^{n+1}\|^2-\|\rho_h^{n+1}\|^2+\|\tilde{\rho}_h^{n+1}-\rho_h^{n+1}\|^2).
\end{equation}

Rewriting \eqref{high:positivity:lag:3} into the following form
\begin{equation}\label{bdf2:stab:20}
3\rho_h^{n+1}-2\delta t(\lambda_h^{n+1}+\xi_h^{n+1})=3\tilde{\rho}_h^{n+1}-2\delta t(\lambda_h^n+\xi_h^n).
\end{equation}
Taking inner product of \eqref{bdf2:stab:20} with itself on both sides, we obtain
\begin{equation}\label{bdf2:stab:30}
\begin{split}
&3\|\rho_h^{n+1}\|^2+\frac 43\delta t^2\|\lambda_h^{n+1}+\xi_h^{n+1}\|^2-4\delta t[\rho_h^{n+1},\xi_h^{n+1}]=3\|\tilde{\rho}_h^{n+1}\|^2\\&+\frac 43\delta t^2\|\lambda_h^{n}+\xi_h^n\|^2-4\delta t[\tilde{\rho}_h^{n+1},\lambda_h^n+\xi_h^n].
\end{split}
\end{equation}
On the left side of \eqref{bdf2:stab:30}, we use KKT condition $\lambda_h^{n+1}\rho_h^{n+1}=0$.

From Lemma \ref{pos:lemma} and notice $\rho_h^{n+1} \ge 0$, we have
\begin{equation}\label{ineq}
-4\delta t[\rho_h^{n+1},\xi_h^{n+1}] \ge 0. 
\end{equation}

Summing up equation \eqref{bdf2:stab:1} and equation \eqref{bdf2:stab:30}, we obtain
\begin{equation}\label{bdf2:stab:4}
\begin{split}
&4\|\rho_h^{n+1}\|^2-4\|\rho_h^n\|^2+\|2\rho_h^{n+1}-\rho_h^n\|^2-\|2\rho_h^{n}-\rho_h^{n-1}\|^2
+4\delta t\|\Grad\tilde{\rho}_h^{n+1}\|^2\\&+\frac 43\delta t^2\|\lambda_h^{n+1}+\xi_h^{n+1}\|^2-\frac 43\delta t^2\|\lambda_h^n+\xi_h^n\|^2-4\delta t[\rho_h^{n+1},\xi_h^{n+1}] =-2\delta t[(\tilde{\rho}_h^{n+1})^2,\Delta  \tilde c_h^{n+1}].
\end{split}
\end{equation}
Applying Young's inequality, for the right hand side of equation \eqref{bdf2:stab:4}, we have
\begin{equation}\label{stab:3:0}
\delta t[(\tilde{\rho}_h^{n+1})^2,\Delta \tilde c_h^{n+1}]\leq \frac{\delta t}{2}(\|(\tilde{\rho}_h^{n+1})^2\|^2+\|\Delta \tilde c_h^{n+1}\|^2).
\end{equation}
We consider the case in 2D, using  Ladyzhenskaya inequality \cite{liu2018positivity}
\begin{equation}\label{La:eq}
\|(\tilde \rho_h^{n+1})^2\|^2 \leq 2\|\tilde \rho_h^{n+1}\|^2\|\Grad \tilde \rho_h^{n+1}\|^2.
\end{equation}

Taking inner product of equation \eqref{high:positivity:lag:1} with $-2\delta t\Delta \tilde c_h^{n+1}$, we obtain
\begin{equation}\label{bdf2:stab:5}
\begin{split}
&\|\Grad \tilde c_h^{n+1}\|^2-\|\Grad c_h^n\|^2 +\|2\Grad \tilde c_h^{n+1}-\Grad c_h^n\|^2-\|2\Grad c_h^{n}-\Grad c_h^{n-1}\|^2
\\&+\|\Grad \tilde c_h^{n+1}-2\Grad c_h^n+\Grad c_h^{n-1}\|^2+2\delta t\|\Delta \tilde c_h^{n+1}\|^2\\& =-2\delta t[2\rho_h^n-\rho_h^{n-1},\Delta \tilde c_h^{n+1} ].
\end{split}
\end{equation}
Using Young inequality, we derive
\begin{equation}\label{yang:eq2}
2\delta t[2\rho_h^n-\rho_h^{n-1},\Delta \tilde c_h^{n+1} ] \leq \delta t(\|\Delta \tilde c_h^{n+1}\|^2+4\|\rho_h^n\|^2+\|\rho_h^{n-1}\|^2).
\end{equation}

Noticing  that $\eta_h^{n+1}$ is a scalar variable and  rewriting \eqref{high:positivity:lag:5} into 
\begin{equation}
c_h^{n+1}-\tilde c_h^{n+1}=\frac{2\delta t}{3}\eta_h^{n+1},
\end{equation}
which implies that $\Grad c_h^{n+1}=\Grad \tilde c_h^{n+1}$. 

Summing up equations \eqref{bdf2:stab:4} and \eqref{bdf2:stab:5} and using \eqref{stab:3:0}, \eqref{ineq}, \eqref{La:eq} and \eqref{yang:eq2}, we obtain
\begin{equation}\label{bdf2:stab:6}
\begin{split}
&4\|\rho_h^{n+1}\|^2-4\|\rho_h^n\|^2+\|2\rho_h^{n+1}-\rho_h^n\|^2-\|2\rho_h^{n}-\rho_h^{n-1}\|^2
+2\delta t(2-\|\tilde{\rho}^{n+1}\|^2)\|\Grad\tilde{\rho}_h^{n+1}\|^2\\&+\frac 43\delta t^2\|\lambda_h^{n+1}+\xi_h^{n+1}\|^2-\frac 43\delta t^2\|\lambda_h^n+\xi_h^n\|^2+\|\Grad  c_h^{n+1}\|^2-\|\Grad c_h^n\|^2 +\|2\Grad  c_h^{n+1}-\Grad c_h^n\|^2\\&-\|2\Grad c_h^{n}-\Grad c_h^{n-1}\|^2
+\|\Grad  c_h^{n+1}-2\Grad c_h^n+\Grad c_h^{n-1}\|^2\leq 
 \delta t(4\|\rho_h^n\|^2+\|\rho_h^{n-1}\|^2).
\end{split}
\end{equation}

Thus, if 
\begin{equation}\label{ass:bdf2}
2-\|\tilde{\rho}_h^{n+1}\|^2>0
\end{equation}
is satisfied, then we obtain
\begin{equation}\label{e:stab:7}
\begin{split}
&4\|\rho_h^{n+1}\|^2+\|2\rho_h^{n+1}-\rho_h^n\|^2+
\frac 43\delta t^2\|\lambda_h^{n+1}+\xi_h^{n+1}\|^2+\|\Grad c_h^{n+1}\|^2 \\&+\|2\Grad c_h^{n+1}-\Grad c_h^n\|^2\leq \|2\rho_h^{n}-\rho_h^{n-1}\|^2+4 (1+\delta t)\|\rho_h^n\|^2 +\delta t\|\rho_h^{n-1}\|^2+ \frac 43\delta t^2\|\lambda_h^n+\xi_h^n\|^2\\&+\|\Grad c_h^n\|^2+\|2\Grad c_h^{n}-\Grad c_h^{n-1}\|^2.
\end{split}
\end{equation}
Summing up \eqref{e:stab:7} from $n=1$ to $n=m-1$, we obtain
\begin{equation}
\begin{split}
&4\|\rho_h^{m}\|^2+\|2\rho_h^{m}-\rho_h^{m-1}\|^2+\frac 43\delta t^2\|\lambda_h^{m}+\xi_h^{m}\|^2
+\|\Grad c_h^{m}\|^2+\|2\Grad c_h^{m}-\Grad c_h^{m-1}\|^2
\\& \leq \delta t\sum\limits_{n=1}^{m-1}(4\|\rho_h^n\|^2+\|\rho_h^{n-1}\|^2)+4\|\rho_h^{1}\|^2+\|2\rho_h^{1}-\rho_h^{0}\|^2+\|\Grad c_h^{1}\|^2 +\|2\Grad c_h^{1}-\Grad c_h^0\|^2.
\end{split}
\end{equation}

Then by Gronwall Lemma \ref{Gron2}, if $m\delta t \leq T$, we derive
\begin{equation}\label{e:stab:8}
\begin{split}
&4\|\rho_h^{m}\|^2+\|2\rho_h^{m}-\rho_h^{m-1}\|^2+\frac 43\delta t^2\|\lambda_h^{m}+\xi_h^m\|^2+\|\Grad c_h^{m}\|^2+\|2\Grad c_h^{m}-\Grad c_h^{m-1}\|^2\\&\leq  e^{5T}(4\|\rho_h^{1}\|^2+\|2\rho_h^{1}-\rho_h^{0}\|^2+\|\Grad c_h^{1}\|^2 +\|2\Grad c_h^{1}-\Grad c_h^0\|^2)\leq e^{5T}C_1,
\end{split}
\end{equation}
where $C_1$ is   constant defined by 
\begin{equation}
C_1=4\|\rho_h^{1}\|^2+\|2\rho_h^{1}-\rho_h^{0}\|^2+\|\Grad c_h^{1}\|^2 +\|2\Grad c_h^{1}-\Grad c_h^0\|^2.
\end{equation}
We rewrite \eqref{high:positivity:lag:3} into 
\begin{equation}\label{eq:bdf}
\tilde \rho_h^{n+1}=\rho_h^{n+1}-\frac{2\delta t}{3}(\lambda_h^{n+1}+\xi_h^{n+1}-\lambda_h^n-\xi^n_h).
\end{equation}
From \eqref{eq:bdf}, we obtain 
\begin{equation}\label{eq:bdf:2}
\|\tilde \rho_h^{n+1}\| \leq \|\rho_h^{n+1}\| + \frac{2\delta t}{3}(\|\lambda_h^{n+1}+\xi_h^{n+1}\|+\|\lambda_h^n+\xi^n_h\|).
\end{equation}
Then using Triangular inequality and Young's inequality, from \eqref{eq:bdf:2}  we derive
\begin{equation}
\|\tilde \rho_h^{n+1}\|^2 \leq 3 \|\rho_h^{n+1}\|^2+\frac{4\delta t^2}{3}\|\lambda_h^{n+1}+\xi_h^{n+1}\|^2
+\frac{4\delta t^2}{3}\|\lambda_h^{n}+\xi_h^{n}\|^2\leq 2e^{5T}C_1.
\end{equation}
Then if we have
\begin{equation}
4\|\rho_h^{1}\|^2+\|2\rho_h^{1}-\rho_h^{0}\|^2+\|\Grad c_h^{1}\|^2 +\|2\Grad c_h^{1}-\Grad c_h^0\|^2
\leq e^{-5T},
\end{equation}
then we can obtain \eqref{ass:bdf2}. Then the proof is finished.

\end{proof}

\subsection{Stability of Crank-Nicolson scheme}

\begin{theorem}\label{cn:stab:result}
Given  $\eps=\gamma=\chi=\mu=1$ and the final time $T$, then for $m\delta t\leq T$, the second-order, energy dissipative, mass conserving and  positivity preserving scheme \eqref{en:positivity:lag:1}-\eqref{en:positivity:lag:7} for Keller-Segel equation \eqref{de:keller:1}-\eqref{de:keller:3}  is stable if the small data condition 
 \begin{equation}
\| \rho_h^{1}\|^2+\|\Grad c_h^{1} \|^2 \leq \frac 23 e^{-\frac{5}{2}T},
\end{equation}
is satisfied.
\end{theorem}
\begin{proof}
Taking inner product of \eqref{en:positivity:lag:1} with $-2\delta t  \Delta \tilde c_h^{n+\frac 12}$ and using Young inequality, we obtain
\begin{equation}\label{p:cn:stab:1}
\begin{split}
&\|\Grad \tilde c_h^{n+1} \|^2 -\|\Grad c_h^{n}\|^2+2\delta t\| \Delta  \tilde c_h^{n+\frac 12}\|^2 = -2\delta t[\frac 32 \rho_h^{n}-\frac 12\rho_h^{n-1}, \Delta \tilde c_h^{n+\frac 12}  ]
\\&\leq \delta t\|\Delta \tilde c_h^{n+\frac 12} \|^2+\delta t(\frac 94\|\rho_h^{n}\|^2+\frac 14\|\rho_h^{n-1}\|^2).
\end{split}
\end{equation}
Taking inner product of \eqref{en:positivity:lag:2} with $2\delta t \tilde \rho_h^{n+\frac 12}$, we obtain
\begin{equation}\label{p:cn:stab:2}
\begin{split}
\|\tilde \rho_h^{n+1}\|^2-\|\rho_h^n\|^2 + 2\delta t\|\Grad \tilde \rho_h^{n+\frac 12}\|^2= 
-\delta t[(\tilde \rho_h^{n+\frac 12})^2,\Delta \tilde c_h^{n+\frac 12}] + 2\delta t [\tilde \rho_h^{n+\frac 12},  \lambda_h^n(\bz)+\xi_h^n].    
\end{split}
\end{equation}
Rewriting \eqref{en:positivity:lag:3} into 
\begin{equation}\label{p:cn:stab:4}
\rho_h^{n+1}(\bz)-\frac{\delta t}{2}(\lambda_h^{n+1}(\bz)+\xi_h^{n+1})=\tilde{\rho}_h^{n+1}(\bz)-\frac{\delta t}{2}(\lambda_h^n(\bz)+\xi_h^n).
\end{equation}
Taking inner product of \eqref{cn:stab:3} with itself, we derive
\begin{equation}\label{p:cn:stab:5}
\begin{split}
&\|\rho_h^{n+1}(\bz)\|^2 + \frac{\delta t^2}{4}\|\lambda_h^{n+1}(\bz)+\xi_h^{n+1}\|^2
-\delta t[\rho_h^{n+1}(\bz),\lambda_h^{n+1}(\bz) +\xi_h^{n+1}]
\\& =\|\tilde{\rho}_h^{n+1}(\bz)\|^2 + \frac{\delta t^2}{4}\|\lambda_h^n(\bz)+\xi_h^n\|^2
-\delta t[\tilde \rho_h^{n+1}(\bz), \lambda_h^n(\bz)+\xi_h^n]. 
\end{split}
\end{equation}
Applying Young's inequality, for the right hand side of equation \eqref{p:cn:stab:2}, we have
\begin{equation}
\delta t[(\tilde{\rho}_h^{n+\frac 12})^2,\Delta \tilde c_h^{n+\frac 12}]\leq \frac{\delta t}{2}(\|(\tilde{\rho}_h^{n+\frac 12})^2\|^2+\|\Delta \tilde c_h^{n+\frac 12}\|^2).
\end{equation}
Using  Ladyzhenskaya inequality, we have
\begin{equation}
\|(\tilde \rho_h^{n+\frac 12})^2\|^2 \leq 2\|\tilde \rho_h^{n+\frac 12}\|^2\|\Grad \tilde \rho_h^{n+\frac 12}\|^2.
\end{equation}
Summing up \eqref{p:cn:stab:2} with \eqref{p:cn:stab:5}, we obtain
\begin{equation}\label{p:cn:stab:6}
\begin{split}
&\| \rho_h^{n+1}\|^2-\|\rho_h^n\|^2 + (2-\|\tilde \rho_h^{n+\frac 12}\|^2)\delta t\|\Grad \tilde \rho_h^{n+\frac 12}\|^2+ \frac{\delta t^2}{4}\|\lambda_h^{n+1}(\bz)+\xi_h^{n+1}\|^2- \frac{\delta t^2}{4}\|\lambda_h^n(\bz)+\xi_h^n\|^2 \\&\leq   \frac{\delta t}{2}\|\Delta \tilde c_h^{n+\frac 12}\|^2+\delta t[\rho_h^{n+1}(\bz),\lambda_h^{n+1}(\bz) +\xi_h^{n+1}]+\delta t[\rho_h^n(\bz), \lambda_h^{n}(\bz) +\xi_h^{n}] \leq \frac{\delta t}{2}\|\Delta \tilde c_h^{n+\frac 12}\|^2,
\end{split}
\end{equation}
where we used KKT condition $\lambda_h^{n+1}\rho_h^{n+1}=0$, $\lambda_h^{n}\rho_h^{n}=0$, and from 
Lemma \ref{pos:lemma}, we have
\begin{equation}
[\rho_h^{n+1}(\bz), \lambda_h^{n+1}(\bz) +\xi_h^{n+1}]+[\rho_h^n(\bz), \lambda_h^{n}(\bz) +\xi_h^{n}]
=[\rho_h^{n+1}(\bz), \xi_h^{n+1}]+[\rho_h^n(\bz), \xi_h^n]\leq 0. 
\end{equation}
Summing up \eqref{p:cn:stab:1} with \eqref{p:cn:stab:6}, we obtain
\begin{equation}\label{p:cn:stab:7}
\begin{split}
&\| \rho_h^{n+1}\|^2-\|\rho_h^n\|^2 + (2-\|\tilde \rho_h^{n+\frac 12}\|^2)\delta t\|\Grad \tilde \rho_h^{n+\frac 12}\|^2+ \frac{\delta t^2}{4}\|\lambda_h^{n+1}(\bz)+\xi_h^{n+1}\|^2- \frac{\delta t^2}{4}\|\lambda_h^n(\bz)+\xi_h^n\|^2 \\&+\|\Grad c_h^{n+1} \|^2 -\|\Grad c_h^{n}\|^2+\frac{\delta t}{2}\| \Delta  \tilde c_h^{n+\frac 12}\|^2\leq \delta t(\frac 94\|\rho_h^{n}\|^2+\frac 14\|\rho_h^{n-1}\|^2).
\end{split}
\end{equation}
where we also  used  $\Grad c_h^{n+1}=\Grad \tilde c_h^{n+1}$ (Notice that from \eqref{en:positivity:lag:6}, we have $\Grad c_h^{n+1}-\Grad \tilde c_h^{n+1}=\frac{\delta t}{\eps}\Grad \eta_h^{n+1}=0$).
Thus, if 
\begin{equation}\label{eq:n}
2-\|\tilde{\rho}_h^{n+\frac 12}\|^2>0,
\end{equation}
is satisfied, then summing up from $n=1$ to $n=m-1$ for  \eqref{p:cn:stab:7}, we have
\begin{equation}\label{p:cn:stab:8}
\begin{split}
\| \rho_h^{m}\|^2+ \frac{\delta t^2}{4}\|\lambda_h^{m}(\bz)+\xi_h^{m+1}\|^2 +\|\Grad c_h^{m} \|^2\leq \delta t\sum\limits_{n=1}^{m-1}(\frac 94\|\rho_h^{n}\|^2+\frac 14\|\rho_h^{n-1}\|^2) +\| \rho_h^{1}\|^2+\|\Grad c_h^{1} \|^2 .
\end{split}
\end{equation}
By Gronwall Lemma \ref{Gron2}, if $m\delta t \leq T$, we have
\begin{equation}
\| \rho_h^{m+1}\|^2+ \frac{\delta t^2}{4}\|\lambda_h^{m+1}(\bz)+\xi_h^{m+1}\|^2 +\|\Grad c_h^{m+1} \|^2 \leq e^{\frac 52 T}(\| \rho_h^{1}\|^2+\|\Grad c_h^{1} \|^2). 
\end{equation}

We rewrite \eqref{en:positivity:lag:3} into 
\begin{equation}\label{p:eq:bdf}
\tilde \rho_h^{n+\frac 12}=\rho_h^{n+1}-\frac{\delta t}{4}(\lambda_h^{n+1}+\xi_h^{n+1}-\lambda_h^n-\xi^n_h).
\end{equation}
From \eqref{p:eq:bdf}, we obtain 
\begin{equation}
\|\tilde \rho_h^{n+\frac 12}\| \leq \|\rho_h^{n+1}\| + \frac{\delta t}{4}(\|\lambda_h^{n+1}+\xi_h^{n+1}\|+\|\lambda_h^n+\xi^n_h\|).
\end{equation}
Then using Triangular inequality and Young's inequality, from \eqref{eq:bdf:2}  we derive
\begin{equation}
\|\tilde \rho_h^{n+\frac 12}\|^2   \leq 3\|\rho_h^{n+1}\|^2+\frac{3\delta t^2}{16}\|\lambda_h^{n+1}+\xi_h^{n+1}\|^2
+\frac{3\delta t^2}{16}\|\lambda_h^{n}+\xi_h^{n}\|^2\leq 3e^{\frac{5}{2}T}(\| \rho_h^{1}\|^2+\|\Grad c_h^{1} \|^2).
\end{equation}
Then if the condition
\begin{equation}
\| \rho_h^{1}\|^2+\|\Grad c_h^{1} \|^2 \leq \frac 23 e^{-\frac{5}{2}T},
\end{equation}
is satisfied, then we can recover \eqref{eq:n}. The proof is completed.

\end{proof}

\subsection{Stability of  bound-preserving schemes}
In this subsection, we shall present stability results for energy dissipative, mass conserving,  bound-preserving schemes \eqref{keller:lag:1}-\eqref{keller:lag:6} for type-II keller-Segel equation \eqref{keller:1}-\eqref{keller:3}.

\begin{theorem}
Setting  $\eps=\mu=\gamma=\chi=1$, then the numerical scheme \eqref{keller:lag:1}-\eqref{keller:lag:6} with $k=2$ is unconditional  stable in the sense that
\begin{equation}
\begin{split}
&4\| \rho_h^{m}\|^2+\|2 \rho_h^{m}- \rho_h^{m-1}\|^2
+\frac 43\delta t^2\|\lambda_h^{m}g'(\rho_h^{m})+\xi_h^m\|^2
+\| \Grad c_h^{m}\|^2+\|2 \Grad c_h^{m}- \Grad c_h^{m-1}\|^2\\&+2\delta t\sum\limits_{n=1}^{m-1}\|\Delta  c_h^{n+1}\|^2
\leq C (4\| \rho_h^{0}\|^2 +\|2 \rho_h^{1}- \rho_h^{0}\|^2+\|\Grad c_h^{0}\|^2+ \|2 \Grad c_h^{1}-\Grad c_h^{0}\|^2),
\end{split}
\end{equation}
for $\forall m \ge 2$ and $m\leq \frac{T}{\delta t}$, where $C$ is a positive constant depending on $M$ and $T$. 
\end{theorem}
\begin{proof}
Taking inner product of equation \eqref{keller:lag:2} with $4\delta t\tilde{\rho}_h^{n+1}$, we obtain
\begin{equation}\label{keller:bdf2:stab:1}
\begin{split}
&[3\tilde{\rho}_h^{n+1}-4 \rho_h^n+ \rho_h^{n-1},2\tilde{ \rho}_h^{n+1}]
+4\delta t \|\Grad \tilde{\rho}_h^{n+1}\|^2=4\delta t[\eta(C_2(\rho_h^{n}))\Grad C_2(c_h^{n}),\Grad \tilde{\rho}_h^{n+1}]
\\&+4\delta t[\lambda_h^n g'(\rho_h^n)+\xi_h^n,\tilde{\rho}_h^{n+1}],
\end{split}
\end{equation}
  where $C_2(\rho_h^{n})=2\rho_h^n-\rho_h^{n-1}$ from \eqref{eq:positivity:bdf2}.

For the first term in \eqref{keller:bdf2:stab:1}, we use \eqref{bdf2:stab:2},  \eqref{bdf2:stab:2:e},  \eqref{bdf2:stab:3} and   \eqref{bdf2:stab:3:e}. 
Consider last term in \eqref{bdf2:stab:3}, using \eqref{keller:lag:3}, then we have
\begin{equation}\label{larger}
\begin{split}
&4[ \rho_h^{n+1}- \rho_h^n,\tilde{ \rho}_h^{n+1}- \rho_h^{n+1}]=-\frac{8\delta t}{3}[\rho_h^{n+1}- \rho_h^n,\lambda_h^{n+1}g'(\rho_h^{n+1})-\lambda_h^ng'(\rho_h^n)+\xi_h^{n+1}-\xi_h^n]\\&
=-\frac{8\delta t}{3}[ \rho_h^{n+1}- \rho_h^n,\lambda_h^{n+1}g'(\rho_h^{n+1})] -\frac{8\delta t}{3}[ \rho_h^{n}- \rho_h^{n+1},\lambda_h^{n}g'(\rho_h^{n})]=I_1 +I_2,
\end{split}
\end{equation}
where we used the mass conservation
\begin{equation}
[\rho_h^{n+1}- \rho_h^n, \xi_h^{n+1}-\xi_h^n]=(\xi_h^{n+1}-\xi_h^n)[\rho_h^{n+1}- \rho_h^n, 1]=0.
\end{equation}
Consider terms  $I_1$ and $I_2$ respectively,  using KKT-condition $\lambda_h^{n+1}g(\rho_h^{n+1})=0$,  we obtain 
\begin{equation}\label{stab:i1}
\begin{split}
I_1&= -\frac{8\delta t}{3}[ \rho_h^{n+1}- \rho_h^n,\lambda_h^{n+1}g'(\rho_h^{n+1})] 
\\&=-\frac{8\delta t}{3}[ \rho_h^{n+1}- \rho_h^n,\lambda_h^{n+1}g'(\rho_h^{n+1})] +\frac{8\delta t}{3}[\lambda_h^{n+1},g(\rho_h^{n+1})]
\\&= -\frac{8\delta t}{3} [\lambda_h^{n+1}, -(\rho_h^{n+1})^2+2\rho_h^n\rho_h^{n+1}-M\rho_h^n ]
\\&=\frac{8\delta t}{3} [ \lambda_h^{n+1}, (\rho_h^{n+1}-\rho_h^n)^2]-\frac{8\delta t}{3} [ \lambda^{n+1}_h, \rho_h^n(\rho_h^n-M)]\ge 0,
\end{split}
\end{equation}
where we use $0\leq \rho_h^n \leq M$ and $\lambda_h^{n+1}\ge 0$.

Similarly, for the term $I_2$,  using KKT-condition $\lambda_h^{n}g(\rho_h^{n})=0$, we have 
\begin{equation}\label{stab:i2}
\begin{split}
I_2&= -\frac{8\delta t}{3}[ \rho_h^{n}- \rho_h^{n+1},\lambda_h^{n}g'(\rho_h^{n})]
\\&=-\frac{8\delta t}{3}[ \rho_h^{n}- \rho_h^{n+1},\lambda_h^{n}g'(\rho_h^{n})] +\frac{8\delta t}{3}[\lambda_h^{n},g(\rho_h^{n})]
\\&= -\frac{8\delta t}{3} [\lambda_h^{n},-(\rho_h^{n})^2+2\rho_h^n\rho_h^{n+1}-M\rho_h^{n+1} ]
\\&=\frac{8\delta t}{3}[ \lambda_h^{n},(\rho_h^{n+1}-\rho_h^n)^2]-[\lambda_h^n, \rho_h^{n+1}(\rho_h^{n+1}-M)]\ge 0,
\end{split}
\end{equation}
since $0\leq \rho_h^{n+1} \leq M$ and $\lambda_h^{n}\ge 0$.

Combining inequalities \eqref{stab:i1} and \eqref{stab:i2}, we obtain
\begin{equation}\label{bdf:in}
\begin{split}
&4[ \rho_h^{n+1}- \rho_h^n,\tilde{ \rho}_h^{n+1}- \rho_h^{n+1}]=-\frac{8\delta t}{3}[ \rho_h^{n+1}- \rho_h^n,\lambda_h^{n+1}g'(\rho_h^{n+1})-\lambda_h^ng'(\rho_h^n)]\ge 0.
\end{split}
\end{equation}

Combing \eqref{e:bdf2:f} with \eqref{bdf:in}, \eqref{bdf2:stab:2},  \eqref{bdf2:stab:2:e},  \eqref{bdf2:stab:3} and   \eqref{bdf2:stab:3:e}, we obtain 
\begin{equation}\label{bdf2:first}
\begin{split}
&(3\|\tilde \rho_h^{n+1}\|^2-3\|\rho_h^{n+1}\|^2)+\|\rho_h^{n+1}\|^2-\|\rho_h^n\|^2+\|2\rho_h^{n+1}-\rho_h^n\|^2-\|2\rho_h^{n}-\rho_h^{n-1}\|^2\\&\leq 
[3\tilde{\rho}_h^{n+1}-4 \rho_h^n+ \rho_h^{n-1},2\tilde{ \rho}_h^{n+1}]. 
\end{split}
\end{equation}

Rewrite \eqref{keller:lag:3} into the following form
\begin{equation}\label{bound:eq:proj}
3 \rho_h^{n+1}-2\delta t(\lambda_h^{n+1}g'(\rho_h^{n+1})+\xi_h^{n+1})=3\tilde{\rho}_h^{n+1}-2\delta t(\lambda_h^ng'(\rho_h^n)+\xi_h^n).
\end{equation}
Taking inner product  of the equation \eqref{bound:eq:proj} with itself on both sides, dividing by $3$ on its both sides, we derive
\begin{equation}\label{bdf2:ineq:first}
\begin{split}
&3\| \rho_h^{n+1}\|^2-4\delta t[\rho_h^{n+1},\lambda_h^{n+1}g'(\rho_h^{n+1})+\xi_h^{n+1}]+\frac 43\delta t^2\|\lambda_h^{n+1}g'(\rho_h^{n+1})+\xi_h^{n+1}\|^2 \\&= 3\|\tilde{\rho}_h^{n+1}\|^2-4\delta t[\tilde{\rho}_h^{n+1},\lambda_h^ng'(\rho_h^n)+\xi_h^n] + \frac 43\delta t^2\|\lambda_h^ng'(\rho_h^n)+\xi_h^n\|^2.
\end{split}
\end{equation}

Taking the discrete inner product of \eqref{re:2:bound} with $1$ on both sides and  noticing $[\rho_h^{n+1},1]=[\rho_h^n,1]$,  we obtain
	\begin{equation}
		[\lambda_h^{n+1}g'(\rho_h^{n+1})+\xi_h^{n+1}, 1]=0,
	\end{equation}
	which implies that 
	\begin{equation}\label{mass:stab:5}
		\xi_h^{n+1}=-\frac{[\lambda_h^{n+1}g'(\rho_h^{n+1}),1]}{|\Omega|}
		=-\frac{[\lambda_h^{n+1},M-2\rho_h^{n+1}]}{|\Omega|}.
	\end{equation}

Using  the fact that $\lambda_h^{n+1}(\bz)g(\rho_h^{n+1}(\bz))=0$ and \eqref{mass:stab:5}, we have
	\begin{equation*}\label{mass:stab:6}
		\begin{split}
			-4\delta t[\rho_h^{n+1},&\lambda_h^{n+1}g'(\rho_h^{n+1})+\xi_h^{n+1}]
			=-4\delta t[\lambda_h^{n+1},\rho_h^{n+1}g'(\rho_h^{n+1})-g(\rho_h^{n+1})]-4\delta t\xi_h^{n+1}[\rho_h^{n+1},1]
			\\&=-4\delta t[\lambda_h^{n+1}, -(\rho_h^{n+1})^2]+\frac{4\delta t}{|\Omega|}[\lambda_h^{n+1},M-2\rho_h^{n+1}][\rho_h^{n+1},1]
			\\&=-4\delta t[\lambda_h^{n+1},-(\rho_h^{n+1}-\frac{[\rho_h^{n+1},1]}{|\Omega|})^2+(\frac{[\rho_h^{n+1},1]}{|\Omega|})(\frac{[\rho_h^{n+1},1]}{|\Omega|}-M)].
		\end{split}
	\end{equation*}
	Since $0\le \rho_h^{n+1} \le M$, we have
	\begin{equation}
		(\frac{[\rho_h^{n+1},1]}{|\Omega|})(\frac{[\rho_h^{n+1},1]}{|\Omega|}-M) \le 0,
	\end{equation}
	which, together with  $\lambda_h^{n+1} \ge 0$, implies that
	\begin{equation}\label{mass:leq}
		-4\delta t[\rho_h^{n+1},\lambda_h^{n+1}g'(\rho_h^{n+1})+\xi_h^{n+1}] \ge 0.
	\end{equation}

Using Young inequality, we derive
\begin{equation}\label{ineq:2}
4\delta t[\eta(C_2(\rho_h^{n}))\Grad C_2(c_h^{n}),\Grad \tilde{\rho}_h^{n+1}]\leq 4\delta t\|\Grad \tilde{\rho}_h^{n+1}\|^2 + 4\delta t M^2\|2\Grad  c_h^{n}-\Grad c_h^{n-1}\|^2, 
\end{equation}
where we use inequality: $-M\leq C_2(\rho_h^{n})=2\rho_h^n-\rho_h^{n-1}\leq 2M$,
\begin{equation}
-2M\leq \eta(C_2(\rho_h^{n}))=\frac{C_2(\rho_h^n)(M-C_2(\rho_h^n))}{M} \leq \frac M4.
\end{equation}

Combining equation \eqref{bdf2:first} with  \eqref{bdf2:ineq:first} and using \eqref{mass:leq}, \eqref{ineq:2}, we derive
\begin{equation}\label{eq:rho:stab}
\begin{split}
&4(\| \rho_h^{n+1}\|^2-\| \rho_h^n\|^2)+\|2 \rho_h^{n+1}- \rho_h^n\|^2-\|2 \rho_h^n- \rho_h^{n-1}\|^2
\\&+\frac 43\delta t^2(\|\lambda_h^{n+1}g'(\rho_h^{n+1})+\xi_h^{n+1}\|^2-\|\lambda_h^ng'(\rho_h^n)+\xi_h^n\|^2)
\leq   4\delta t M^2\|2\Grad  c_h^{n}-\Grad c_h^{n-1}\|^2.
\end{split}
\end{equation}

Taking inner product of equation \eqref{keller:lag:1} with $-4\delta t \Delta \tilde c_h^{n+1}$, we obtain
\begin{equation}\label{eq:c:Stab}
\begin{split}
[3 \Grad \tilde c_h^{n+1}-4\Grad c_h^n+\Grad c_h^{n-1},2\Grad \tilde c_h^{n+1}]+4\delta t \|\Delta \tilde  c_h^{n+1}\|^2=-4\delta t[C_2(\rho_h^n), \Delta \tilde c_h^{n+1}].
\end{split}
\end{equation}

Similar to positivity-preserving schemes \eqref{high:positivity:lag:1}-\eqref{high:positivity:lag:6}, we also have 
$\Grad c_h^{n}=\Grad \tilde c_h^n$ for $\forall n>0$. 

Then,  for the first  term in \eqref{eq:c:Stab}, we have
\begin{equation}\label{non:eq:1}
\begin{split}
&2[3 \Grad c_h^{n+1}-4 \Grad c_h^n+ \Grad c_h^{n-1}, \Grad c_h^{n+1}]=
\| \Grad c_h^{n+1}\|^2-\| \Grad c_h^n\|^2\\&+\|2 \Grad c_h^{n+1}- \Grad c_h^n\|^2-\|2 \Grad  c_h^n- \Grad c_h^{n-1}\|^2+\| \Grad c_h^{n+1}-2 \Grad c_h^n+ \Grad c_h^{n-1}\|^2.
\end{split}
\end{equation}

Using Young inequality, we derive
\begin{equation}\label{non:eq:2}
4\delta t[C_2(\rho_h^n), \Delta \tilde c_h^{n+1}] \leq 2\delta t(\|C_2(\rho_h^n)\|^2+\|\Delta \tilde c_h^{n+1}\|^2).
\end{equation}

Using inequalities \eqref{non:eq:1} and \eqref{non:eq:2}, we rewrite \eqref{eq:c:Stab} into
\begin{equation}\label{stab2:new}
\begin{split}
&\|\Grad c_h^{n+1}\|^2-\|\Grad c_h^n\|^2 + \|2\Grad c_h^{n+1}-\Grad c_h^n\|^2
-\|2\Grad c_h^{n}-\Grad c_h^{n-1}\|^2
\\&+2\delta t\|\Delta  c_h^{n+1}\|^2\leq 2\delta t\|C_2(\rho_h^n)\|^2\leq 2\delta t \|2\rho_h^n-\rho_h^{n-1}\|^2.
\end{split}
\end{equation}

Summing up \eqref{stab2:new} with   \eqref{eq:rho:stab}, we obtain
\begin{equation}\label{stab2:fina:3}
\begin{split}
&4(\| \rho_h^{n+1}\|^2-\| \rho_h^n\|^2)+\|2 \rho_h^{n+1}- \rho_h^n\|^2-\|2 \rho_h^n- \rho_h^{n-1}\|^2
+\frac 43\delta t^2(\|\lambda_h^{n+1}g'(\rho_h^{n+1})+\xi_h^{n+1}\|^2\\&-\|\lambda_h^ng'(\rho_h^n)+\xi_h^n\|^2)
+\|\Grad c_h^{n+1}\|^2-\|\Grad c_h^n\|^2 + \|2\Grad c_h^{n+1}-\Grad c_h^n\|^2
-\|2\Grad c_h^{n}-\Grad c_h^{n-1}\|^2\\&+2\delta t\|\Delta  c_h^{n+1}\|^2
\leq  2\delta t\|2\rho_h^n-\rho_h^{n-1}\|^2+4\delta t M^2\|2\Grad  c_h^{n}-\Grad c_h^{n-1}\|^2.
\end{split}
\end{equation}

Finally summing up \eqref{stab2:fina:3} from $n=1$ to $n=m-1$, then
after dropping some uneccessary terms, we obtain
\begin{equation}\label{stab:final:2}
\begin{split}
&4\| \rho_h^{m}\|^2+\|2 \rho_h^{m}- \rho_h^{m-1}\|^2
+\frac 43\delta t^2\|\lambda_h^{m}g'(\rho_h^{m})+\xi_h^m\|^2
+ \| \Grad c_h^{m}\|^2+\|2 \Grad c_h^{m}- \Grad c_h^{m-1}\|^2\\&+2\delta t\sum\limits_{n=1}^{m-1}\|\Delta  c_h^{n+1}\|^2
\leq 2 \delta t\sum\limits_{n=1}^{m-1}(\|2\rho_h^n-\rho_h^{n-1}\|^2 +2M^2\|2\Grad  c_h^{n}-\Grad c_h^{n-1}\|^2 )+ 
4\| \rho_h^{0}\|^2 +\|2 \rho_h^{1}- \rho_h^{0}\|^2\\&+\|\Grad c_h^{0}\|^2+ \|2\Grad c_h^{1}- \Grad c_h^{0}\|^2.
\end{split}
\end{equation}
Using discrete Gronwall Lemma \eqref{Gron2}, we obtain the stability results.
\end{proof}

\subsection{Stability of Crank-Nicolson scheme} 
Now we consider the stability result for the  energy dissipative and bound-preserving scheme  \eqref{en:keller:lag:1}-\eqref{en:keller:lag:6}. 

\begin{theorem}
Given  $\eps=\mu=\gamma=\chi=1$ and $m\delta t\leq T$, then the energy dissipative, mass conservative and bound-preserving scheme \eqref{en:keller:lag:1}-\eqref{en:keller:lag:6}  is stable in the sense that
\begin{equation}
\begin{split}
&\| \rho_h^{m}\|^2 + \|\Grad c_h^{m} \|^2 + \delta t\sum\limits_{n=0}^{m-1}| \Delta  \tilde c_h^{n+\frac 12}\|^2 + \frac{\delta t^2}{4}\|\lambda_h^{m}(\bz)g'(\rho_h^{m}(\bz))+\xi_h^{m}\|^2+\delta t\sum\limits_{n=0}^{m-1}| \Grad \tilde \rho_h^{n+\frac 12}\|^2\\&\leq C(\|\rho_h^0\|^2 + \|\Grad c_h^{0}\|^2),
\end{split}
\end{equation}
for $\forall m \ge 2$ and $m\leq \frac{T}{\delta t}$, where $C$ is a positive constant depending on $M$ and $T$. 
\end{theorem}
\begin{proof}
In the following proof, notice that $\eps=\mu=\gamma=\chi=1$. 
Taking inner product of \eqref{en:keller:lag:1} with $-2\delta t  \Delta \tilde c_h^{n+\frac 12}$ and using Young inequality, we obtain
\begin{equation}\label{cn:stab:1}
\begin{split}
&\|\Grad \tilde c_h^{n+1} \|^2 -\|\Grad c_h^{n}\|^2+2\delta t\| \Delta  \tilde c_h^{n+\frac 12}\|^2 = -2\delta t(\frac 32 \rho_h^{n}-\frac 12\rho_h^{n-1}, \Delta \tilde c_h^{n+\frac 12}  )
\\&\leq \delta t\|\Delta \tilde c_h^{n+\frac 12} \|^2+\delta t(\frac 94\|\rho_h^{n}\|^2+\frac 14\|\rho_h^{n-1}\|^2).
\end{split}
\end{equation}
Taking inner product of \eqref{en:keller:lag:2} with $2\delta t \tilde \rho_h^{n+\frac 12}$, we obtain
\begin{equation}\label{cn:stab:2}
\begin{split}
\|\tilde \rho_h^{n+1}\|^2-\|\rho_h^n\|^2 + 2\delta t\|\Grad \tilde \rho_h^{n+\frac 12}\|^2&= 2\delta t
(\eta(\frac 32 \rho_h^{n}(\bz)-\frac 12\rho_h^{n-1}(\bz))\Grad (\frac 32c_h^{n}-\frac 12c_h^{n-1}) , \Grad \tilde \rho_h^{n+\frac 12}) \\& + 2\delta t (\tilde \rho_h^{n+\frac 12},  \lambda_h^n(\bz)g'(\rho_h^n(\bz)) +\xi_h^n).    
\end{split}
\end{equation}

Notice that $0\leq \rho_h^{n-1},\rho_h^n \leq M$, we have
\begin{equation}\label{rho:max}
-\frac{3M}{4}\leq \eta(\frac 32\rho_h^{n} -\frac 12\rho_h^{n-1})=\frac{(\frac 32\rho_h^n-\frac 12\rho_h^{n-1})(M-(\frac 32 \rho_h^n-\frac 12\rho_h^{n-1}))}{M} \leq \frac M4, 
\end{equation}
where we used the  inequality: $-\frac{M}{2}\leq \frac 32\rho_h^n-\frac 12\rho_h^{n-1}\leq \frac 32M$.

Using Cauchy-Schwartz inequality, for the first term in the right hand side of \eqref{cn:stab:2}, we have
\begin{equation}\label{non:ineq:en}
\begin{split}
&2\delta t(\eta(\frac 32 \rho_h^{n}(\bz)-\frac 12\rho_h^{n-1}(\bz))\Grad  (\frac 32c_h^{n}-\frac 12c_h^{n-1})(\bz) , \Grad \tilde \rho_h^{n+\frac 12}) \leq 3M\delta t\|\Grad  (\frac 32c_h^{n}-\frac 12c_h^{n-1}) \|\|\Grad \tilde \rho_h^{n+\frac 12}\|
\\& \leq \delta t(\|\Grad \tilde \rho_h^{n+\frac 12}\|^2 +\frac{9M^2}{4}\|\Grad  (\frac 32c_h^{n}-\frac 12c_h^{n-1}) \|^2 ). 
\end{split}
\end{equation}

Using \eqref{non:ineq:en}, equation \eqref{cn:stab:2} can be reformulated as 
\begin{equation}\label{cn:stab:3}
\begin{split}
\|\tilde \rho_h^{n+1}\|^2-\|\rho_h^n\|^2 + \delta t\|\Grad \tilde \rho_h^{n+\frac 12}\|\leq \frac{9M^2\delta t}{4}\|\Grad  (\frac 32c_h^{n}-\frac 12c_h^{n-1}) \|^2  + \delta t (\tilde \rho_h^{n+1}+\rho_h^n,  \lambda_h^n(\bz)g'(\rho_h^n(\bz)) +\xi_h^n).    
\end{split}
\end{equation}

Rewriting \eqref{en:keller:lag:3} into 
\begin{equation}\label{cn:stab:4}
\rho_h^{n+1}(\bz)-\frac{\delta t}{2}(\lambda_h^{n+1}(\bz)g'(\rho_h^{n+1}(\bz))+\xi_h^{n+1})=\tilde{\rho}_h^{n+1}(\bz)-\frac{\delta t}{2}(\lambda_h^n(\bz)g'(\rho_h^n(\bz))+\xi_h^n).
\end{equation}
Taking inner product of \eqref{cn:stab:3} with itself, we derive
\begin{equation}\label{cn:stab:5}
\begin{split}
&\|\rho_h^{n+1}(\bz)\|^2 + \frac{\delta t^2}{4}\|\lambda_h^{n+1}(\bz)g'(\rho_h^{n+1}(\bz))+\xi_h^{n+1}\|^2
-\delta t(\rho_h^{n+1}(\bz),\lambda_h^{n+1}(\bz)g'(\rho_h^{n+1}(\bz)) +\xi_h^{n+1})
\\& =\|\tilde{\rho}_h^{n+1}(\bz)\|^2 + \frac{\delta t^2}{4}\|\lambda_h^n(\bz)g'(\rho_h^n(\bz))+\xi_h^n\|^2
-\delta t(\tilde \rho_h^{n+1}(\bz), \lambda_h^n(\bz)g'(\rho_h^n(\bz))+\xi_h^n). 
\end{split}
\end{equation}
Summing  up \eqref{cn:stab:3} with \eqref{cn:stab:5} and using \eqref{mass:leq}, we obtain
\begin{equation}\label{cn:stab:6}
\begin{split}
&\| \rho_h^{n+1}\|^2-\|\rho_h^n\|^2 + \frac{\delta t^2}{4}\|\lambda_h^{n+1}(\bz)g'(\rho_h^{n+1}(\bz))+\xi_h^{n+1}\|^2-\frac{\delta t^2}{4}\|\lambda_h^n(\bz)g'(\rho_h^n(\bz))+\xi_h^n\|^2\\&+ \delta t\|\Grad \tilde \rho_h^{n+\frac 12}\|\leq \frac{9M^2\delta t}{4}\|\Grad  (\frac 32c_h^{n}-\frac 12c_h^{n-1}) \|^2.   
\end{split}
\end{equation}
Summing  up  \eqref{cn:stab:1} with \eqref{cn:stab:6}, we obtain
\begin{equation}\label{cn:stab:7}
\begin{split}
&\| \rho_h^{n+1}\|^2-\|\rho_h^n\|^2 +\|\Grad \tilde c_h^{n+1} \|^2 -\|\Grad c_h^{n}\|^2+\delta t\| \Delta  \tilde c_h^{n+\frac 12}\|^2\\&+ \frac{\delta t^2}{4}\|\lambda_h^{n+1}(\bz)g'(\rho_h^{n+1}(\bz))+\xi_h^{n+1}\|^2-\frac{\delta t^2}{4}\|\lambda_h^n(\bz)g'(\rho_h^n(\bz))+\xi_h^n\|^2\\&+ \delta t\|\Grad \tilde \rho_h^{n+\frac 12}\|\leq \frac{9M^2\delta t}{4}\|\Grad  (\frac 32c_h^{n}-\frac 12c_h^{n-1})\|^2      +\delta t(\frac 94\|\rho_h^{n}\|^2+\frac 14\|\rho_h^{n-1}\|^2).
\end{split}
\end{equation}
Taking inner product of \eqref{en:keller:lag:5} with $\Delta (c_h^{n+1}+\tilde c_h^{n+1})$ and taking integration by part, we obtain
\begin{equation}\label{eq:grad}
\|\Grad c_h^{n+1}\|^2=\|\Grad \tilde c_h^{n+1}\|^2.
\end{equation}
Using \eqref{eq:grad}, then equation \eqref{cn:stab:7} can be rewritten into 
\begin{equation}\label{cn:stab:8}
\begin{split}
&\| \rho_h^{n+1}\|^2-\|\rho_h^n\|^2 +\|\Grad c_h^{n+1} \|^2 -\|\Grad c_h^{n}\|^2+\delta t\| \Delta  \tilde c_h^{n+\frac 12}\|^2\\&+ \frac{\delta t^2}{4}\|\lambda_h^{n+1}(\bz)g'(\rho_h^{n+1}(\bz))+\xi_h^{n+1}\|^2-\frac{\delta t^2}{4}\|\lambda_h^n(\bz)g'(\rho_h^n(\bz))+\xi_h^n\|^2\\&+ \delta t\|\Grad \tilde \rho_h^{n+\frac 12}\|\leq \frac{9M^2\delta t}{4}(\frac 94\|\Grad c_h^{n} \|^2  + \frac 14 \|\Grad c_h^{n-1} \|^2) +\delta t(\frac 94\|\rho_h^{n}\|^2+\frac 14\|\rho_h^{n-1}\|^2).
\end{split}
\end{equation}
Summing up equation \eqref{cn:stab:8} from $n=0$ to $m-1$, we obtain
\begin{equation}\label{cn:stab:9}
\begin{split}
&\| \rho_h^{m}\|^2 + \|\Grad c_h^{m} \|^2 + \delta t\sum\limits_{n=0}^{m-1}| \Delta  \tilde c_h^{n+\frac 12}\|^2 + \frac{\delta t^2}{4}\|\lambda_h^{m}(\bz)g'(\rho_h^{m}(\bz))+\xi_h^{m}\|^2+\delta t\sum\limits_{n=0}^{m-1}| \Grad \tilde \rho_h^{n+\frac 12}\|^2 \\&
 \leq \|\rho_h^0\|^2 + \|\Grad c_h^{0}\|^2+\delta t\sum\limits_{n=0}^{m-1}\{\frac{81M^2}{16}\|\Grad c_h^{n} \|^2  + \frac{9M^2}{16}\|\Grad c_h^{n-1} \|^2+\frac 94\|\rho_h^{n}\|^2+\frac 14\|\rho_h^{n-1}\|^2\}.
\end{split}
\end{equation}
Notice that $\lambda_h^0=\xi_h^0=0$. Then apply the discrete Gronwall Lemma \ref{Gron2}  to \eqref{cn:stab:9}, we obtain the desired result.

\end{proof}

\section{Error estimate}
We shall carry out a  complete error analysis for a second-order bound preserving scheme with a hybrid spectral discretization that we shall describe below.

  Let $L_N$ be the Legendre polynomial of degree $N$, and $\{x_k\}_{0\le k\le N}$ be the roots of $(1-x^2)L_N'(x)$, i.e., the   Legendre-Gauss-Lobatto points.
We set $\Sigma_N=\{x_k\}_{1\le k\le N-1}$ and  $\bar \Sigma_N=\{x_k\}_{0\le k\le N}$  if $d=1$,  $\Sigma_N=\{(x_k,x_i)\}_{1\le k,i\le N-1}$  and $\bar \Sigma_N=\{(x_k,x_i)\}_{0\le k,i\le N}$ if $d=2$ and   $\Sigma_N=\{(x_k,x_i,x_j)\}_{1\le k,i,j\le N-1}$ and $\bar \Sigma_N=\{(x_k,x_i,x_j)\}_{0\le k,i,j\le N}$ if $d=3$.   We introduce the notation $f \lesssim g$  which implies that $f\leq c g$ where $c$ is a independent constant.

We define the interpolation operator $\hat I_N: C(\Omega)\rightarrow P_N$ by $(\hat I_{N}u)(\bz)=u(\bz)$ for all $\bz\in \bar\Sigma_N$. Then, we  have  \cite{shen2011spectral}
\begin{equation}\label{inter}
	\|v-\hat I_N v\|_{H^s} \lesssim   N^{s-r}\|v\|_{H^r}, \quad \forall v\in H^r(\Omega)\cap X,\, (s=0,\, 1).
\end{equation}
Furthermore, we have the following nice property 
\begin{equation}\label{pro2}
		0\leq \hat I_N u(\bz) \leq M, \quad  \forall  \bz \in \Sigma_N,  \text{ if } 0\le u\le M,
	\end{equation}
	where $\hat I_N$ is the interpolation operator defined above. 
	
It is easy to check that $\hat I_N \rho$ can't preserve the discrete mass, i.e. $[\hat I_N \rho(t^n) ,1]\neq [\rho^n_N,1]$ in general. For the purpose of error analysis, we also introduce the following  "biased" error functions, see also in \cite{tong2024positivity}. We define
\begin{equation}
I_N \rho(t^n)=(1+\eps_N^n)\hat{I}_N \rho(t^n),\quad I_N c(t^n) = \hat I_N c(t^n),
\end{equation}
with 
\begin{equation}
\eps_N^n=\frac{M_0-[\hat I_N\rho(t^n),1]}{[\hat I_N\rho(t^n),1]},\quad M_0=[\hat I_N\rho(t^0),1].
\end{equation}
For initialization, we can obtain  $[\rho_N^0,1]=[\hat I_N\rho(t^0),1]$ by setting $\rho_N^0(\bz)=\frac{M_0}{[\hat I_N\rho(t^0,\bz),1]}\rho(t^0,\bz)$ for  $\forall \bz \in \Sigma_N$. We can easily check that $[I_N \rho(t^n),1]=M_0$ for $0\leq n\leq \frac{T}{\delta t}$. Assume that the density $\rho$ holds the same  regularity as in Theorem \ref{err:th}, we have $|\eps_N^n| \lesssim N^{-l}$ ($l\ge 2$) for $0\leq n\leq \frac{T}{\delta t}$, see \cite{shen2011spectral}.

Denoting that
\begin{equation} \label{error denoting}
\begin{split}
	&	\ole_{\rho, N}^{n+1} = \rho(t^{n+1}) - I_N \rho(t^{n+1}) , \, \he_{\rho,N}^{n+1} = I_N \rho(t^{n+1}) - \rho_N^{n+1} , \, \te_{\rho,N}^{n+1} = I_N \rho(t^{n+1}) - \tilde\rho_N^{n+1},\\
	&\ole_{c,N}^{n+1} = c(t^{n+1}) - I_N c(t^{n+1}) , \, \he_{c,N}^{n+1} = I_N c(t^{n+1}) - c_N^{n+1} , \, \te_{c,N}^{n+1} = I_N c(t^{n+1}) - \tilde c_N^{n+1}.
		\end{split}
	\end{equation}
By introducing  the  "biased" error functions for density $\rho$, using the notations in \eqref{error denoting} for the numerical scheme \eqref{en:keller:lag:1}-\eqref{en:keller:lag:6} we can easily derive that 
\begin{equation}\label{e:con}
[\ole_{\rho,N}^{n+1},1]=[\he_{\rho,N}^{n+1},1]=0.	
\end{equation}
We also derive that from \eqref{error denoting}
\begin{equation}
\begin{split}
	&\rho(t^{n+1}) - \rho_N^{n+1}=\bar e_{\rho,N}^{n+1}+\hat e_{\rho,N}^{n+1},\;\quad \rho(t^{n+1})-\tilde{\rho}_N^{n+1}=\bar e_{\rho,N}^{n+1}+\tilde e_{\rho,N}^{n+1},\\
	& c(t^{n+1}) - c_N^{n+1}=\bar e_{c,N}^{n+1}+\hat e_{c,N}^{n+1},\; \quad c(t^{n+1})-\tilde{c}_N^{n+1}=\bar e_{c,N}^{n+1}+\tilde e_{c,N}^{n+1}.
	\end{split}
\end{equation}
Let $t^k=k\delta t$ and $k=0,1,2,\cdots, \frac{T}{\delta t}$, $t^{k+\frac12}=\frac 12(t^{k+1}+t^k)$ and $f^{n+\frac 12}=\frac{f^{n+1}+f^n}{2}$ for any function $f$. We denote
\begin{equation}
\begin{split}
&e_{\lambda,N}^{n+1} = \lambda(t^{n+1})g'(\rho(t^{n+1})) - \lambda_N^{n+1}g'(\rho_N^{n+1});\quad 
M_N^{n+\frac 12} = \frac 32 \ole^{n}_{\rho,N}-\frac 12 \ole^{n-1}_{\rho,N};\\
& J_N^{n+\frac 12}= \partial_t c(t^{n+\frac 12}) - \frac{c(t^{n+1})-c(t^n)}{\delta t};\;\quad 
I_N^{n+\frac 12}=\rho(t^{n+\frac 12})-(\frac 32\rho(t^n)-\frac 12\rho(t^{n-1}));\\
&K_N^{ n+\frac 12}  = \frac{\ole_{c,N}^{n+1}-\ole_{c,N}^{n}}{\delta t}; \qquad 
T_N^{n+\frac 12} = \Delta\big(\rho(t^{n+\frac 12})  -\frac{\rho (t^{n+1})+\rho(t^n)}{2}\big);\\
&G_N^{n+\frac 12} = \Delta\big(c(t^{n+\frac 12})  -\frac{c (t^{n+1})+c(t^n)}{2}\big);\qquad
H_N^{ n+\frac 12}  = \frac{\Delta(\ole_{c,N}^{n+1}+\ole_{c,N}^{n})}{2}; \\
&B_N^{n+\frac 12}=\partial_t\rho(t^{n+\frac 12}) - \frac{\rho(t^{n+1})-\rho(t^n)}{\delta t},\quad A_N^{ n+\frac 12}  = \frac{\ole_{\rho,N}^{n+1}-\ole_{\rho,N}^{n}}{\delta t}.
\end{split}
\end{equation}
By using the Taylor expansion with integral residual, it is easy to show that \cite{shen2010numerical}
\begin{equation}\label{estimates}
\begin{split}
&\|J_N^{n+\frac 12}\|^2 \lesssim \delta t^3\int_{t^n}^{t^{n+1}}\|c_{ttt}\|^2dt,\quad \|G_N^{n+\frac 12}\|^2 \lesssim \delta t^3\int_{t^n}^{t^{n+1}}\|c_{tt}\|_{H^2}^2dt,\\
&\|T_N^{n+\frac 12}\|_{H^{-1}}^2 \lesssim \delta t^3\int_{t^n}^{t^{n+1}}\|\rho_{tt}\|_{H^1}^2dt,\quad 
\|I_N^{n+\frac 12}\|^2 \lesssim \delta t^3\int_{t^n}^{t^{n+1}}\|\rho_{tt}\|^2dt,\\
&\|B_N^{n+\frac 12}\|_{H^s}^2 \lesssim \delta t^3\int_{t^n}^{t^{n+1}}\|\rho_{ttt}\|_{H^s}^2dt,\quad s=-1,0,1,2,\\
&\|A_N^{n+\frac 12}\|^2\lesssim \frac{1}{\delta t}\int_{t^n}^{t^{n+1}}\|(I-I_N)\rho_t(t)\|^2dt,\quad 
\|K_N^{n+\frac 12}\|^2\lesssim \frac{1}{\delta t}\int_{t^n}^{t^{n+1}}\|(I-I_N)c_t(t)\|^2dt.
\end{split}
\end{equation}

\begin{theorem}\label{err:th}
 Given   $\eps=\mu=\gamma=\chi=1$,  $T\ge 0$ and  $l\ge 2$, for the Keller-Segel equation \eqref{full:dis:2} we assume that  the exact solution $\rho(\bx,t)$ has regularity of 
	\begin{equation*}
	\begin{split}
	   L^2([0,T],H^l(\Omega)) \cap C^1([0,T],H^l(\Omega))\cap C^2([0,T],H^1(\Omega))\cap L^{\infty}([0,T],L^{\infty}(\Omega))\cap C^3([0,T],H^{-1}(\Omega)), 
	  \end{split}
	  \end{equation*}
	  and the exact solution $c(\bx,t)$ has the regularity of 
	  	\begin{equation*}
	   L^2([0,T],H^l(\Omega)) \cap C^1([0,T],H^l(\Omega))\cap C^2([0,T],H^2(\Omega))\cap L^{\infty}([0,T],W^{1,\infty}(\Omega))\cap C^3([0,T],L^2(\Omega)).
	  \end{equation*}
	  Then the numerical scheme \eqref{en:keller:lag:1}-\eqref{en:keller:lag:6} admits the following error estimate:
\begin{equation}\label{error:re}
\begin{split}
&\|\rho(t^m)-\rho_N^m\|^2+ \frac{\delta t^2}{4}\|\lambda_N^{m}g'(\rho_N^{m})+\xi_N^{m}\|^2
+\|\Grad(c(t^m)-c_N^m) \|^2\\&+\frac{\delta t}{4}\sum\limits_{n=0}^{m-1}\|\Delta (c(t^{n+1})-c_N^{n+1})+\Delta (c(t^{n})-c_N^{n})\|^2\lesssim \delta t^4+N^{4-2l} , \quad \forall 0\leq m \leq \frac{T}{\delta t}.
\end{split}
\end{equation}	
	\end{theorem}
\begin{proof}
Subtracting the  first  equation  of \eqref{keller:2} from the numerical scheme \eqref{en:keller:lag:1} , we obtain
\begin{equation}\label{keller:error:01}
\begin{split}
J_N^{n+\frac 12}+\frac{\tilde{e}^{n+1}_{c,N}-\he^{n}_{c,N}}{\delta t} + K_N^{ n+\frac 12}  = \Delta \frac{\tilde{e}^{n+1}_{c,N}+\he^{n}_{c,N}}{2} + \frac 32 \he^{n}_{\rho,N}-\frac 12 \he^{n-1}_{\rho,N} + I_N^{n+\frac 12}+G_N^{n+\frac 12}+H_N^{n+\frac 12}+M_N^{n+\frac 12}.
\end{split}
\end{equation}
From \eqref{en:keller:lag:5}, we derive
\begin{equation}\label{keller:error:02}
\frac{\he_{c,N}^{n+1}-\tilde{e}_{c,N}^{n+1}}{\delta t}=\eta_N^{n+ \frac 12}.
\end{equation}
For density $\rho$, we have
\begin{equation}\label{keller:error:03}
\begin{split}
\frac{\tilde{e}_{\rho,N}^{n+1}-\he_{\rho,N}^n}{\delta t} &
-\Delta\frac{\tilde{e}_{\rho,N}^{n+1}+\he_{\rho,N}^{n}}{2} +\mathcal{N}_N^{n+\frac 12}=-\lambda_N^ng'(\rho_N^n)-\xi_N^n-  R_N^{n+\frac 12} ,
\end{split}
\end{equation}
where
\begin{equation}\label{R:de}
R_N^{n+\frac 12}=B_N^{n+\frac 12}-\Delta\frac{\ole_{\rho,N}^{n+1}+\ole_{\rho,N}^{n}}{2} +\frac{\ole_{\rho,N}^{n+1}-\ole_{\rho,N}^n}{\delta t} -T_N^{n+\frac 12}.
\end{equation}
The third term in \eqref{keller:error:03} is 
\begin{equation}\label{N:eq}
\mathcal{N}_N^{n+\frac 12}=\Grad\cdot(\eta(\rho)\Grad c)(t^{n+\frac 12})-\Grad\cdot\big(\eta(\frac 32\rho_N^{n}-\frac 12\rho_N^{n-1})\Grad (\frac 32 c_N^{n}-\frac 12c_N^{n-1})\big)
=K_N^1+K_N^2+K_N^3,
\end{equation}
in which 
\begin{equation}
\begin{split}
&K_N^1 =\Grad\cdot(\eta(\rho(t^{n+\frac 12}))\Grad c(t^{n+\frac 12})-  \Grad\cdot\big(\eta(\frac 32\rho(t^{n})-\frac 12\rho(t^{n-1}))\Grad (\frac 32 c(t^{n})-\frac 12c(t^{n-1}))\big),\\
& K_N^2=\Grad\cdot\big(\eta(\frac 32\rho(t^{n})-\frac 12\rho(t^{n-1}))\Grad (\frac 32 c(t^{n})-\frac 12c(t^{n-1}))\big)\\&-\Grad\cdot\big(\eta(\frac 32I_N\rho(t^{n})-\frac 12I_N\rho(t^{n-1}))\Grad (\frac 32 I_Nc(t^{n})-\frac 12I_Nc(t^{n-1}))\big),\\
& K_N^3=\Grad\cdot\big(\eta(\frac 32I_N\rho(t^{n})-\frac 12I_N\rho(t^{n-1}))\Grad (\frac 32 I_Nc(t^{n})-\frac 12I_Nc(t^{n-1}))\big)\\&-\Grad\cdot\big(\eta(\frac 32\rho_N^n-\frac 12\rho_N^{n-1})\Grad (\frac 32 c_N^n-\frac 12c_N^{n-1})\big). 
\end{split}
\end{equation}
For correction step \eqref{en:keller:lag:3}, we have 
\begin{equation}\label{keller:error:04}
\frac{\he_{\rho,N}^{n+1}-\tilde{e}_{\rho,N}^{n+1}}{\delta t}=-\frac{\lambda_N^{n+1}g'(\rho_N^{n+1})+\xi_N^{n+1}-(\lambda_N^ng'(\rho_N^n)+\xi_N^n)}{2}.
\end{equation}
Taking inner product of \eqref{keller:error:01} with $-\delta t \Delta (\tilde{e}^{n+1}_{c,N}+\he^{n}_{c,N})$, we derive 
\begin{equation}\label{err:c:1}
\begin{split}
&\|\Grad\tilde{e}^{n+1}_{c,N} \|^2-\|\Grad \he^{n}_{c,N}\|^2+\frac{\delta t}{2}\|\Delta (\tilde{e}^{n+1}_{c,N}+\he^{n}_{c,N})\|^2=-\delta t\big(\frac 32 \he^{n}_{\rho,N}-\frac 12 \he^{n-1}_{\rho,N} , \Delta (\tilde{e}^{n+1}_{c,N}+\he^{n}_{c,N}) \big )\\&-\delta t\big(I_N^{n+\frac 12}+G_N^{n+\frac 12}+H_N^{n+\frac 12}+M_N^{n+\frac 12}-J_N^{n+\frac 12}-K_N^{n+\frac 12},\Delta (\tilde{e}^{n+1}_{c,N}+\he^{n}_{c,N}) \big ).
\end{split}
\end{equation}
For the first term in the right hand side of \eqref{err:c:1}, we have
\begin{equation}\label{err:c:2}
-\delta t\big(\frac 32 \he^{n}_{\rho,N}-\frac 12 \he^{n-1}_{\rho,N} , \Delta (\tilde{e}^{n+1}_{c,N}+\he^{n}_{c,N}) \big )\leq\frac{\delta t}{8}\|\Delta (\tilde{e}^{n+1}_{c,N}+\he^{n}_{c,N})\|^2 + 2\delta t\|\frac 32 \he^{n}_{\rho,N}-\frac 12 \he^{n-1}_{\rho,N} \|^2.
\end{equation}
For the second term in the right hand side of \eqref{err:c:1}, we have
\begin{equation}\label{err:c:3}
\begin{split}
&\delta t\big(I_N^{n+\frac 12}+G_N^{n+\frac 12}+H_N^{n+\frac 12}+M_N^{n+\frac 12}-J_N^{n+\frac 12}-K_N^{n+\frac 12},\Delta (\tilde{e}^{n+1}_{c,N}+\he^{n}_{c,N}) \big )
\\&\leq \frac{\delta t}{8}\|\Delta (\tilde{e}^{n+1}_{c,N}+\he^{n}_{c,N})\|^2 +2 \delta t(\|(I_N^{n+\frac 12}\|^2+\|G_N^{n+\frac 12}\|^2+\|H_N^{n+\frac 12}\|^2)
\\& + 2\delta t(\|(M_N^{n+\frac 12}\|^2+\|J_N^{n+\frac 12}\|^2+\|K_N^{n+\frac 12}\|^2).
\end{split}
\end{equation}
From \eqref{err:c:1}, \eqref{err:c:2}, \eqref{err:c:3}, we obtain
\begin{equation}\label{err:c:4}
\begin{split}
&\|\Grad\tilde{e}^{n+1}_{c,N} \|^2-\|\Grad \he^{n}_{c,N}\|^2+\frac{\delta t}{4}\|\Delta (\tilde{e}^{n+1}_{c,N}+\he^{n}_{c,N})\|^2\leq  2\delta t\|\frac 32 \he^{n}_{\rho,N}-\frac 12 \he^{n-1}_{\rho,N} \|^2
\\& +2 \delta t(\|(I_N^{n+\frac 12}\|^2+\|G_N^{n+\frac 12}\|^2+\|H_N^{n+\frac 12}\|^2)
 + 2\delta t(\|(M_N^{n+\frac 12}\|^2+\|J_N^{n+\frac 12}\|^2+\|K_N^{n+\frac 12}\|^2).
\end{split}
\end{equation}
Taking inner product of  \eqref{keller:error:02} with $\Delta (\tilde{e}^{n+1}_{c,N}+\he^{n}_{c,N})$,  we arrive at
\begin{equation}\label{c:eq:D}
\|\Grad\tilde{e}^{n+1}_{c,N} \|^2-\|\Grad \he^{n}_{c,N}\|^2=\delta t(\eta_N^{n+\frac 12},\Delta (\tilde{e}^{n+1}_{c,N}+\he^{n}_{c,N}))=0.
\end{equation}
Then using \eqref{c:eq:D}, we rewrite \eqref{err:c:4} into 
\begin{equation}\label{err:c:5}
\begin{split}
&\|\Grad\he^{n+1}_{c,N} \|^2-\|\Grad \he^{n}_{c,N}\|^2+\frac{\delta t}{4}\|\Delta (\he^{n+1}_{c,N}+\he^{n}_{c,N})\|^2\leq  2\delta t\|\frac 32 \he^{n}_{\rho,N}-\frac 12 \he^{n-1}_{\rho,N} \|^2
\\& + 2\delta t(\|(I_N^{n+\frac 12}\|^2+\|G_N^{n+\frac 12}\|^2+\|H_N^{n+\frac 12}\|^2)
 + 2\delta t(\|(M_N^{n+\frac 12}\|^2+\|J_N^{n+\frac 12}\|^2+\|K_N^{n+\frac 12}\|^2).
\end{split}
\end{equation}
Taking inner product of \eqref{keller:error:03} with $\delta t (\tilde{e}^{n+1}_{\rho,N}+\he^{n}_{\rho,N})$, we obtain
\begin{equation}\label{rho:eq:1}
\begin{split}
&\|\tilde{e}^{n+1}_{\rho,N}\|^2-\|\he^{n}_{\rho,N}\|^2+ \frac{\delta t}{2}\|\Grad(\tilde{e}^{n+1}_{\rho,N}+\he^{n}_{\rho,N})\|^2+ \delta t(\mathcal{N}_N^{n+\frac 12},\tilde{e}^{n+1}_{\rho,N}+\he^{n}_{\rho,N} )\\&= 
-\delta t(\lambda_N^ng'(\rho_N^n)+\xi_N^n, \tilde{e}^{n+1}_{\rho,N}+\he^{n}_{\rho,N})-\delta t(R_N^{n+\frac 12} , \tilde{e}^{n+1}_{\rho,N}+\he^{n}_{\rho,N}).
\end{split}
\end{equation}
Rewrite \eqref{keller:error:04} into
\begin{equation}\label{corre:eq}
\he_{\rho,N}^{n+1}+\frac{\delta t}{2}(\lambda_N^{n+1}g'(\rho_N^{n+1})+\xi_N^{n+1}) =\tilde{e}_{\rho,N}^{n+1}+\frac{\delta t}{2}(\lambda_N^{n}g'(\rho_N^{n})+\xi_N^{n}).
\end{equation}
Taking inner product of \eqref{corre:eq} with itself on both sides implies that
\begin{equation}\label{rho:eq:2}
\begin{split}
&\|\he_{\rho,N}^{n+1}\|^2+\frac{\delta t^2}{4}\|\lambda_N^{n+1}g'(\rho_N^{n+1})+\xi_N^{n+1}\|^2+\delta t(\he_{\rho,N}^{n+1}, \lambda_N^{n+1}g'(\rho_N^{n+1})+\xi_N^{n+1})\\&
=\|\tilde{e}_{\rho,N}^{n+1}\|^2+\frac{\delta t^2}{4}\|\lambda_N^{n}g'(\rho_N^{n})+\xi_N^{n}\|^2
+\delta t(\tilde{e}_{\rho,N}^{n+1},\lambda_N^{n}g'(\rho_N^{n})+\xi_N^{n} ).
\end{split}
\end{equation}
We can show that
\begin{equation}\label{rho:in:1}
\begin{split}
[\he_{\rho,N}^{n+1}, \lambda_N^{n+1}g'(\rho_N^{n+1})+\xi_N^{n+1}] =- [\rho_N^{n+1} -I_N \rho(t^{n+1}) , \lambda_N^{n+1}g'(\rho_N^{n+1})+\xi_N^{n+1}] \ge 0,
\end{split}
\end{equation}
where we used the following inequality
\begin{equation}\label{j2}
\begin{split}
& -[\rho_N^{n+1}- I_N \rho(t^{n+1}),\lambda_N^{n+1}g'(\rho_N^{n+1})]
\\&=-[\rho_N^{n+1}- I_N \rho(t^{n+1}),\lambda_N^{n+1}g'(\rho_N^{n+1})] +[\lambda_N^{n+1},g(\rho_N^{n+1})]
\\&= - [\lambda_N^{n+1},-(\rho_N^{n+1})^2+2I_N\rho(t^{n+1})\rho_N^{n+1}-MI_N \rho (t^{n+1}) ]
\\&= [\lambda_h^{n+1},(I_N\rho(t^{n+1})-\rho_N^{n+1})^2]-2 [ \lambda_N^{n+1},I_N\rho(t^{n+1})(I_N \rho(t^{n+1})-M)]\ge 0.
\end{split}
\end{equation}
and the equality from \eqref{e:con}
\begin{equation}
[\he_{\rho,N}^{n+1}, \xi_N^{n+1}]=\xi_N^{n+1}[\he_{\rho,N}^{n+1}, 1]=0.
\end{equation}
Similarly, we have 
\begin{equation}\label{rho:in:2}
(\lambda_N^ng'(\rho_N^n)+\xi_N^n, \he^{n}_{\rho,N}) \ge  0.
\end{equation}
Summing up \eqref{rho:eq:1} with \eqref{rho:eq:2} and using \eqref{rho:in:1} and \eqref{rho:in:2}, dropping some positive terms,  it arrives at
\begin{equation}\label{rho:eq:3}
\begin{split}
&\|\he^{n+1}_{\rho,N}\|^2-\|\he^{n}_{\rho,N}\|^2+ \frac{\delta t}{2}\|\Grad(\tilde{e}^{n+1}_{\rho,N}+\he^{n}_{\rho,N})\|^2+ \delta t(\mathcal{N}_N^{n+\frac 12},\tilde{e}^{n+1}_{\rho,N}+\he^{n}_{\rho,N} )\\&+ \frac{\delta t^2}{4}\|\lambda_N^{n+1}g'(\rho_N^{n+1})+\xi_N^{n+1}\|^2-\frac{\delta t^2}{4}\|\lambda_N^{n}g'(\rho_N^{n})+\xi_N^{n}\|^2\le  
-\delta t(R_N^{n+\frac 12} , \tilde{e}^{n+1}_{\rho,N}+\he^{n}_{\rho,N}).
\end{split}
\end{equation}
For the last term in \eqref{rho:eq:3}, using Young's inequality we have
\begin{equation}
-\delta t(R_N^{n+\frac 12} , \tilde{e}^{n+1}_{\rho,N}+\he^{n}_{\rho,N})
\leq \frac{\delta t}{8}\|\Grad(\tilde{e}^{n+1}_{\rho,N}+\he^{n}_{\rho,N})\|^2+ 2\delta tC_0\|R_N^{n+\frac 12}\|_{H^{-1}}^2,
\end{equation}
where $C_0$ is a constant independent of $\delta t$. 
For the third term in \eqref{rho:eq:3}, using \eqref{N:eq} the following inequality is available:
\begin{equation}
\delta t(\mathcal{N}_N^{n+\frac 12},\tilde{e}^{n+1}_{\rho,N}+\he^{n}_{\rho,N} ) = \delta t(K_N^1+K_N^2+K_N^3, \tilde{e}^{n+1}_{\rho,N}+\he^{n}_{\rho,N}).
\end{equation}
We first consider 
\begin{equation}\label{error:s}
\begin{split}
& \delta t(K_N^3, \tilde{e}^{n+1}_{\rho,N}+\he^{n}_{\rho,N}) =\delta t \Big(\Grad\cdot\big(\eta(\frac 32I_N\rho(t^{n})-\frac 12I_N\rho(t^{n-1}))\Grad (\frac 32 I_Nc(t^{n})-\frac 12I_Nc(t^{n-1}))\big)\\&-\Grad\cdot\big(\eta(\frac 32\rho_N^n-\frac 12\rho_N^{n-1})\Grad (\frac 32 I_Nc(t^{n})-\frac 12I_Nc(t^{n-1}))\big), \tilde{e}^{n+1}_{\rho,N}+\he^{n}_{\rho,N} \Big)\\&+\delta t\Big(\Grad\cdot\big(\eta(\frac 32\rho_N^n-\frac 12\rho_N^{n-1})\Grad (\frac 32 \hat e^n_{c,N}-\frac 12\hat e^{n-1}_{c,N})\big), \tilde{e}^{n+1}_{\rho,N}+\he^{n}_{\rho,N}\Big),
\end{split}
\end{equation}
where we used 
\begin{equation}
\begin{split}
K_N^3&= \Grad\cdot\big(\eta(\frac 32I_N\rho(t^{n})-\frac 12I_N\rho(t^{n-1}))\Grad (\frac 32 I_Nc(t^{n})-\frac 12I_Nc(t^{n-1}))\big)\\&-\Grad\cdot\big(\eta(\frac 32\rho_N^n-\frac 12\rho_N^{n-1})\Grad (\frac 32 I_Nc(t^{n})-\frac 12I_Nc(t^{n-1}))\big)\\&+\Grad\cdot\big(\eta(\frac 32\rho_N^n-\frac 12\rho_N^{n-1})\Grad (\frac 32 I_Nc(t^{n})-\frac 12I_Nc(t^{n-1}))\big)\\&-\Grad\cdot\big(\eta(\frac 32\rho_N^n-\frac 12\rho_N^{n-1})\Grad (\frac 32 c_N^n-\frac 12c_N^{n-1})\big).
\end{split}
\end{equation}
Notice that $0\leq \rho_N^n,\rho_N^{n-1}\leq M$ and $0\leq  I_N \rho(t^{n}),I_N \rho(t^{n-1}) \leq M $, using $\eta'(\rho)=\frac{M-2\rho}{M}$ then we have
\begin{equation}\label{ineq:n}
|\eta(\frac 32I_N\rho(t^{n})-\frac 12I_N\rho(t^{n-1}))-\eta(\frac 32\rho_N^n-\frac 12\rho_N^{n-1})|
\leq 2|\frac 32\hat e_{\rho,N}^n-\frac 12\hat e_{\rho,N}^{n-1}|.
\end{equation}
Using \eqref{ineq:n}, taking integration by part, we have
\begin{equation}\label{error:7}
\begin{split}
&\Big(\Grad\cdot\big(\eta(\frac 32I_N\rho(t^{n})-\frac 12I_N\rho(t^{n-1}))\Grad (\frac 32 I_Nc(t^{n})-\frac 12I_Nc(t^{n-1}))\big)\\&-\Grad\cdot\big(\eta(\frac 32\rho_N^n-\frac 12\rho_N^{n-1})\Grad (\frac 32 I_Nc(t^{n})-\frac 12I_Nc(t^{n-1}))\big), \tilde{e}^{n+1}_{\rho,N}+\he^{n}_{\rho,N} \Big)\\&\leq 2\delta t\Big(|\frac 32\hat e_{\rho,N}^n-\frac 12\hat e_{\rho,N}^{n-1}||\Grad (\frac 32 I_Nc(t^{n})-\frac 12I_Nc(t^{n-1}))|, |\Grad (\tilde{e}^{n+1}_{\rho,N}+\he^{n}_{\rho,N})| \Big)
\\& \leq  \frac{\delta t}{8}\|\Grad (\tilde{e}^{n+1}_{\rho,N}+\he^{n}_{\rho,N})\|^2 + 8\delta t(C^{\star})^2\|\frac 32\hat e_{\rho,N}^n-\frac 12\hat e_{\rho,N}^{n-1}\|^2,
\end{split}
\end{equation}
where $C^{\star}=\|\Grad (\frac 32 I_Nc(t^{n})-\frac 12I_Nc(t^{n-1}))\|_{L^\infty}$. 
For the last term in \eqref{error:s}, using \eqref{rho:max} we have
\begin{equation}\label{error:8}
\begin{split}
&\delta t\Big(\Grad\cdot\big(\eta(\frac 32\rho_N^n-\frac 12\rho_N^{n-1})\Grad (\frac 32 \hat e^n_{c,N}-\frac 12\hat e^{n-1}_{c,N})\big), \tilde{e}^{n+1}_{\rho,N}+\he^{n}_{\rho,N}\Big)
\\&=-\delta t\Big( \eta(\frac 32\rho_N^n-\frac 12\rho_N^{n-1})\Grad (\frac 32 \hat e^n_{c,N}-\frac 12\hat e^{n-1}_{c,N}) , \Grad (\tilde{e}^{n+1}_{\rho,N}+\he^{n}_{\rho,N})\Big)
\\& \leq \frac{\delta t}{8}\|\Grad (\tilde{e}^{n+1}_{\rho,N}+\he^{n}_{\rho,N})\|^2 + \frac{9M^2\delta t}{8}\|\Grad (\frac 32 \hat e^n_{c,N}-\frac 12\hat e^{n-1}_{c,N})\|^2.
\end{split}
\end{equation}
Combining \eqref{error:7} and \eqref{error:8}, we obtain
\begin{equation}\label{error:9}
\begin{split}
 &\delta t(K_N^3, \tilde{e}^{n+1}_{\rho,N}+\he^{n}_{\rho,N})  \leq \frac{\delta t}{4}\|\Grad (\tilde{e}^{n+1}_{\rho,N}+\he^{n}_{\rho,N})\|^2 + 8\delta t(C^{\star})^2\|\frac 32\hat e_{\rho,N}^n-\frac 12\hat e_{\rho,N}^{n-1}\|^2
\\& + \frac{9M^2\delta t}{8}\|\Grad (\frac 32 \hat e^n_{c,N}-\frac 12\hat e^{n-1}_{c,N})\|^2.
\end{split}
\end{equation}
Using Taylor expansion we derive
\begin{equation}\label{error:10}
\begin{split}
&\delta t(K_N^1, \tilde{e}^{n+1}_{\rho,N}+\he^{n}_{\rho,N}) = -\Big(\eta(\rho(t^{n+\frac 12}))\Grad c(t^{n+\frac 12})-  \eta(\frac 32\rho(t^{n})-\frac 12\rho(t^{n-1}))\Grad (\frac 32 c(t^{n})-\frac 12c(t^{n-1})),\\& \Grad(\tilde{e}^{n+1}_{\rho,N}+\he^{n}_{\rho,N})\Big)
\leq C_1\delta t^5+\frac{\delta t}{8}\|\Grad(\tilde{e}^{n+1}_{\rho,N}+\he^{n}_{\rho,N})\|^2.
\end{split}
\end{equation}
Similar for the following term, we have
\begin{equation}\label{error:11}
\begin{split}
&\delta t(K_N^2, \tilde{e}^{n+1}_{\rho,N}+\he^{n}_{\rho,N}) =-\Big(\eta(\frac 32\rho(t^{n})-\frac 12\rho(t^{n-1}))\Grad (\frac 32 c(t^{n})-\frac 12c(t^{n-1}))\\&-\eta(\frac 32I_N\rho(t^{n})-\frac 12I_N\rho(t^{n-1}))\Grad (\frac 32 I_Nc(t^{n})-\frac 12I_Nc(t^{n-1})), \Grad(\tilde{e}^{n+1}_{\rho,N}+\he^{n}_{\rho,N})\Big)
\\&=\Big(\eta(\frac 32\rho(t^{n})-\frac 12\rho(t^{n-1}))\Grad (\frac 32 c(t^{n})-\frac 12c(t^{n-1}))\\&-\eta(\frac 32I_N\rho(t^{n})-\frac 12I_N\rho(t^{n-1}))\Grad (\frac 32 c(t^{n})-\frac 12c(t^{n-1})), \Grad(\tilde{e}^{n+1}_{\rho,N}+\he^{n}_{\rho,N})\Big)
\\&+\Big(\eta(\frac 32I_N\rho(t^{n})-\frac 12I_N\rho(t^{n-1}))\Grad (\frac 32 c(t^{n})-\frac 12c(t^{n-1}))\\&-\eta(\frac 32I_N\rho(t^{n})-\frac 12I_N\rho(t^{n-1}))\Grad (\frac 32 I_Nc(t^{n})-\frac 12 I_N c(t^{n-1})), \Grad(\tilde{e}^{n+1}_{\rho,N}+\he^{n}_{\rho,N})\Big)
\\& \leq \frac{\delta t}{8}\|\Grad(\tilde{e}^{n+1}_{\rho,N}+\he^{n}_{\rho,N})\|^2+C_2\delta t\|\frac 32\bar{e}_{\rho,N}^n-\frac 12\bar{e}_{\rho,N}^{n-1}\|^2+ C_3\delta t\|\frac 32\bar{e}_{c,N}^n-\frac 12\bar{e}_{c,N}^{n-1}\|^2.
\end{split}
\end{equation}
Combining \eqref{error:9}, \eqref{error:10}, \eqref{error:11} with \eqref{rho:eq:3}, we have
\begin{equation}\label{error:12}
\begin{split}
&\|\he^{n+1}_{\rho,N}\|^2-\|\he^{n}_{\rho,N}\|^2+ \frac{\delta t^2}{4}\|\lambda_N^{n+1}g'(\rho_N^{n+1})+\xi_N^{n+1}\|^2-\frac{\delta t^2}{4}\|\lambda_N^{n}g'(\rho_N^{n})+\xi_N^{n}\|^2\\&\le  
C_1\delta t^5+ C_2\delta t\|\frac 32\bar{e}_{\rho,N}^n-\frac 12\bar{e}_{\rho,N}^{n-1}\|^2+C_3 \delta t\|\frac 32\bar{e}_{c,N}^n-\frac 12\bar{e}_{c,N}^{n-1}\|^2
\\&+ 2C_0\delta t\|R_N^{n+\frac 12}\|_{H^{-1}}^2 +\frac{9M^2\delta t}{8}\|\Grad (\frac 32 \hat e^n_{c,N}-\frac 12\hat e^{n-1}_{c,N})\|^2+8\delta t(C^{\star})^2\|\frac 32\hat e_{\rho,N}^n-\frac 12\hat e_{\rho,N}^{n-1}\|^2.
\end{split}
\end{equation}
Summing up \eqref{error:12} with \eqref{err:c:5}, we obtain
\begin{equation}\label{error:13}
\begin{split}
&\|\he^{n+1}_{\rho,N}\|^2-\|\he^{n}_{\rho,N}\|^2+ \|\Grad\he^{n+1}_{c,N} \|^2-\|\Grad \he^{n}_{c,N}\|^2+\frac{\delta t}{4}\|\Delta (\he^{n+1}_{c,N}+\he^{n}_{c,N})\|^2+\frac{\delta t^2}{4}\|\lambda_N^{n+1}g'(\rho_N^{n+1})\\&+\xi_N^{n+1}\|^2-\frac{\delta t^2}{4}\|\lambda_N^{n}g'(\rho_N^{n})+\xi_N^{n}\|^2\le  
C_1\delta t^5+ C_2\delta t\|\frac 32\bar{e}_{\rho,N}^n-\frac 12\bar{e}_{\rho,N}^{n-1}\|^2+C_3 \delta t\|\frac 32\bar{e}_{c,N}^n-\frac 12\bar{e}_{c,N}^{n-1}\|^2
\\&+ 2C_0\delta t\|R_N^{n+\frac 12}\|_{H^{-1}}^2 +\frac{9M^2\delta t}{8}\|\Grad (\frac 32 \hat e^n_{c,N}-\frac 12\hat e^{n-1}_{c,N})\|^2+2\delta t(4(C^{\star})^2+1)\|\frac 32\hat e_{\rho,N}^n-\frac 12\hat e_{\rho,N}^{n-1}\|^2
\\&+2\delta t(\|(I_N^{n+\frac 12}\|^2+\|G_N^{n+\frac 12}\|^2+\|H_N^{n+\frac 12}\|^2)
 + 2\delta t(\|(M_N^{n+\frac 12}\|^2+\|J_N^{n+\frac 12}\|^2+\|K_N^{n+\frac 12}\|^2).
\end{split}
\end{equation}
Summing up \eqref{error:13} from $n=0$ to $n=m-1=\frac{T}{\delta t}-1$, we derive
\begin{equation}\label{error:14}
\begin{split}
&\|\he^{m}_{\rho,N}\|^2+ \frac{\delta t^2}{4}\|\lambda_N^{m}g'(\rho_N^{m})+\xi_N^{m}\|^2
+\|\Grad\he^{m}_{c,N} \|^2+\frac{\delta t}{4}\sum\limits_{n=0}^{m-1}\|\Delta (\he^{n+1}_{c,N}+\he^{n}_{c,N})\|^2\\&\leq \delta t\sum\limits_{n=0}^{m-1}\Big\{  C_1\delta t^4+ C_2\|\frac 32\bar{e}_{\rho,N}^n-\frac 12\bar{e}_{\rho,N}^{n-1}\|^2+C_3 \|\frac 32\bar{e}_{c,N}^n-\frac 12\bar{e}_{c,N}^{n-1}\|^2
\\&+ 2C_0\|R_N^{n+\frac 12}\|_{H^{-1}}^2 +\frac{9M^2}{8}\|\Grad (\frac 32 \hat e^n_{c,N}-\frac 12\hat e^{n-1}_{c,N})\|^2+2(4(C^{\star})^2+1)\|\frac 32\hat e_{\rho,N}^n-\frac 12\hat e_{\rho,N}^{n-1}\|^2
\\&+2(\|(I_N^{n+\frac 12}\|^2+\|G_N^{n+\frac 12}\|^2+\|H_N^{n+\frac 12}\|^2)
 + 2(\|(M_N^{n+\frac 12}\|^2+\|J_N^{n+\frac 12}\|^2+\|K_N^{n+\frac 12}\|^2).\Big\}.
\end{split}
\end{equation}
Applying the discrete Gronwall Lemma to \eqref{error:14} and noticing the following estimates from \eqref{inter} and \eqref{estimates}
\begin{equation}\label{error:15}
\begin{split}
&\|\bar{e}_{c,N}^n\|^2 \lesssim N^{-2l},\;\quad \|\bar{e}_{\rho,N}^n\|^2\lesssim N^{-2l},\quad \|\Grad \bar{e}_{c,N}^n\|^2 \lesssim N^{2-2l},\;\quad \|\Grad \bar{e}_{\rho,N}^n\|^2\lesssim N^{2-2l},\\
&\|H_N^{n+\frac 12}\|^2 \lesssim N^{4-2l},\;\quad \|M_N^{n+\frac 12}\|^2 \lesssim N^{-2l},\; \quad \|K_N^{n+\frac 12}\|^2\lesssim N^{-2l},\\
&\|J_N^{n+\frac 12}\|^2\lesssim \delta t^4,\;\quad \|I_N^{n+\frac 12}\|^2\lesssim \delta t^4,\;\quad \|G_N^{n+\frac 12}\|^2 \lesssim \delta t^4.
\end{split}
\end{equation}
From \eqref{R:de} and \eqref{estimates}, we also have 
\begin{equation}\label{error:16}
\begin{split}
&\delta t\sum\limits_{n=0}^{m-1}\|R_N^{n+\frac 12}\|_{H^{-1}}^2\leq \delta t\sum\limits_{n=0}^{m-1}\|B_N^{n+\frac 12}\|_{H^{-1}}^2+\delta t\sum\limits_{n=0}^{m-1}\|\Delta\frac{\ole_{\rho,N}^{n+1}+\ole_{\rho,N}^{n}}{2}\|_{H^{-1}}^2 \\&+\delta t\sum\limits_{n=0}^{m-1}\|\frac{\ole_{\rho,N}^{n+1}-\ole_{\rho,N}^n}{\delta t}\|_{H^{-1}}^2 +\delta t\sum\limits_{n=0}^{m-1}\|T_N^{n+\frac 12}\|_{H^{-1}}^2 \lesssim \delta t^4+ N^{4-2l}.
\end{split}
\end{equation}
Finally, combining \eqref{error:14}-\eqref{error:16}, we have
\begin{equation}
\|\he^{m}_{\rho,N}\|^2+ \frac{\delta t^2}{4}\|\lambda_N^{m}g'(\rho_N^{m})+\xi_N^{m}\|^2
+\|\Grad\he^{m}_{c,N} \|^2+\frac{\delta t}{4}\sum\limits_{n=0}^{m-1}\|\Delta (\he^{n+1}_{c,N}+\he^{n}_{c,N})\|^2
\lesssim \delta t^4+ N^{4-2l}.
\end{equation}
Then we obtain the desired error estimate \eqref{error:re}.

\end{proof}

\section{Numerical simulations}
In this section, we carry out  various numerical experiments to demonstrate the performance of proposed positivity schemes  \eqref{en:positivity:lag:1}-\eqref{en:positivity:lag:7}, \eqref{high:positivity:lag:1}-\eqref{high:positivity:lag:6} and bound-preserving schemes \eqref{keller:lag:1}-\eqref{keller:lag:6}, \eqref{en:keller:lag:1}-\eqref{en:keller:lag:6}  for Keller-Segel equations \eqref{de:keller:1}-\eqref{de:keller:3} and \eqref{keller:1}-\eqref{keller:3}. In all numerical simulations, we use periodic boundary condition and implement Fourier-Galerkin methods  in space.

\subsection{Convergence rate}
We first test the convergence rate in time for positivity schemes \eqref{en:positivity:lag:1}-\eqref{en:positivity:lag:7}, \eqref{high:positivity:lag:1}-\eqref{high:positivity:lag:6} and bound-preserving schemes \eqref{keller:lag:1}-\eqref{keller:lag:6} in the domain $[0,2\pi]^2$. The initial conditions  are given by:
\begin{equation}
\begin{split}
&\rho(\bx,0)=10 e^{-10(x-\pi)^2-10(y-\pi)^2},\\
&c(\bx,0)=10 e^{-(x-\pi)^2/2-(y-\pi)^2/2}.
\end{split}
\end{equation} 
We choose $\mu=\gamma=\eps=1$ and $\chi=2$ for two  types of Keller-Segel equations and $N=64$ Fourier collocations in each direction. The reference solution $\rho_e$ is computed by second-order positivity scheme \eqref{high:positivity:lag:1}-\eqref{high:positivity:lag:6} with small time step $\delta t=10^{-6}$.

\begin{figure}[htbp]
\centering
\subfigure[$\rho$ positivity-preserving.]{
\includegraphics[width=0.40\textwidth,clip==]{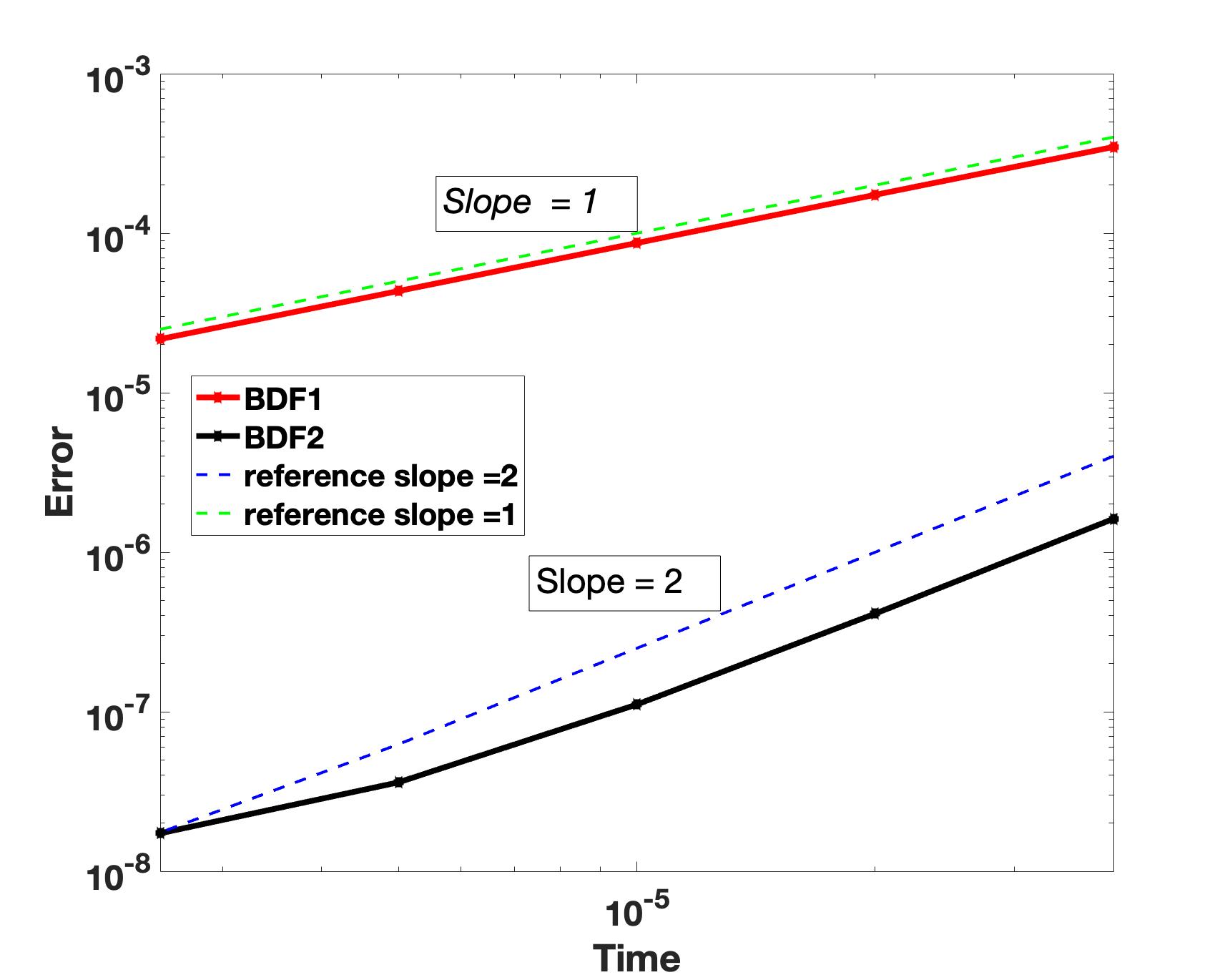}}
\subfigure[$c$ positivity-preserving.]{
\includegraphics[width=0.40\textwidth,clip==]{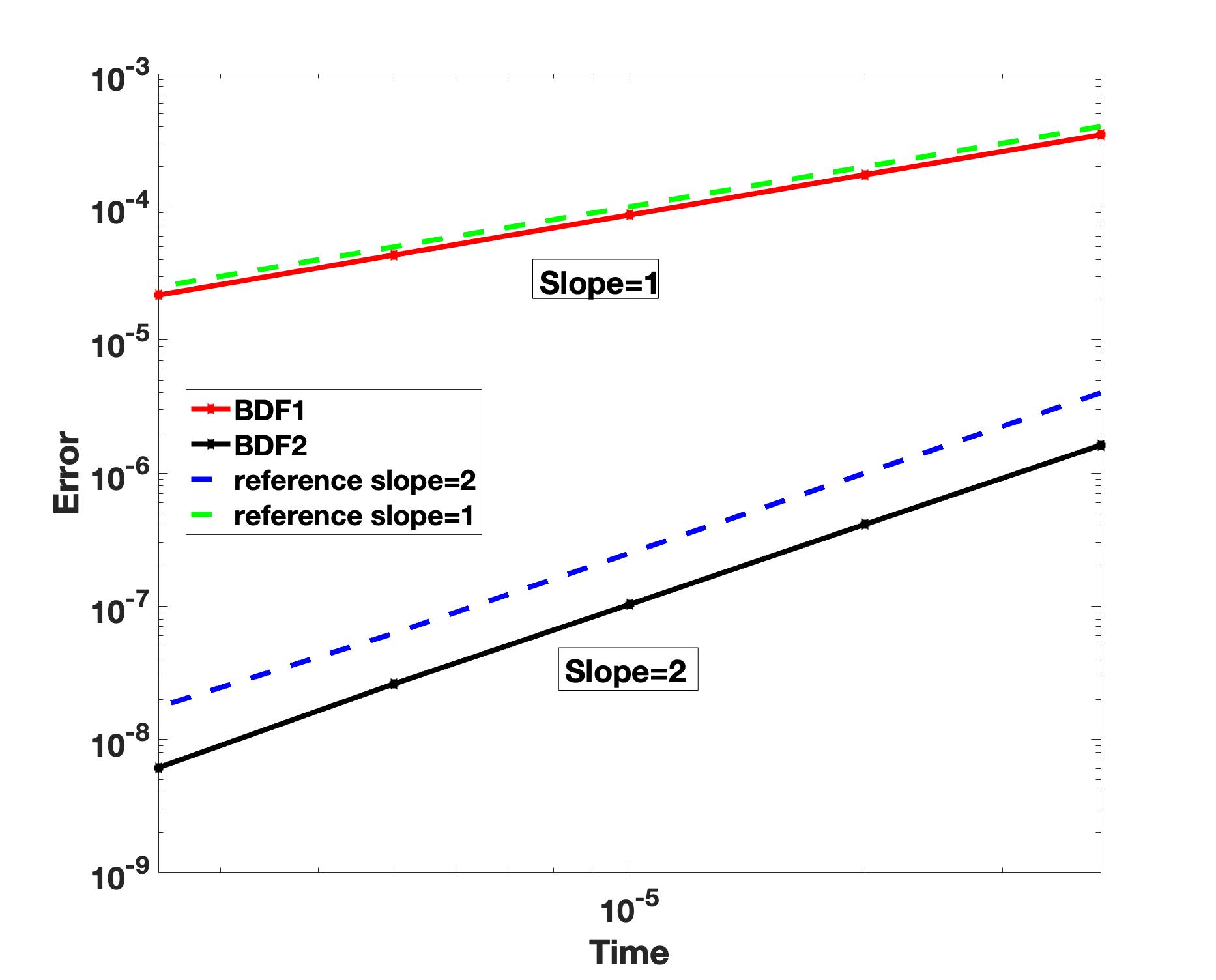}}
\caption{Accuracy test: The $L^{\infty}$ errors  of density  $\rho(\bx,t)$ and concentration $c(\bx,t)$ at $t=0.01$ for the first type Keller-Segel equations \eqref{de:keller:1}-\eqref{de:keller:3}   computed by positivity schemes \eqref{high:positivity:lag:1}-\eqref{high:positivity:lag:4} with $k=1,2$.}\label{convergence}
\end{figure}

\begin{figure}[htbp]
\centering
\subfigure[$\rho$ bound-preserving]{
\includegraphics[width=0.40\textwidth,clip==]{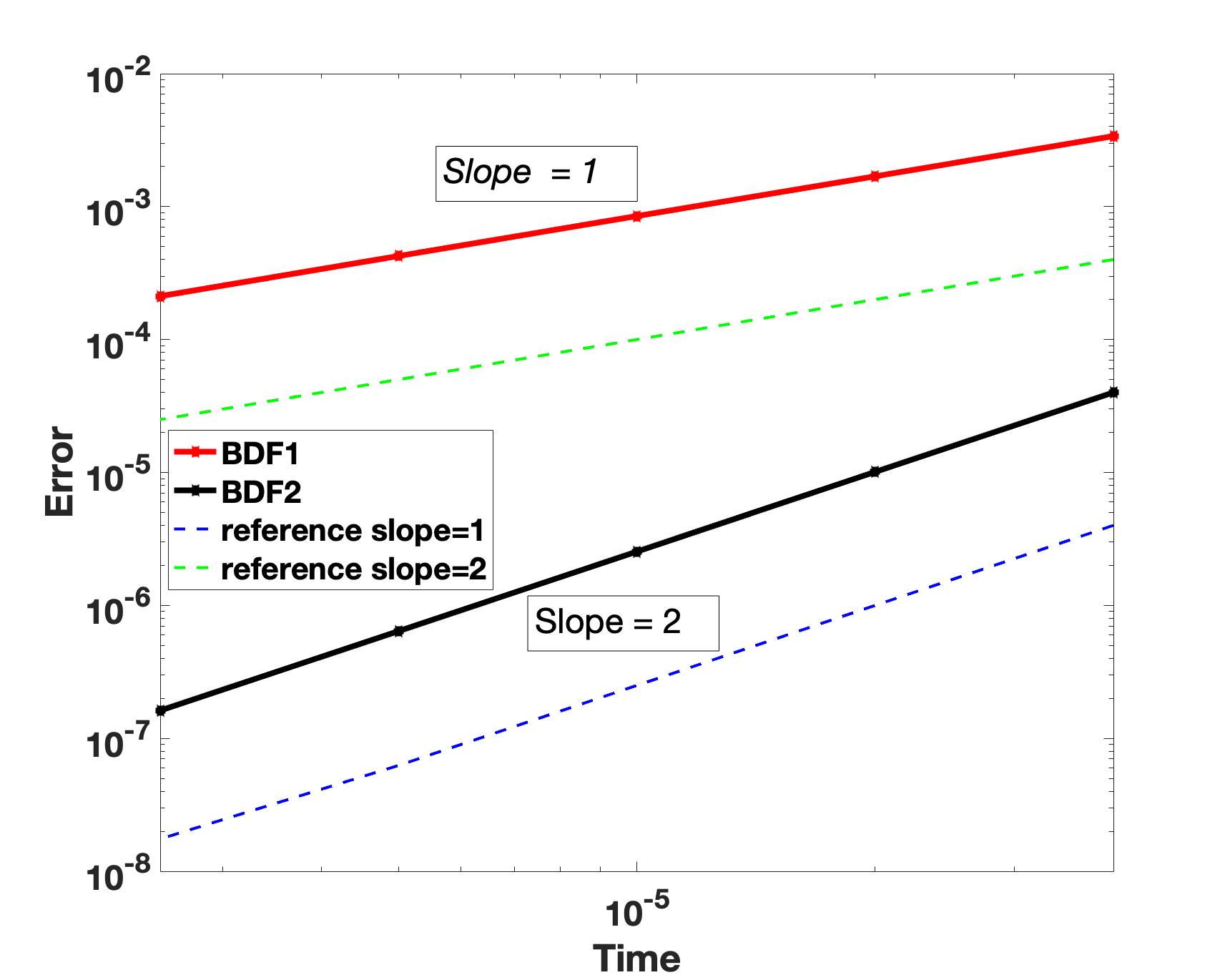}}
\subfigure[ $c$ bound-preserving.]{
\includegraphics[width=0.40\textwidth,clip==]{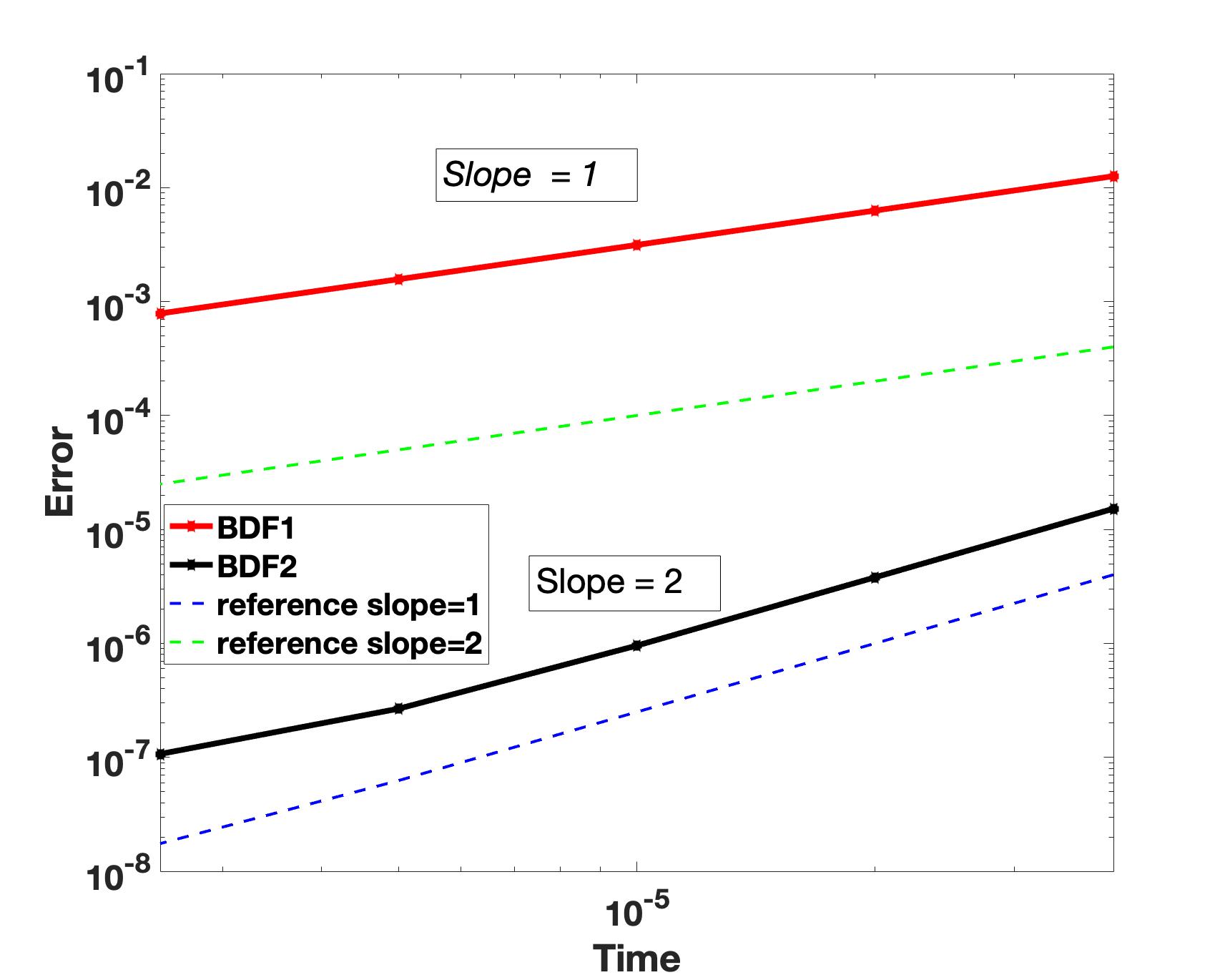}}
\caption{Accuracy test: The $L^{\infty}$ errors  of density  $\rho(\bx,t)$ and concentration $c(\bx,t)$ at $t=0.01$ for the second  type Keller-Segel equations \eqref{keller:1}-\eqref{keller:3}   computed by bound-preserving schemes \eqref{keller:lag:1}-\eqref{keller:lag:6} with $k=1,2$.}\label{convergence2}
\end{figure}

\begin{figure}[htbp]
\centering
\subfigure[$\rho$ at $t=0.05$.]{
\includegraphics[width=0.4\textwidth,clip==]{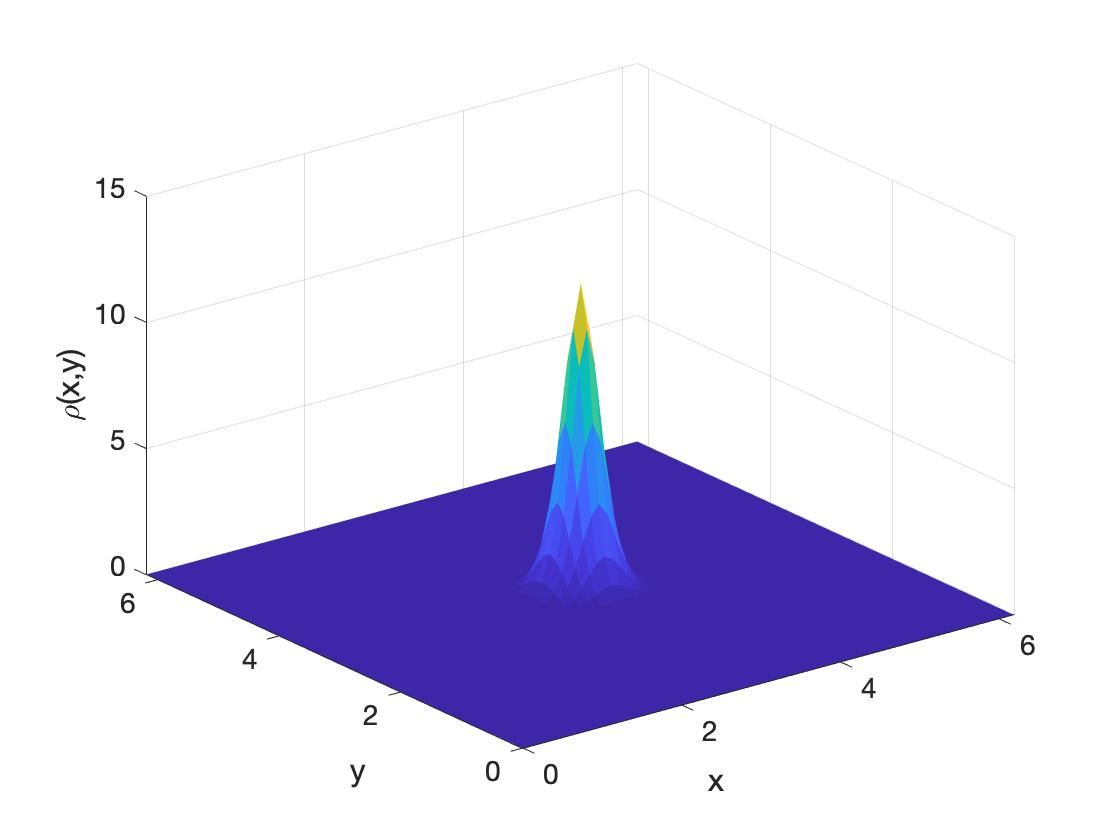}}
\subfigure[$\rho$ at $t=0.2$.]{
\includegraphics[width=0.4\textwidth,clip==]{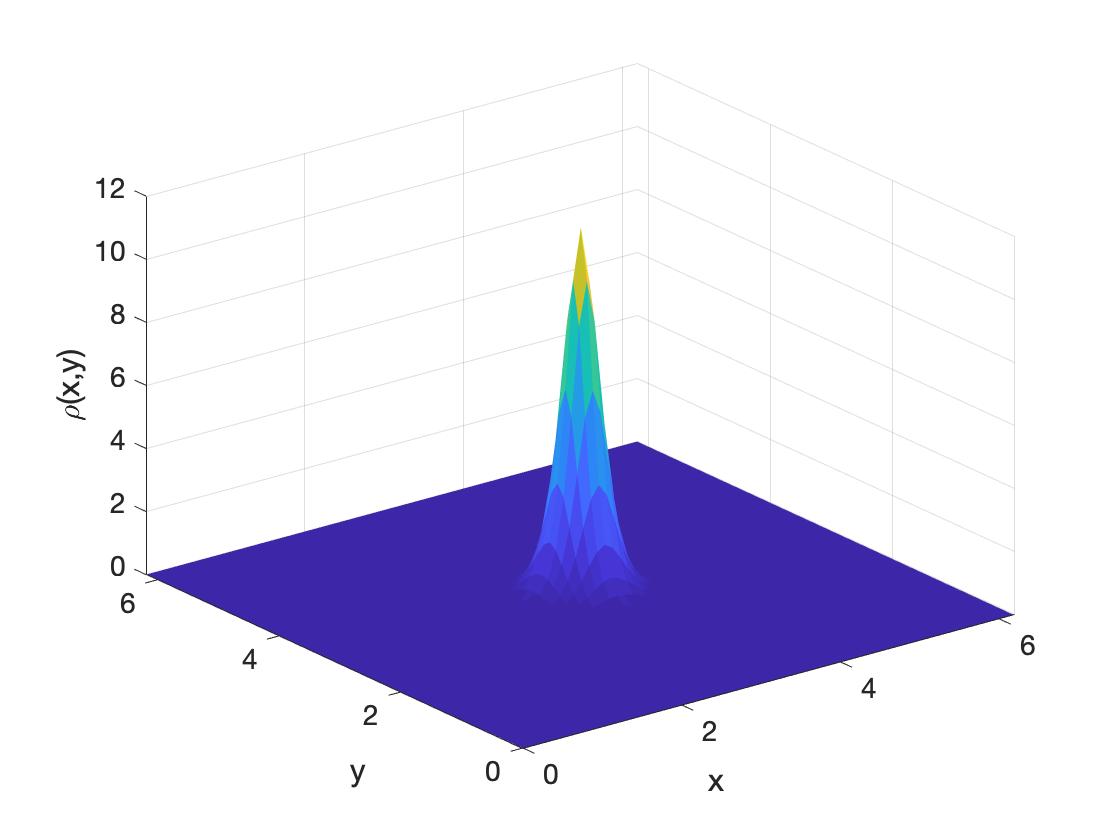}}
\subfigure[$\lambda$ at $t=0.05$.]{
\includegraphics[width=0.4\textwidth,clip==]{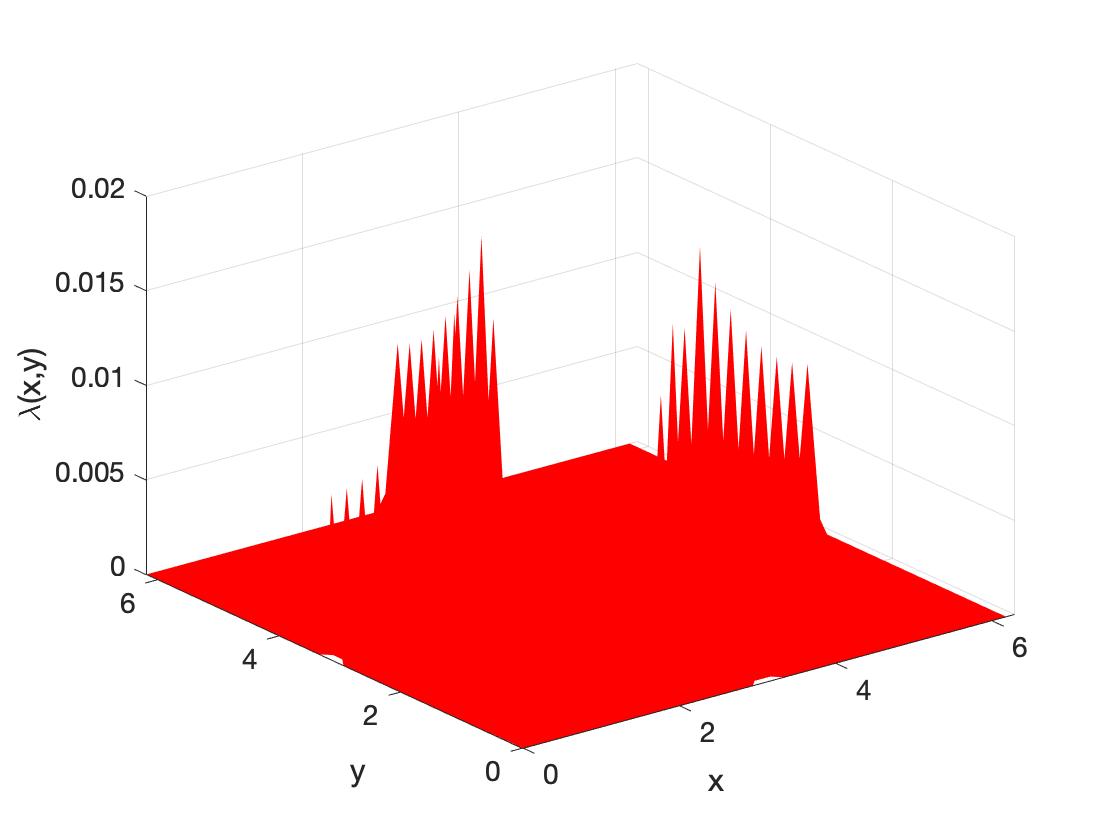}}
\subfigure[$\lambda$ at $t=0.2$.]{
\includegraphics[width=0.4\textwidth,clip==]{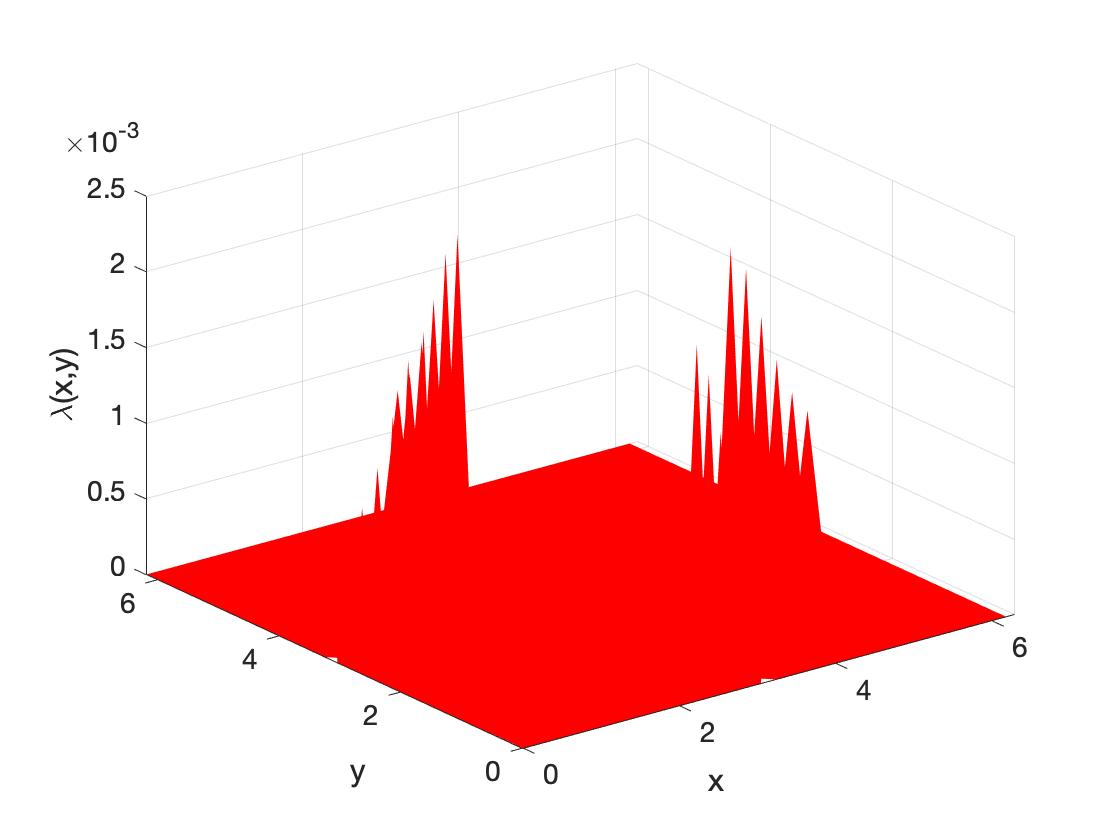}}
\caption{Numerical solutions $\rho(\bx,t)$ and Lagrange multiplier $\lambda(\bx,t)$ at $t=0,005,0.02$ computed by positivity preserving scheme \eqref{high:positivity:lag:1}-\eqref{high:positivity:lag:4} with $k=2$ and time step $\delta t=10^{-4}$.}\label{positivity}
\end{figure}

From Fig.\;\ref{convergence}, we observe that  first-order and second-order convergence rates of $\rho$ and $c$  are achieved for positivity schemes  \eqref{high:positivity:lag:1}-\eqref{high:positivity:lag:6} with $k=1,2$. To validate the positivity property of scheme \eqref{high:positivity:lag:1}-\eqref{high:positivity:lag:6}, we compute numerical results with larger time step $\delta t=10^{-4}$. In Fig.\;\ref{positivity},  we depict  numerical solutions $\rho(\bx,t)$  at $t=0.05$ and $t=0.2$,  respectively. It is observed that  all values of $\rho(\bx,t)$ at collocation points are positive. The corresponding Lagrange multiplier $\lambda(\bx,t)$ at $t=0.05$ and $t=0.2$ are also shown in Fig.\;\ref{positivity}.

To demonstrate the convergence rate of density $\rho(\bx,t)$ and concentration $c(\bx,t)$ for bound-preserving schemes \eqref{keller:lag:1}-\eqref{keller:lag:6}, we choose initial condition as
\begin{equation}
\begin{split}
&\rho(\bx,0)=10 e^{-10(x-\pi)^2-10(y-\pi)^2},\\
&c(\bx,0)=30 e^{-(x-\pi)^2/2-(y-\pi)^2/2}.
\end{split}
\end{equation} 
For computations, we set coefficients to be  $\eps=0.01$, $\mu=\gamma=\chi=1$ and $M=100$.  In Fig.\;\ref{convergence2}, the first-order and second-order convergence rate for density $\rho(\bx,t)$  and concentration $c(\bx,t)$ are shown by using  bound-preserving schemes \eqref{keller:lag:1}-\eqref{keller:lag:6} with $k=1,2$. Both the reference solutions for $\rho(\bx,t)$ and $c(\bx,t)$ are computed by second-order \eqref{keller:lag:1}-\eqref{keller:lag:6} with extremely small time step $\delta t=10^{-6}$.

\begin{figure}[htbp]
\centering
\subfigure[$c$ at $t=0.001$.]{
\includegraphics[width=0.4\textwidth,clip==]{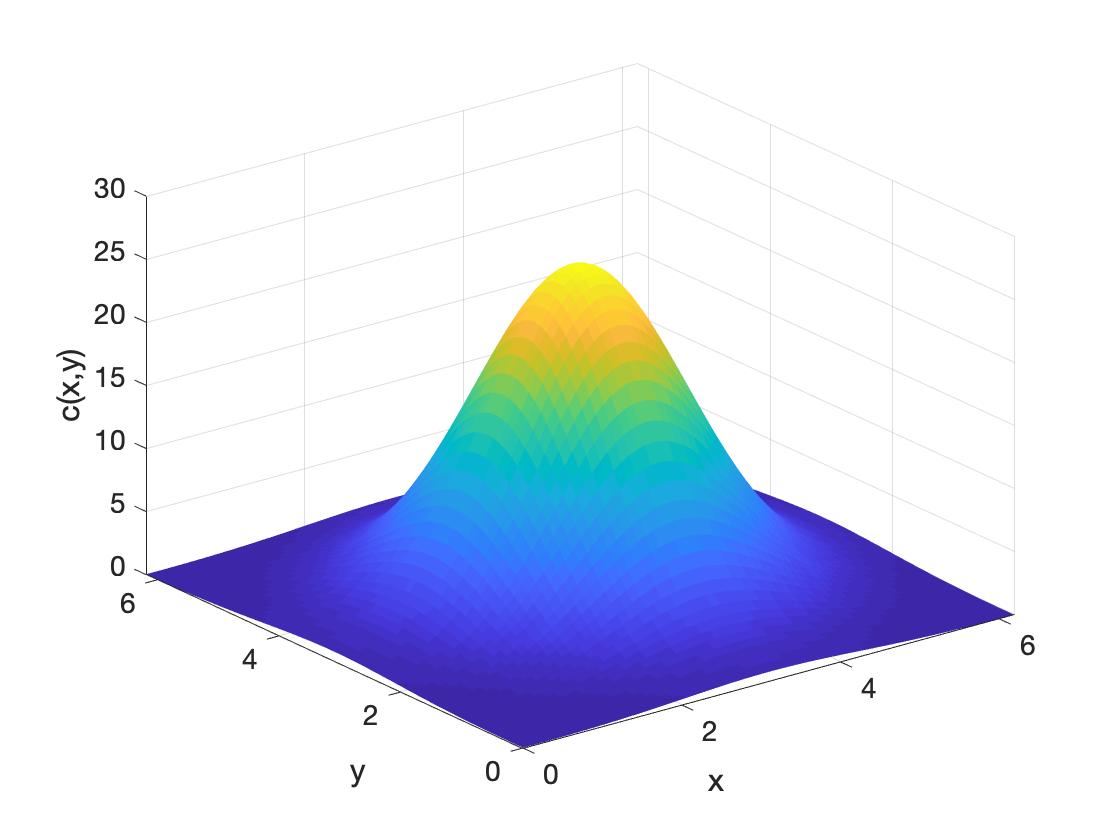}}
\subfigure[$c$ at $t=0.1$.]{
\includegraphics[width=0.4\textwidth,clip==]{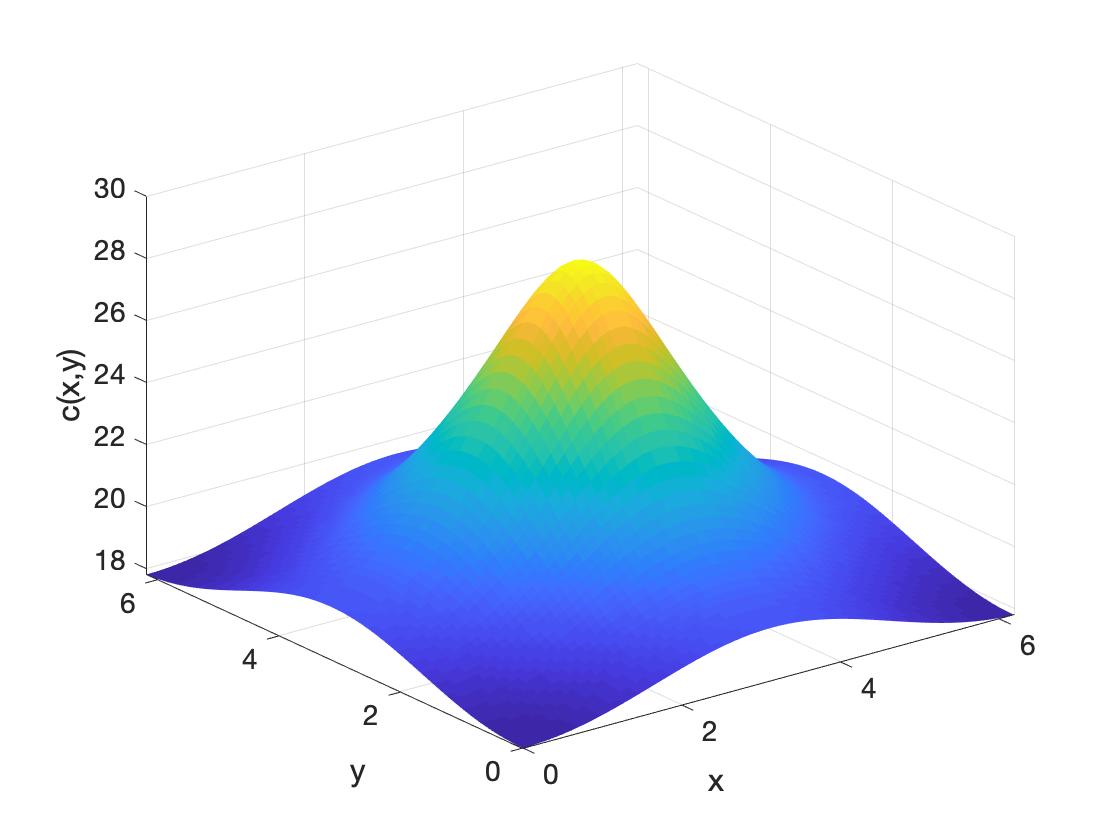}}
\subfigure[$\rho$ at $t=0.1$.]{
\includegraphics[width=0.4\textwidth,clip==]{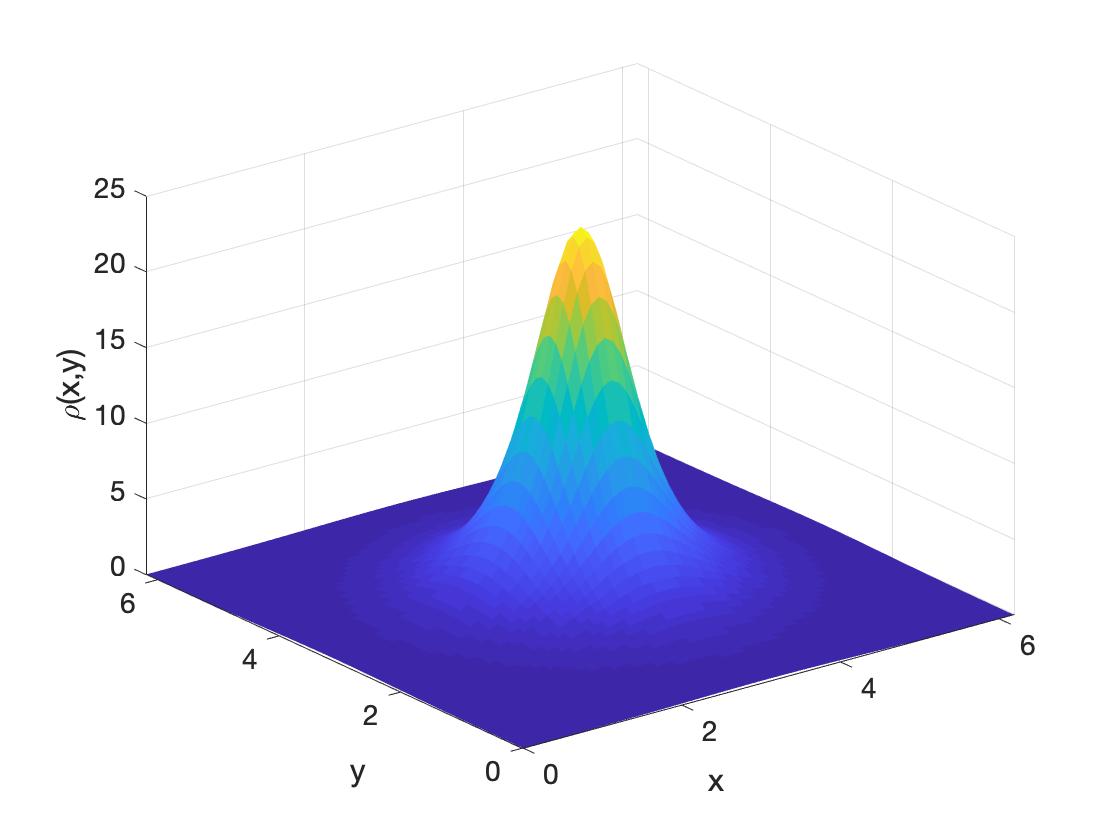}}
\subfigure[$\rho$ at $t=2$.]{
\includegraphics[width=0.4\textwidth,clip==]{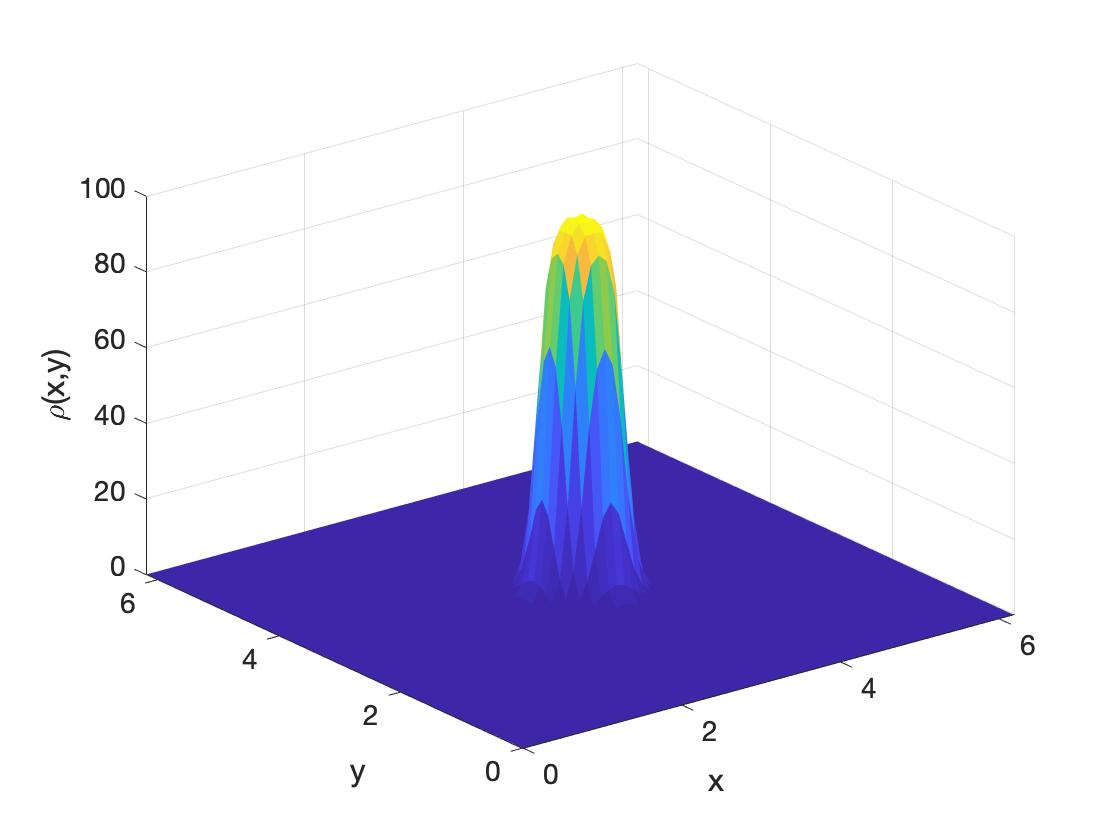}}
\subfigure[$\lambda$ at $t=2$]{
\includegraphics[width=0.4\textwidth,clip==]{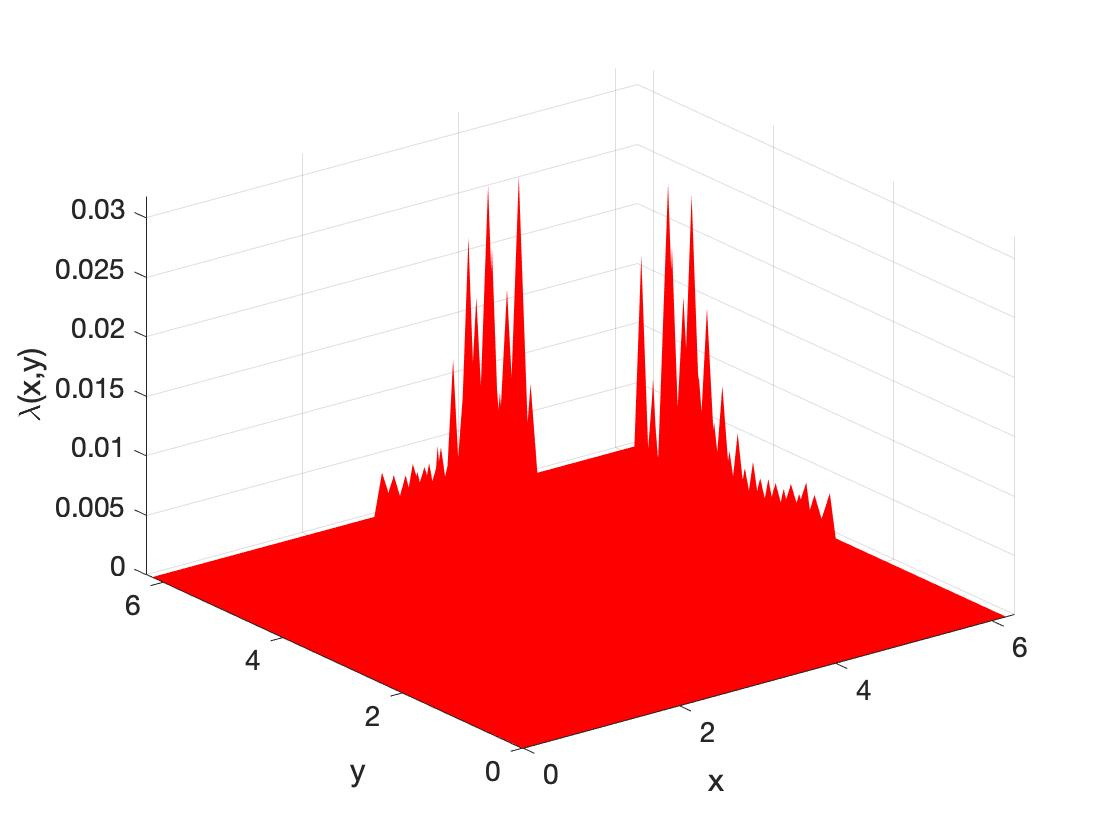}}
\subfigure[$\max{\rho}, \min{\rho}$]{
\includegraphics[width=0.4\textwidth,clip==]{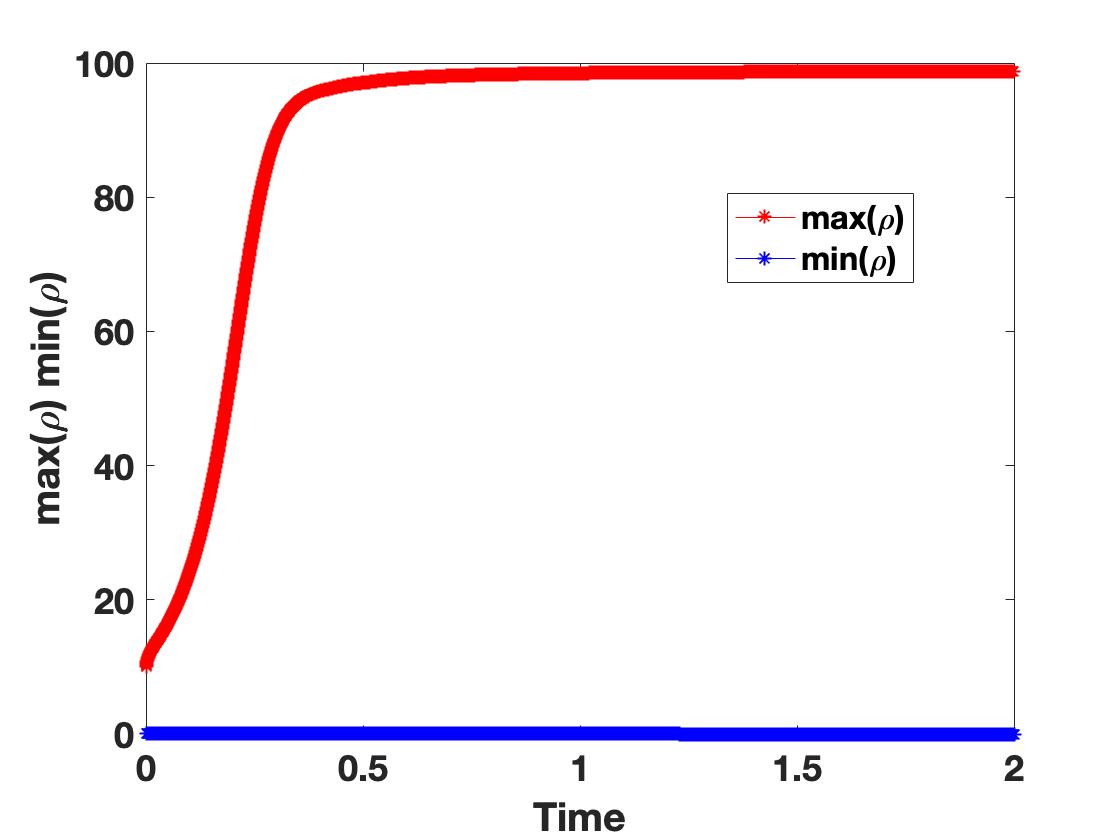}}
\caption{(a)-(d): Concentration $c$ and density $\rho$ computed by second-order bound-preserving scheme \eqref{high:positivity:lag:1}-\eqref{high:positivity:lag:4} with time step $\delta t=10^{-4}$. (e) Lagrange multiplier $\lambda$.  (f): Maximum and minimum value of $\rho$. }\label{bound:convergence}
\end{figure}

We apply scheme \eqref{high:positivity:lag:1}-\eqref{high:positivity:lag:6} with time step $\delta t=10^{-4}$ and $k=2$ to show positivity and bound-preserving properties. In Fig.\;\ref{bound:convergence}, we present concentration $c(\bx,t)$ at $t=0.001, 0.1$ and density $\rho$ at $t=0.1, 2$. We can observe that concentration  $c(\bx,t)$ is positive at all collocations   and density $\rho(\bx,t)$ is located in interval  $[0,100]$ which can be observed in Fig.\;\ref{bound:convergence}.(f).

\subsection{Blow up}
\begin{figure}[htbp]
\centering
\subfigure[$\rho$ at $t=0$.]{
\includegraphics[width=0.40\textwidth,clip==]{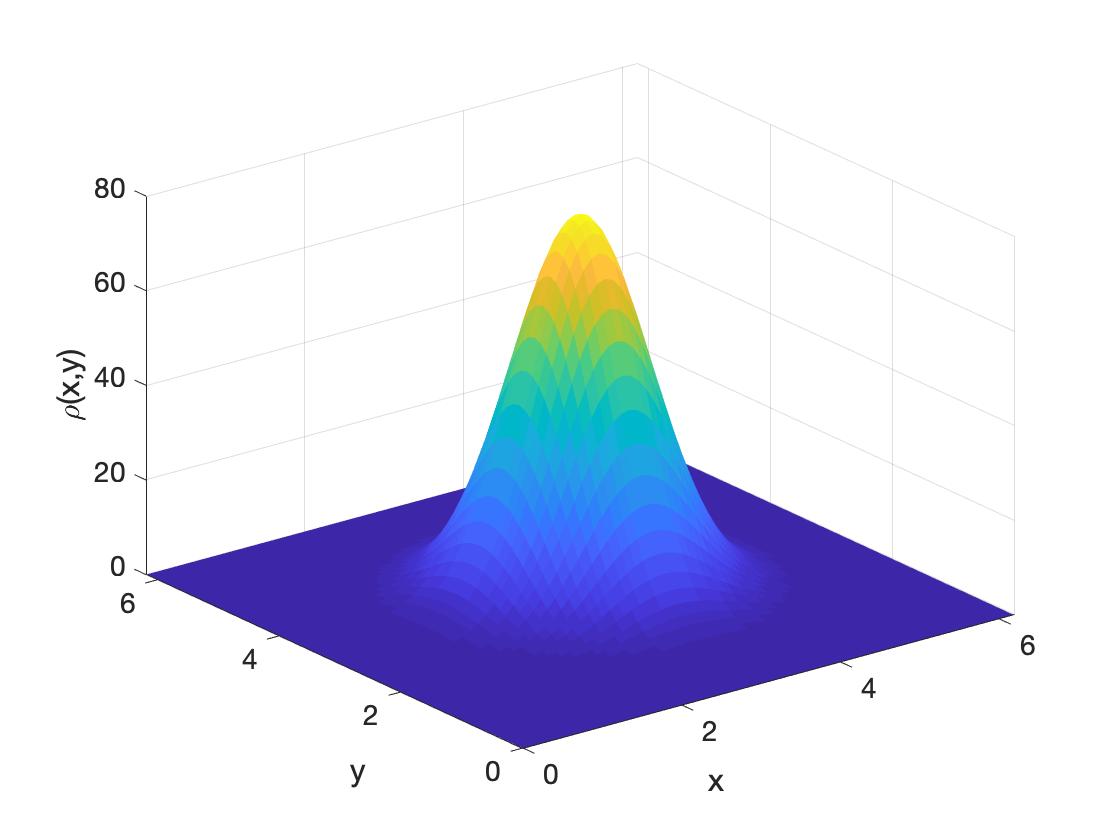}}
\subfigure[$\rho$ at $t=0.005$.]{
\includegraphics[width=0.40\textwidth,clip==]{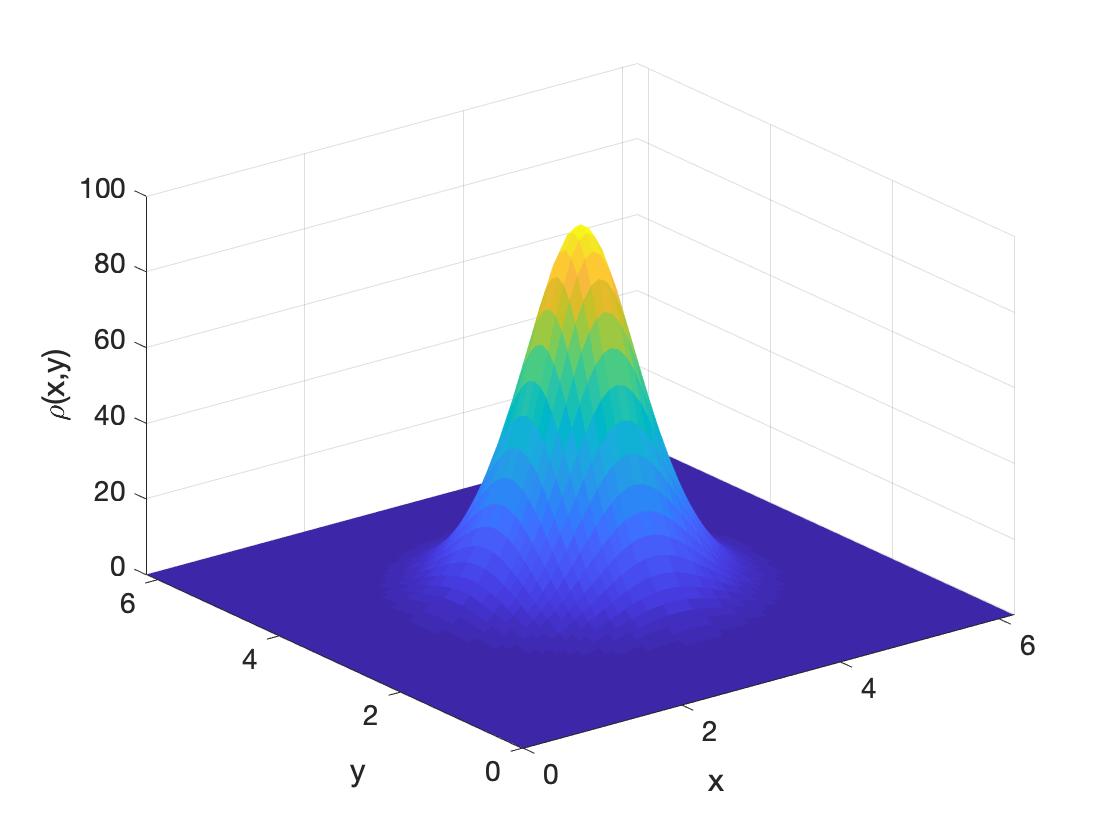}}
\subfigure[$\rho$ at $t=0.01$.]{
\includegraphics[width=0.40\textwidth,clip==]{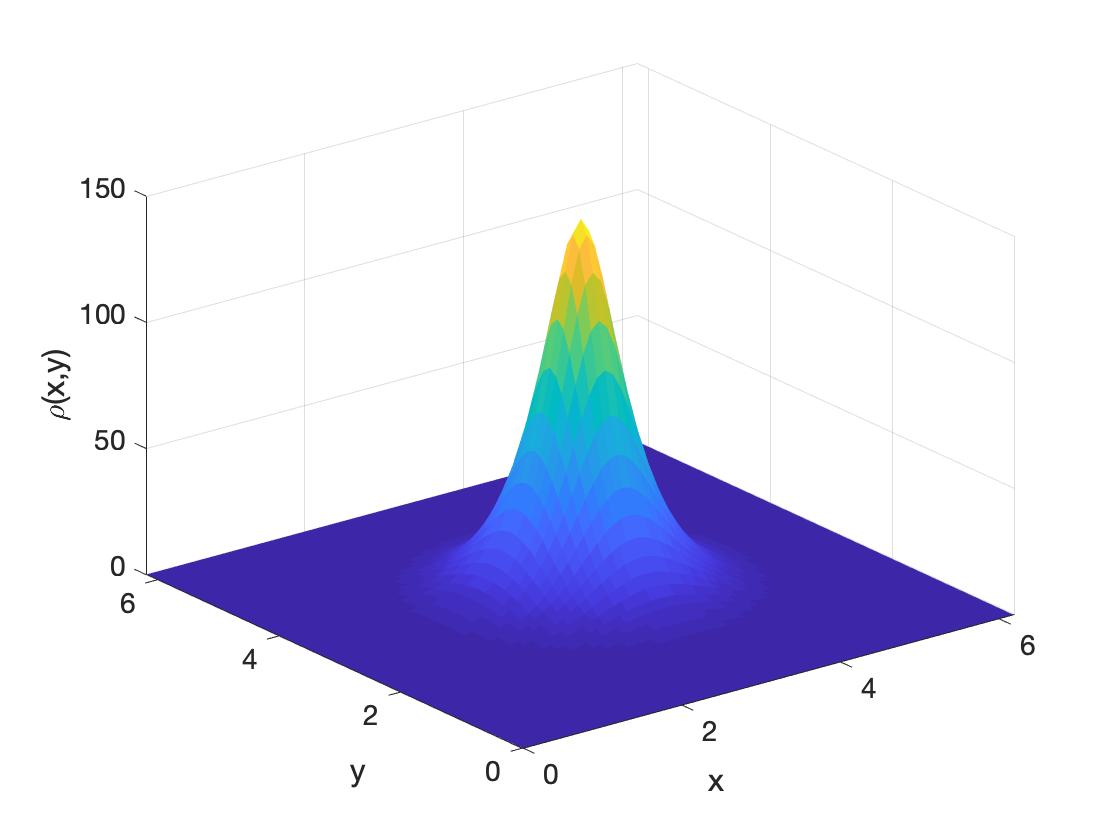}}
\subfigure[$\rho$ at $t=0.02$.]{
\includegraphics[width=0.40\textwidth,clip==]{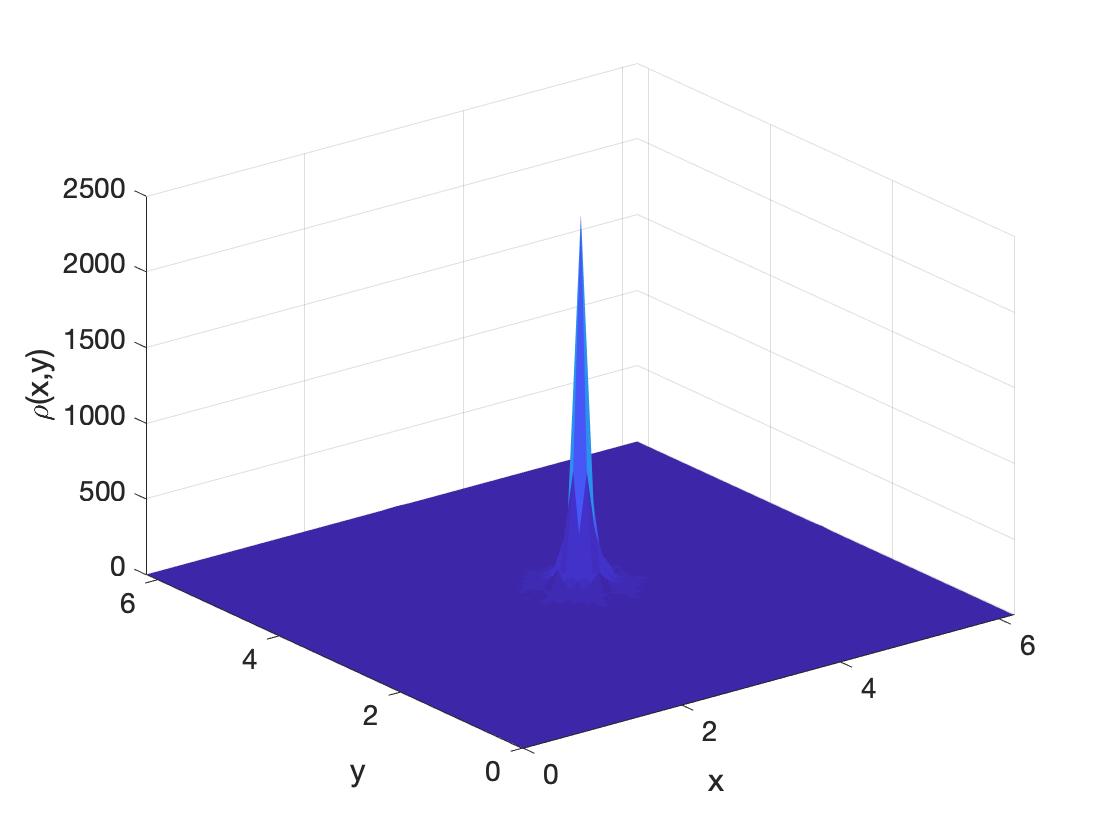}}
\caption{Numerical solution $\rho$ at $t=0,0.005,0.01,0.02$ computed by positivity preserving scheme \eqref{high:positivity:lag:1}-\eqref{high:positivity:lag:4} with $k=2$ and  time step $\delta t=10^{-4}$.}\label{blow_up}
\end{figure}

In this subsection, we focus on the case when $\rho$ blow up \cite{carrillo2012cross,calvez2006volume} for type-I Keller-Segel equation \eqref{de:keller:1}-\eqref{de:keller:3}  and show our positivity  schemes \eqref{high:positivity:lag:1}-\eqref{high:positivity:lag:6}  can simulate the  phenomenon of blow up \cite{zhou2017finite}.  The initial  value is given by 
\begin{equation}
\begin{split}
&\rho(\bx,0)=80 e^{-(x-\pi)^2-(y-\pi)^2},\\
&c(\bx,0)=30 e^{-(x-\pi)^2-(y-\pi)^2},
\end{split}
\end{equation}
in the domain $[0,2\pi]^2$.
In this example, we let $\eps=0.01$, $\mu=\gamma=\chi=1$ and plot profiles of $\rho$ at different time in Fig.\;\ref{blow_up} which display an intensity of blowing up. Since we can observe a sharp profile of $\rho$ at $t=0.02$ appears. Values of $\rho$ at any collocation points  are all positive for different time depicted in Fig.\;\ref{blow_up} which are consistent with KKT-conditions for positivity schemes.

\subsection{Comparison with a generic scheme}
In this subsection, we mainly compare our bound-preserving schemes \eqref{high:positivity:lag:1}-\eqref{high:positivity:lag:6} with a generic semi-implicit scheme \eqref{semi-implicit}  which takes on the form 
\begin{equation}\label{semi-implicit}
\begin{split}
&\eps\frac{3 c_h^{n+1}-4c_h^n+c_h^{n-1}}{2\delta t} = \mu \Delta c_h^{n+1} + 2\rho_h^n-\rho_h^{n-1},\\
&\frac{3\rho_h^{n+1}-4\rho_h^n+\rho_h^{n-1}}{2\delta t}
=\gamma \Delta \rho_h^{n+1}-\chi\Grad\cdot(\eta(2\rho_h^{n}-\rho_h^{n-1})\Grad c^{n+1}_h),
\end{split}
\end{equation}
in the domain $[0,2\pi]^2$. 

\begin{figure}[htbp]
\centering
\subfigure[$\rho$ at $t=1$.]{
\includegraphics[width=0.4\textwidth,clip==]{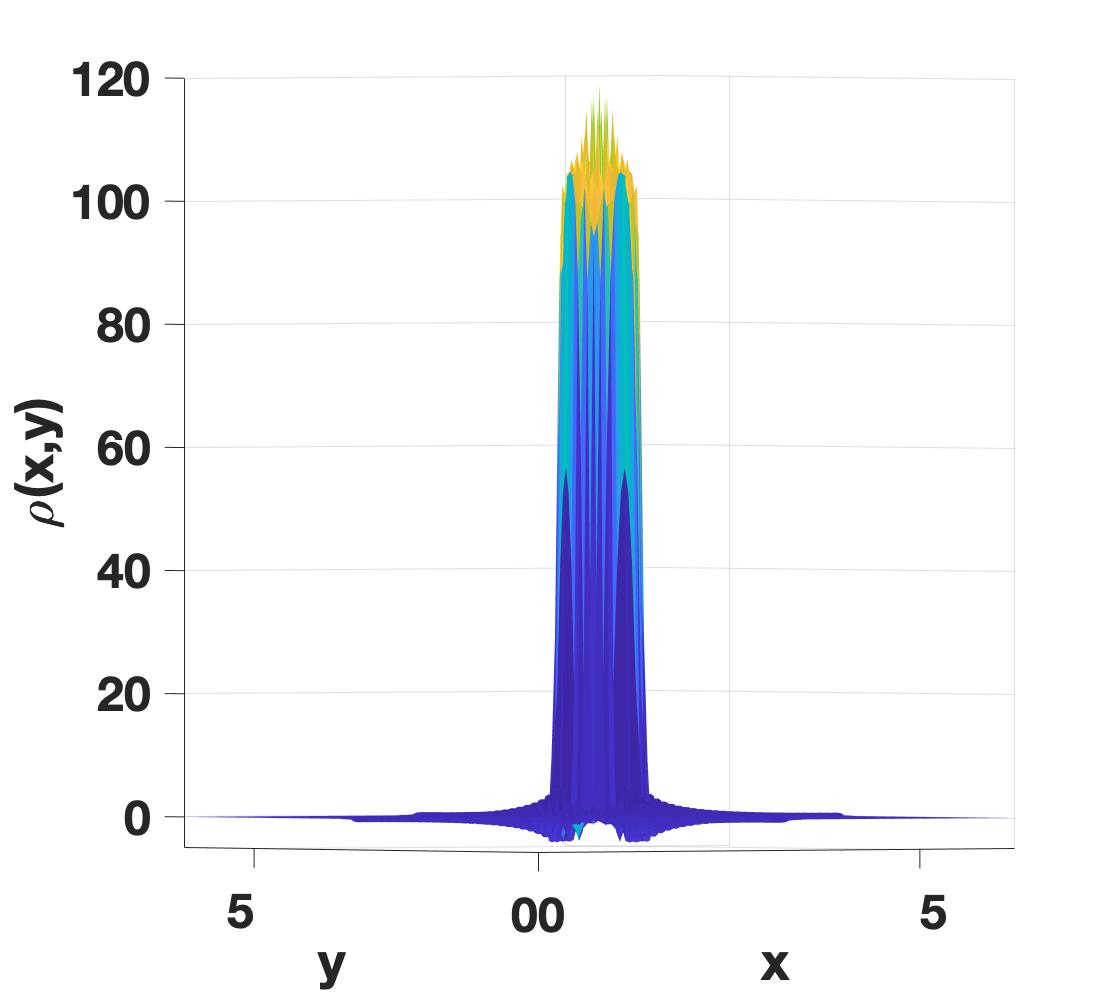}}
\subfigure[Without Lagrange multiplier.]{
\includegraphics[width=0.4\textwidth,clip==]{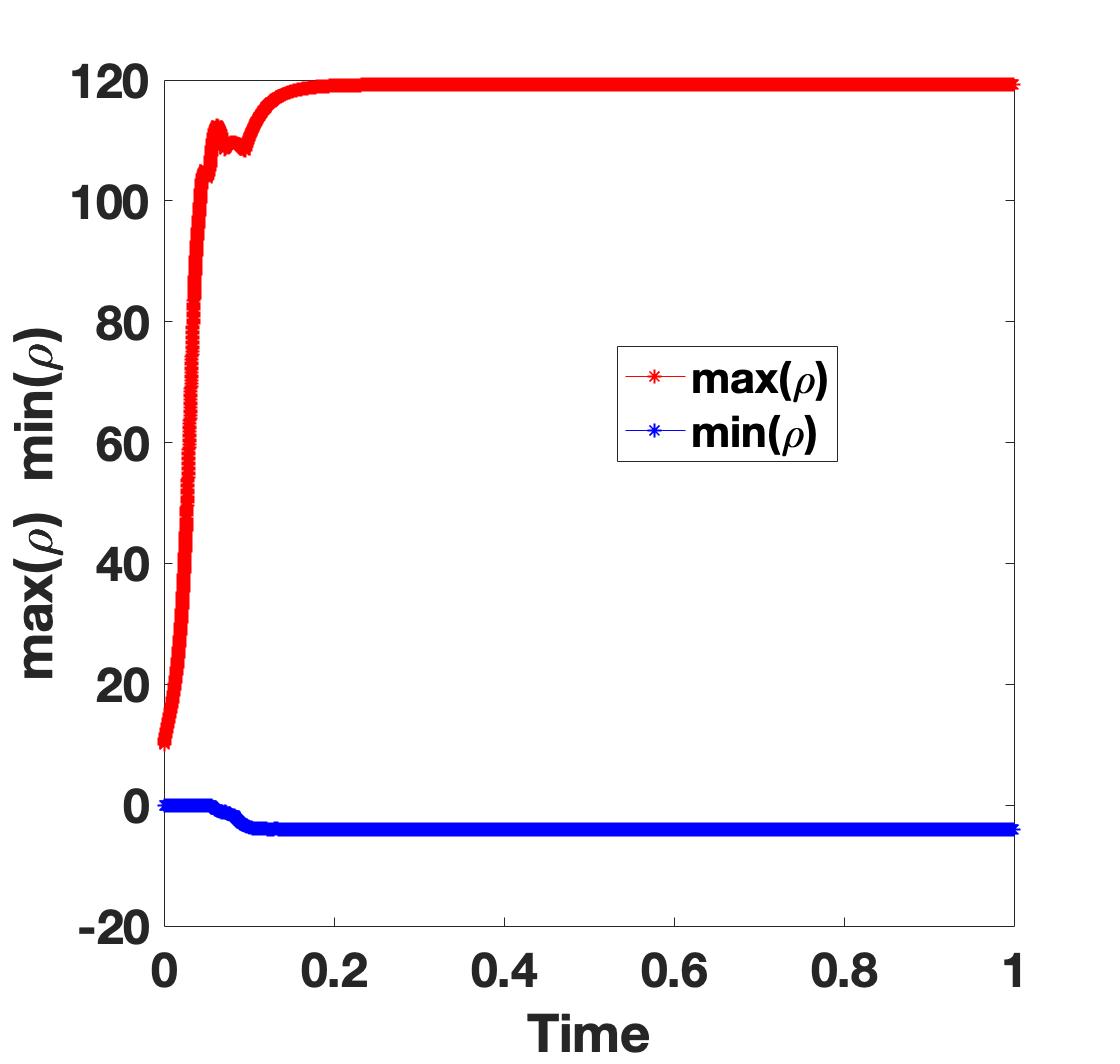}}
\caption{(a): Numerical solution $\rho$ at $t=1$ computed by usual semi-implicit scheme \eqref{semi-implicit} with time step $\delta t=10^{-4}$. (b): Evolutions of $\max(\rho)$ and $\min(\rho)$ with respective to time using semi-implicit scheme.}\label{generic_bound}
\end{figure}

For comparison, we choose the same time step $\delta t=10^{-4}$ and same coefficients of Keller-Segel \eqref{keller:1}-\eqref{keller:3} $\gamma=\mu=1$, $\chi=2$, $\eps=0.01$ and $M=100$.  $128^2$ Fourier points are implemented in space. The initial condition is set  to be  
\begin{equation}
\begin{split}
& \rho(\bx,0)=10e^{-(x-\pi)^2/2-(y-\pi)^2/2},\\
& c(\bx,0)= 30e^{-(x-\pi)^2/2-(y-\pi)^2/2}.
\end{split}
\end{equation}
In Fig.\;\ref{generic_bound}, we depict $\rho$ at $t=1$ and evolutions of maximum and minimum values of $\rho$ with respective to time.  We observe that values of $\rho$ will be negative and larger than $M$ at some collocation points by using generic semi-implicit scheme \eqref{semi-implicit}. The range of $\rho$ will not located inside interval $[0,M]$ from the evolutions. 

\begin{figure}[htbp]
\centering
\subfigure[$\rho$ at $t=0$.]{
\includegraphics[width=0.4\textwidth,clip==]{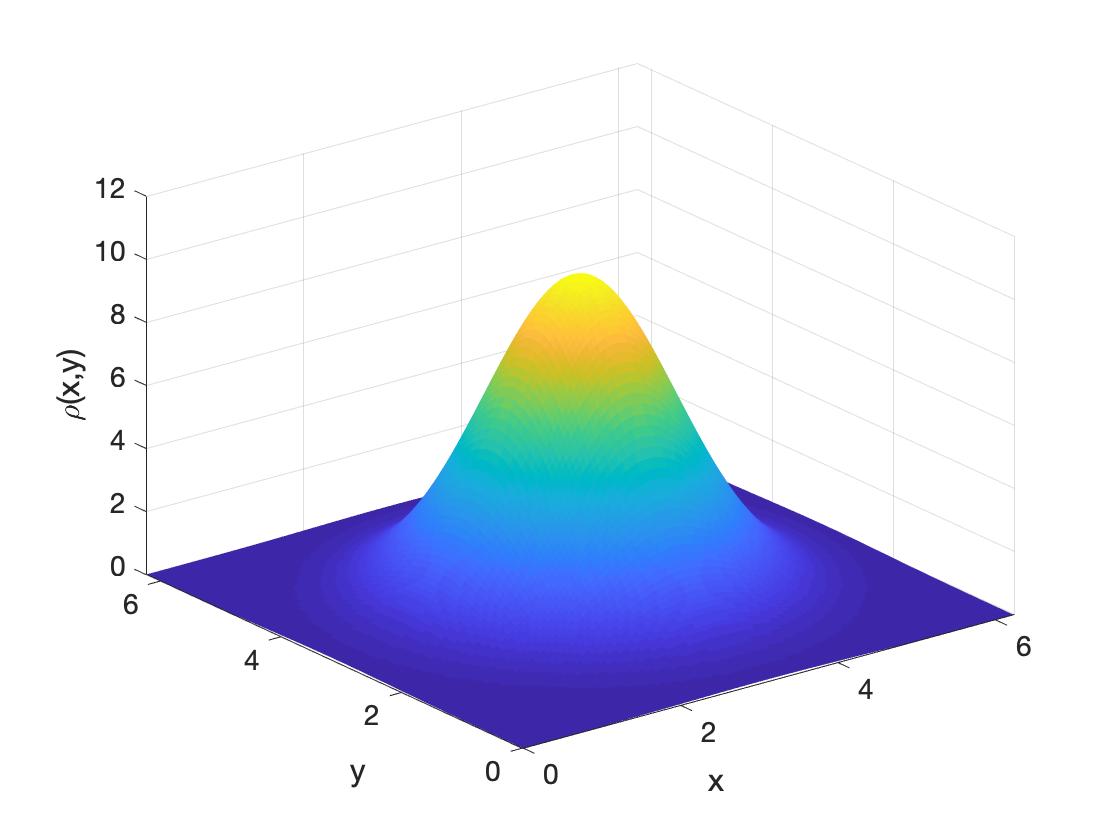}}
\subfigure[$\rho$ at $t=0.01$.]{
\includegraphics[width=0.4\textwidth,clip==]{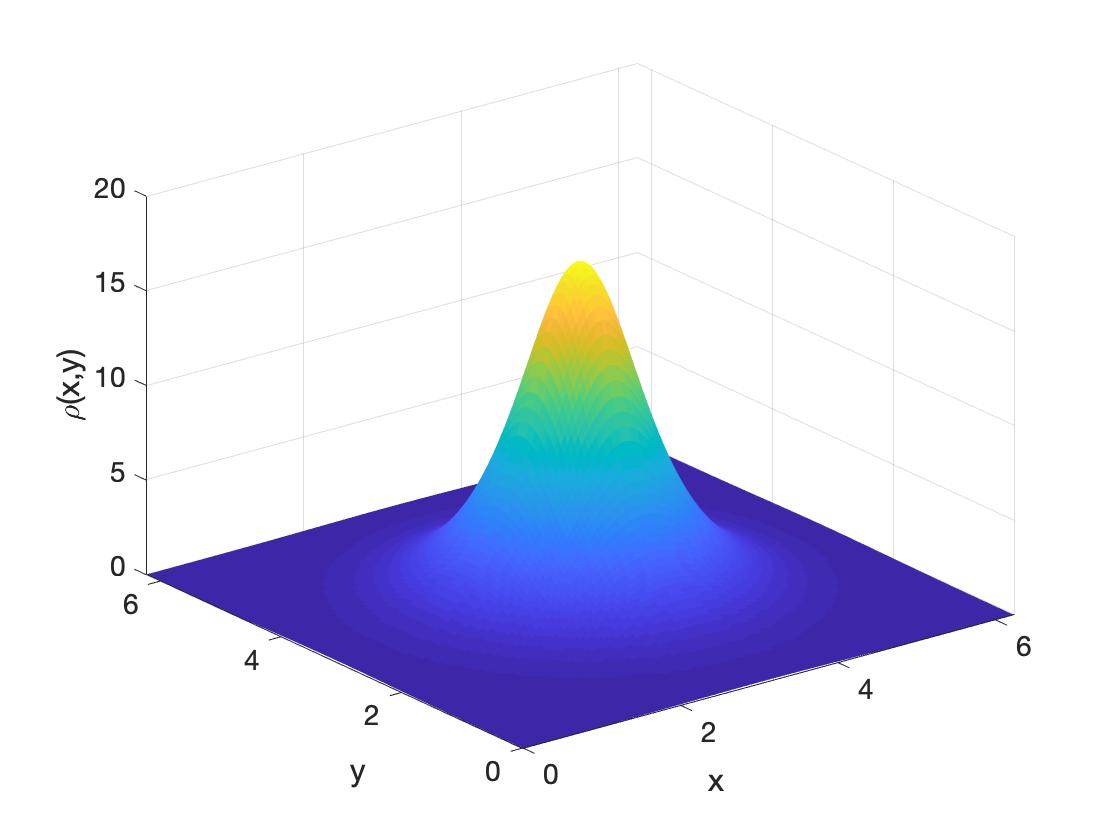}}
\subfigure[$\rho$ at $t=0.05$.]{
\includegraphics[width=0.4\textwidth,clip==]{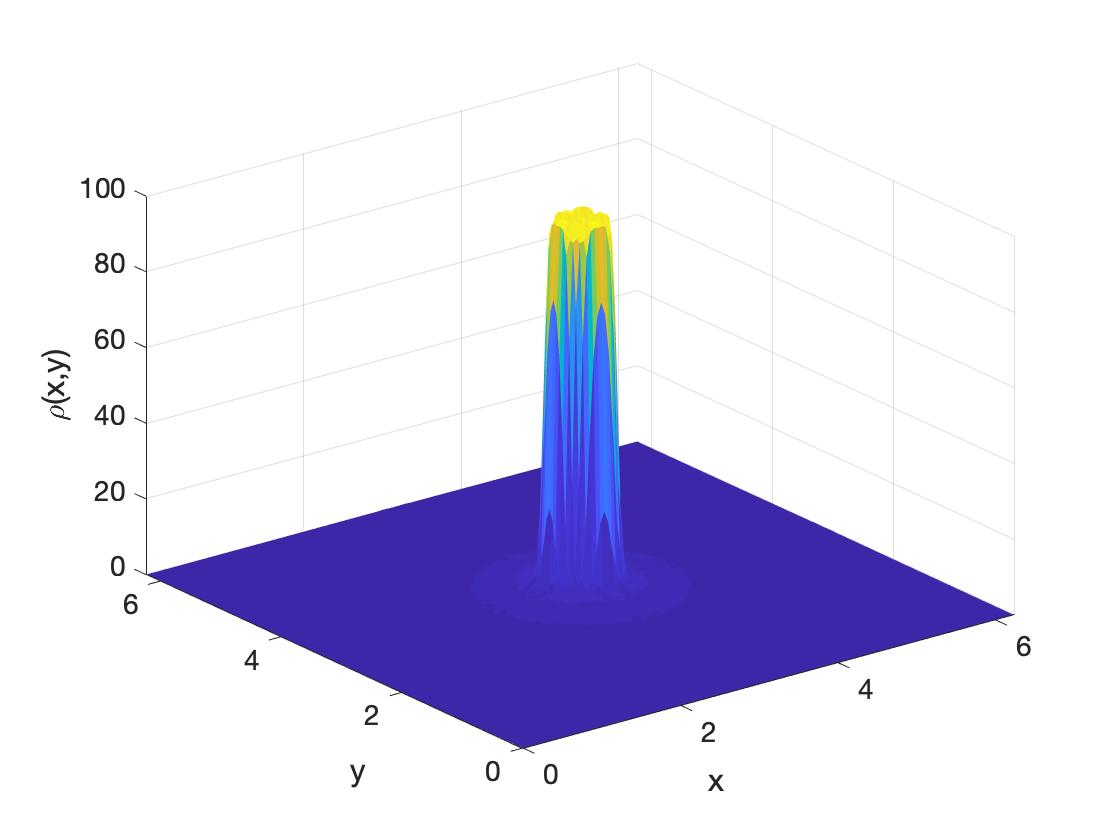}}
\subfigure[$\rho$ at $t=1$.]{
\includegraphics[width=0.4\textwidth,clip==]{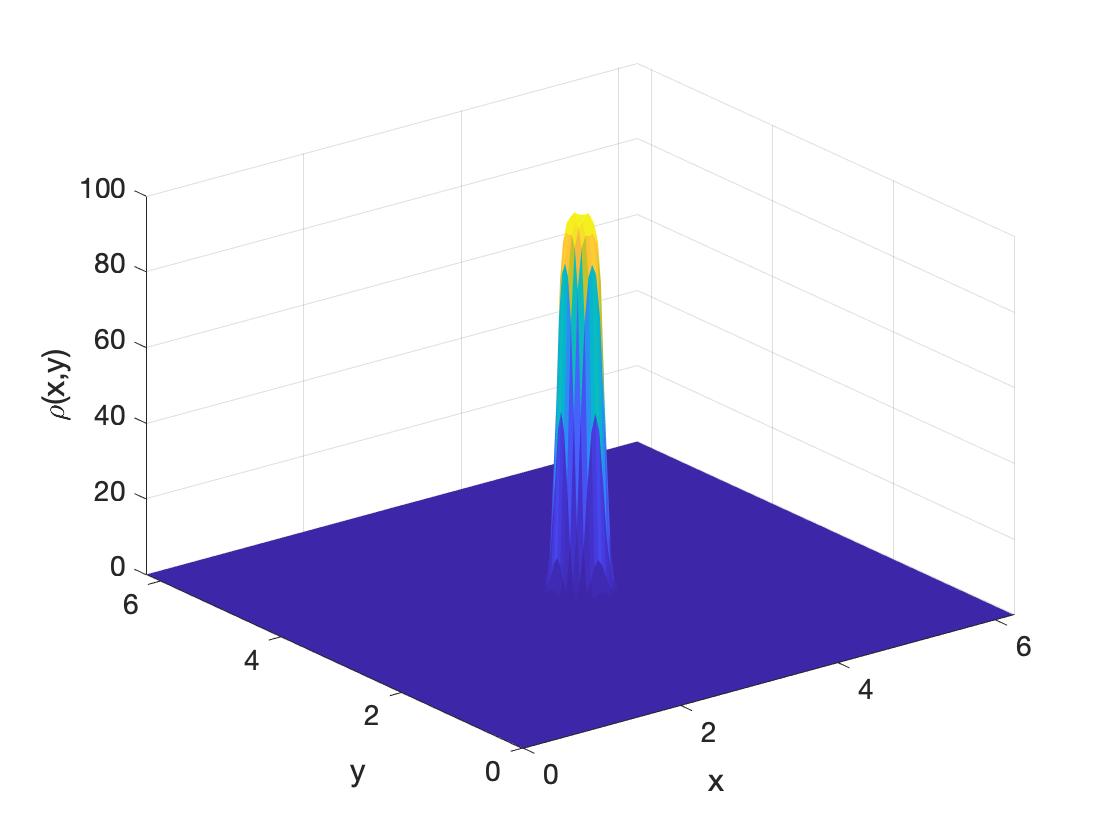}}
\subfigure[$\lambda$ at $t=0.005$.]{
\includegraphics[width=0.4\textwidth,clip==]{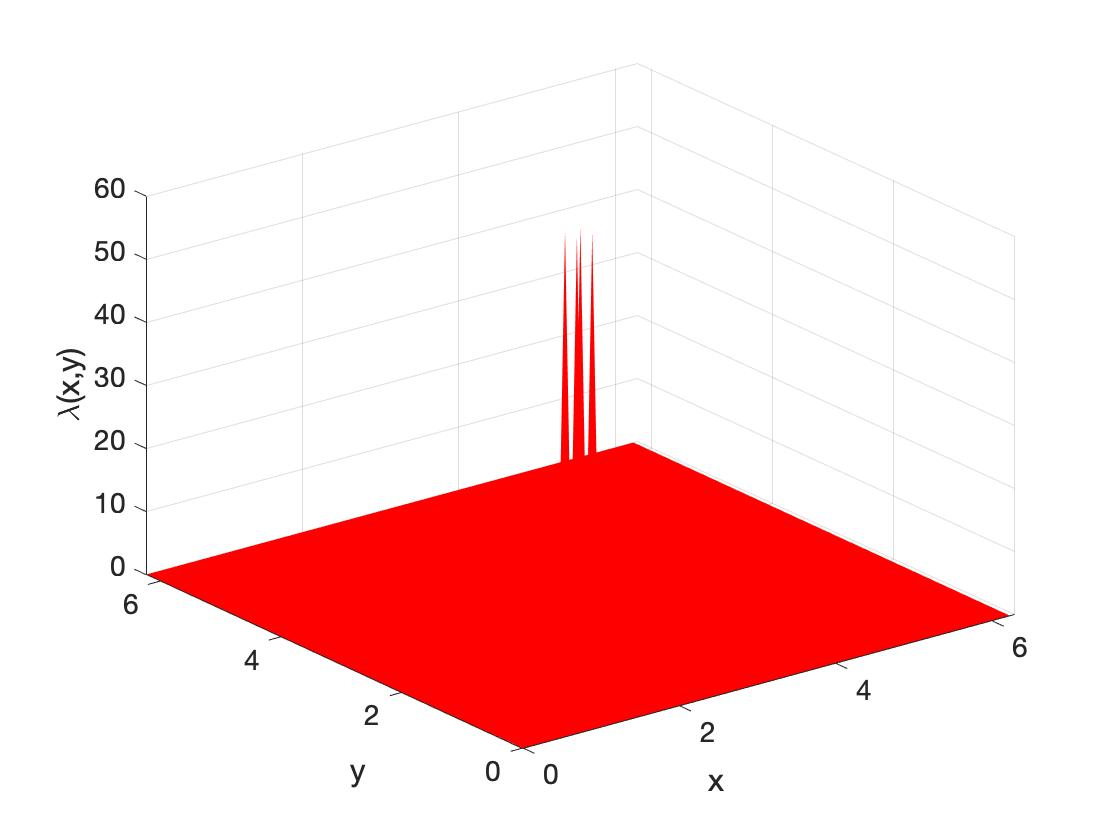}}
\subfigure[$\max(\rho$ and  $\min(\rho)$.]{
\includegraphics[width=0.40\textwidth,clip==]{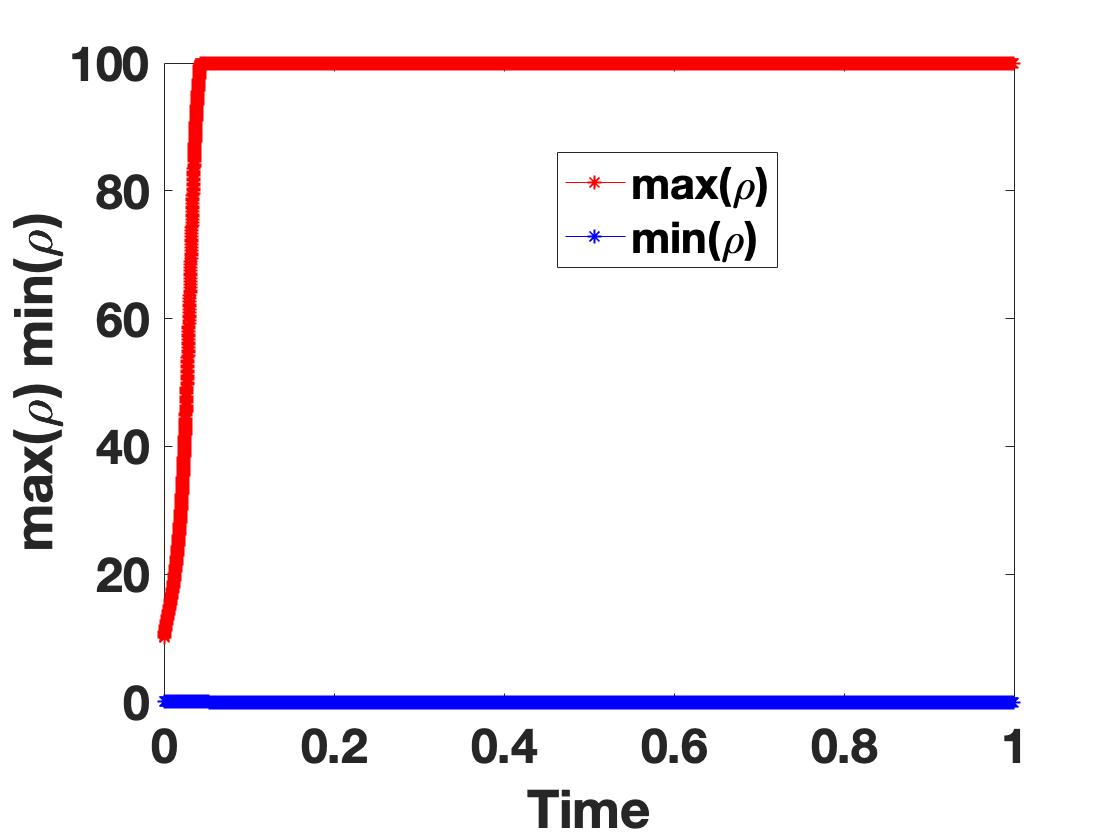}}
\caption{(a)-(d): Numerical solution $\rho$ at $t=0,0.01,0.05,1$ computed by bound preserving scheme \eqref{keller:lag:1}-\eqref{keller:lag:4} with $k=2$ and  time step $\delta t=10^{-4}$. (e): Lagrange multiplier $\lambda$ at $t=0.005$. (f): Evolutions of $\max(\rho)$ and $\min(\rho)$ using bound-preserving scheme.}\label{bound}
\end{figure}

 Observed from Fig.\;\ref{bound}, we find the maximum  and minimum value of $\rho$ always stay within interval $[0,100]$ by using bound-preserving scheme \eqref{keller:lag:1}-\eqref{keller:lag:6} with $k=2$. The profiles of $\rho$ at various time in Fig.\;\ref{bound} are consistent with the property of bound-preserving. Lagrange multiplier $\lambda$  is non-zero and positive at $t=0.05$ which is known from KKT-conditions.

\section{Concluding remarks}
We propose in this paper a  improved  Lagrange multiplier approach and a new energy correction approach  to construct high-order, energy dissipative, positivity/bound-preserving schemes for two  types  of Keller-Segel equations. More precisely, we introduce  space-time Lagrange multiplier functions to enforce the positivity and bound-preserving constraints, and express the expanded system using the KKT conditions; then using a generic  IMEX  scheme to compute a predictive numerical solution which does not satisfy positivity or bound-preserving constraints,  followed by a correction step to produce a  numerical solution with  positivty or bound-preserving constraints. A energy correction approach is introduced in the final step to preserve a discrete energy dissipative law. 

An important feathure is that the new schemes proposed in this paper can be performed with same computational cost compared with the generic  IMEX  scheme.  Both positivity-preserving and bound-preserving  schemes we proposed can be proved to be unique solvable in which Lagrange multipliers can be uniquely determined. 

We established  stability results of the new class of structure-preserving  schemes  for both types of Keller-Segel equations. A error analysis is also presented for a second-order, energy dissipative, bound-preserving scheme for the second type of Keller-Segel equations. Ample numerical experiments are implemented  to validate the accuracy, stability  and properties of positivity and bound-preserving.  We also compare with a generic scheme to show the advantage of our  new proposed schemes for Keller-Segel equations.

\section{Acknowledgment}
The author would like to thank Professor Chun Liu and Professor Jie Shen  for enlightening discussions and their encouragement.

\bibliographystyle{siamplain}
 \bibliographystyle{plain}
\bibliography{references}
\end{document}